\documentclass[3p]{elsarticle}

\usepackage{amsmath}
\usepackage{amsthm}
\usepackage{amssymb}
\usepackage{bm}
\usepackage{accents}

\newtheorem{thm}{Theorem}
\newtheorem{lem}{Lemma}
\newdefinition{rmk}{Remark}
\newproof{pf}{Proof}

\numberwithin{equation}{section}

\graphicspath{{figs/}}

\journal{arXiv}

\begin{document}

\begin{frontmatter}

\title{Decoupling technology for systems of evolutionary equations\tnoteref{label1}}
\tnotetext[label1]{The work was supported by the Russian Science Foundation (grant No. 24-11-00058).}

\author{P.N. Vabishchevich\corref{cor1}\fnref{lab1,lab2}}
\ead{vab@cs.msu.ru}
\cortext[cor1]{Correspondibg author.}

\address[lab1]{Lomonosov Moscow State University, 1, building 52, Leninskie Gory,  119991 Moscow, Russia}

\address[lab2]{North-Eastern Federal University, 58, Belinskogo st, Yakutsk, 677000, Russia}

\begin{abstract}

Numerical methods of approximate solution of the Cauchy problem for coupled systems of evolution equations are considered.
Separating simpler subproblems for individual components of the solution achieves simplification of the problem at a new level in time.
The decoupling method, a significant approach to simplifying the problem, is based on the decomposition of the problem's operator matrix.
The approximate solution is constructed based on the linear composition of solutions to auxiliary problems.
The paper investigates decoupling variants based on extracting the diagonal part of the operator matrix and the lower and upper triangular submatrices.
The study introduces a new decomposition approach, which involves splitting the operator matrix into rows and columns.
The composition stage utilizes various variants of splitting schemes, showcasing the versatility of the approach.
In additive operator-difference schemes, we can distinguish explicit-implicit schemes, factorized schemes for two-component splitting, and regularized schemes for general multicomponent splitting.
The study of stability of two- and three-level decoupling composition schemes is carried out using the theory of stability (correctness) of operator-difference schemes for finite-dimensional Hilbert spaces.
The theoretical results of the decoupling technique for systems of evolution equations are illustrated on a test two-dimensional problem for a coupled system of two diffusion equations with inhomogeneous self- and cross-diffusion coefficients.
\end{abstract}

\end{frontmatter}

\section{Introduction}\label{s-1}

Applied multi-physics models are usually based on systems of coupled nonstationary equations.
They include, for example, parabolic and hyperbolic equations for the desired scalar and vector quantities.
Numerical methods for approximate solutions are constructed using finite-element and finite-volume approximations over space \cite{KnabnerAngermann2003,QuarteroniValli}.
The computational complexity of such problems stems from the fact that we have a coupled system of equations to find a solution at a new level in time.
We have well-developed computational solvers for basic partial derivative equations for particular mono-physical problems.
Decoupling technologies allow us to simplify the problem by solving a sequence of more straightforward problems for separate scalar and vector components of the desired solution.

Applying explicit time approximations provides the most straightforward decomposition of a complex problem into simpler problems.
In this case, we can find an approximate solution relatively quickly at a new time level.
The fundamental disadvantage of this approach is related to the severe constraints on the time step of the grid \cite{Samarskii1989,SamarskiiMatusVabischevich2002}.
We can only count on the conditional convergence of the approximate solution to the exact one.
In computational practice, implicit time approximations \cite{Ascher2008,LeVeque2007} are often used in numerical solutions for nonstationary boundary-value problems for partial derivative equations.
We have a more computationally complex problem on a new level of time. Still, the solution is unconditionally stable concerning the initial data and the right-hand side.
The decoupling method for nonstationary problems consists of decomposing complex problems into simpler ones by choosing particular time approximations.

After the space approximation of the initial boundary value problems for nonstationary partial derivative equations systems, we arrive at the Cauchy problem for systems of evolution equations for scalar and vector quantities.
The system of linear equations of the first order can serve as a basic one. Its singularities are related to the operator matrix for the solution's individual components.
Decoupling technology consists of two stages \cite{DC}.
In the decomposition (analysis) stage, we use an additive representation of this operator matrix, separating individual summands associated with more straightforward problems.
In the composition (analysis) stage, auxiliary simple problems are extracted for the individual summands of the additive representation of the problem operator by selecting approximations.
The approximate solution at a new level in time is sought as a superposition of solutions to auxiliary problems.

Implicit time approximations are used in the approximate solution of auxiliary problems, which yield unconditionally stable operator-difference schemes.
The additive decomposition of the problem operator(s) should provide an acceptable computational challenge at a new level in time.
In the approximate solution of systems of evolutionary equations, the problem's operator matrix can be decoupled in various ways.
Decoupling of the approximate solution is most easily constructed for the diagonal part of the operator matrix.
A second option for computationally convenient additive decoupling is related to the extraction of lower and (or) upper triangular matrices.
The critical points of decoupling technology for systems of parabolic equations on such a basis were considered by us in \cite{DecSys}.
Applying a general approach to constructing decomposition-composition schemes based on the additive representation of the unit operator when solving systems of evolution equations \cite{DC} is related to splitting the operator matrix into separate rows and columns.

At the two-component splitting of the operator matrix of the problem, we can use the simplest inhomogeneous time approximations, which correspond to the use of IMEX methods \cite{Ascher1995,HundsdorferVerwer2003}.
In this case, the operator summand providing decoupling of the approximate solution is taken from the upper level in time, and the other summand --- from the lower level.
In the approximate solution of problems for systems of evolutionary equations, we can use explicit-implicit two- and three-level operator-difference schemes for the first-order evolutionary equation, which are constructed in \cite{Vabishchevich2020}.

In the general case, the composition of the approximate solution is provided by applying additive operator-difference schemes (splitting schemes) \cite{Marchuk1990,VabishchevichAdditive}.
Their study is based on the theory of stability (correctness) of operator-difference schemes \cite{Samarskii1989,SamarskiiMatusVabischevich2002},
which gives necessary and sufficient conditions for the stability of two- and three-level schemes for finite-dimensional Hilbert spaces.
Let us note the main classes of operator-difference schemes that can be used after decomposing the operator matrix of the system of evolution equations \cite{VabishchevichAdditive}.
We can use factorized two-level schemes for two-component splitting, which are operator analogs of classical alternating direction schemes.
We use component-wise splitting schemes (sum approximation schemes) for general multicomponent splitting.
We can distinguish additive-averaging schemes among these schemes, which allow independent solutions of auxiliary problems when transitioning to a new level in time.
Regularized additive operator-difference schemes belong to the class of full approximation schemes. These splitting schemes are constructed based on multiplicative perturbation of separate decomposition operators.

The present paper reflects the research directions mentioned for decoupling methods, using the approximate solution of the Cauchy problem for linear systems of evolution equations for finite-dimensional Hilbert spaces.
The operator matrix is decomposed by extracting the diagonal part, triangular submatrices, and rows and columns.
Unconditional stable two- and three-level additive operator-difference is constructed by decoupling the approximate solution at a new level in time.

The paper is organized as follows.
In Section \ref{s-2}, we formulate the problem of constructing decoupling schemes for approximating nonstationary problems for systems of evolution equations.
We consider the Cauchy problem on a direct sum of finite-dimensional Hilbert spaces for a system of first-order equations with a self-adjoint operator.
Explicit-implicit approximations for extracting the diagonal and triangular parts of the operator matrix of the problem are considered in the \ref{s-3} section.
Two- and three-level composition schemes are constructed and investigated for stability.
In Section \ref{s-4}, we consider multicomponent decoupling schemes for splitting the operator matrix of the problem into rows and columns.
Stability conditions for component-wise splitting and regularized additive operator-difference schemes are formulated.
Section \ref{s-5} is devoted to numerical experiments on the approximate solution of a model two-dimensional problem.
A system of two self- and cross-diffusion equations with discontinuous coefficients is considered.
Numerical solutions are compared using decoupling schemes with standard two-level schemes of first and second-order accuracy. The results are summarized in Section \ref{s-6}.

\section{Problem statement}\label{s-2}

Let $H_{\alpha}, \ \alpha =1,2,\dots,p$ be finite-dimensional real Hilbert (Euclidean) spaces in which the scalar product and norm are $(\cdot,\cdot)_{\alpha}$  and $\|\cdot\|_{\alpha}, \ \alpha =1,2,\dots,p$ respectively.
Individual components of the solution will be denoted by $u_{\alpha}(t), \ \alpha =1,2,\dots,p$ for each $t$ ($0 \leq t \leq T, \ T > 0$). 
A solution is sought for a system of first-order evolution equations:
\begin{equation}\label{2.1}
  \frac{d u_{\alpha}}{d t} +
  \sum_{\beta =1}^{p} A_{\alpha \beta}  u_{\beta} = f_{\alpha}(t),
  \quad \alpha = 1, 2, \dots, p,
  \quad 0 < t \leq T .
\end{equation}
Here $f_{\alpha}(t), \ \alpha =1,2,\dots,p$  are given,
and $A_{\alpha \beta}$ are linear constants (independent of $t$) operators from $H_{\beta}$ to $H_{\alpha}$ 
($A_{\alpha \beta}: H_{\beta} \to H_{\alpha}$)
for all $\alpha =1,2,\dots,p$.
The system of equations (\ref{2.1}) is supplemented with the initial conditions
\begin{equation}\label{2.2}
  u_{\alpha}(0) = u_{\alpha}^0, 
  \quad \alpha =1,2,\dots,p .
\end{equation}
For simplicity, we restrict ourselves to the case when $f_{\alpha}(t) = 0, \ \alpha =1,2,\dots,p$.
The incorporation of the right-hand sides is not essential for achieving the main goal of our study, but it involves many technical details that clutter the text of the paper.

We arrive at systems of evolutionary equations of type (\ref{2.1}) when approximating several applied problems in space.
As an example, we note a system of coupled parabolic equations of the second order \cite{okubo2001diffusion,murray2003mathematical}, which is typical for diffusion problems in multicomponent media, when a solution to the system is sought in the domain $\Omega$
\[
\frac{\partial w_{\alpha}}{\partial t} -
\sum_{\beta =1}^{p} \operatorname{div} (d_{\alpha \beta} ({\bm x}) \operatorname{grad} w_{\beta})  = f_{\alpha}({\bm x},t),
\quad  \alpha = 1, 2, \dots, p ,
\quad  {\bm x} \in \Omega,
\quad  0 < t \leq T .
\]
The coefficients $d_{\alpha \beta}$ describe diffusion processes: self-diffusion --- for $\alpha = \beta$
and cross diffusion --- for $\alpha \neq \beta$.

We will interpret the system of equations (\ref{2.1}) as one
evolutionary equation for vector ${\bm u} = \{u_1, u_2, \dots, u_p \}$:
\begin{equation}\label{2.3}
\frac{d {\bm u}}{d t} + {\bm A} {\bm u} = 0,
\quad 0 < t \leq T,
\end{equation}
where for elements operator matrix ${\bm A}$ we have the representation
\[
{\bm A} = \{A_{\alpha \beta} \},
\quad  \alpha, \beta =1,2,\dots,p .
\]
On the direct sum of spaces
${\bm H} = H_1 \oplus  H_2 \oplus \cdots \oplus H_p$  put
\[
({\bm u}, {\bm v}) = \sum_{\alpha =1}^{p} (u_{\alpha},v_{\alpha})_{\alpha} ,
\quad \|{\bm u} \|^2 = \sum_{\alpha =1}^{p} \|u_{\alpha}\|^2_{\alpha} .
\]
Taking into account (\ref{2.2}), we have
\begin{equation}\label{2.4}
{\bm u}(0) = {\bm u}^0,
\end{equation}
where ${\bm u}^0 = \{u^0_1, u^0_2, \dots, u^0_p \}$.

We will consider the Cauchy problem (\ref{2.3}), (\ref{2.4}) provided that the operator ${\bm A}$ is self-adjoint and positive in ${\bm H}$:
\begin{equation}\label{2.5}
\bm A = \bm A^* > 0 .
\end{equation}
The self-adjoint property is associated with the fulfillment of the equalities
\[
A_{\alpha \beta} = A^*_{\beta \alpha},
\quad  \alpha, \beta =1,2,\dots,p
\]
for the operators of the original system of equations (\ref{2.1}).

For ${\bm D} = {\bm D}^*  > 0$ through ${\bm H}_{{\bm D}}$
denote the space ${\bm H}$ equipped with
scalar product $({\bm y},{\bm w})_{{\bm D}} = ({\bm D}{\bm y},{\bm w})$
and the norm $\|{\bm y}\|_{{\bm D}}= ({\bm D}{\bm y},{\bm y})^{1/2}$.
The simplest a priori estimate for the solution of the Cauchy problem (\ref{2.3}), (\ref{2.4}),
which we will be guided by when studying the corresponding
operator-difference schemes has the form
\begin{equation}\label{2.6}
\| {\bm u} (t) \|_{{\bm D}} \leq  \| {\bm u}^0 \|_{{\bm D}},
\quad 0 < t \leq T ,
\end{equation}
when $\bm D = \bm A^{-1}, \bm I, \bm A$ and ${\bm I}$ is the identity operator in ${\bm H}$.
It provides the stability of solutions of the problem (\ref{2.3}), (\ref{2.4}) to the initial data.

We introduce a uniform, for simplicity, time grid with step $\tau$ and let $\bm y^n = \bm y(t^n), \ t^n = n \tau$, $n = 0,1, \ldots, N, \ N\tau = T$.
For an approximate solution of the differential-operator problem (\ref{2.3}), (\ref{2.4}), we will use the usual weighted schemes \cite{Samarskii1989}.
When using a two-level scheme, equation (\ref{2.3}) is approximated by the difference equation
\begin{equation}\label{2.7}
\frac{{\bm y}^{n+1} - {\bm y}^{n}}{\tau } +
{\bm A} (\sigma {\bm y}^{n+1} + (1-\sigma) {\bm y}^{n} )  = 0,
\quad n = 0,1, \ldots, N-1 ,
\end{equation}
where $\sigma$ is a numeric parameter (weight) usually $0 \le \sigma \le 1$.
For simplicity, we restrict ourselves to the case of the same weight for all equations of the system (\ref{2.1}).
Taking into account (\ref{2.4}), we supplement (\ref{2.7}) with the initial condition
\begin{equation}\label{2.8}
{\bm y}^{0} = {\bm u}^{0}.
\end{equation}
If $\sigma \geq 1/2$, then the operator-difference scheme (\ref{2.7}) is unconditionally stable in  ${\bm H}_{{\bm A}}$ \cite{Samarskii1989,SamarskiiMatusVabischevich2002} and the difference solution satisfies the estimate
\begin{equation}\label{2.9}
\|{\bm y}^{n+1}\|_{\bm D} \leq  \|{\bm u}^{0}\|_{\bm D} ,
\quad n = 0,1, \ldots, N-1 ,
\quad  \bm D = \bm A^{-1}, \bm I, \bm A.
\end{equation}
The estimate (\ref{2.9}) acts as a grid analog of the estimate (\ref{2.6}) and ensures the unconditional stability of the difference scheme with weights (\ref{2.7}), (\ref{2.8}) under the natural constraints $\sigma \geq 1/2$.
Considering the corresponding problem for the error, we make sure that the solution of the operator-difference problem (\ref{2.7}), (\ref{2.8}) converges to the solution
differential-difference problem (\ref{2.3}), (\ref{2.4}) in ${\bm H}_{{\bm D}}$ for $\sigma \geq 1/2$ with $\mathcal{O}((2 \sigma -1)\tau + \tau^2)$.
For $\sigma = 1/2$, we have the second order of convergence in $\tau$.

The transition to a new time level in the scheme (\ref{2.7}) is associated with solving the problem
\[
({\bm I} + \sigma \tau {\bm A}) {\bm y}^{n+1} = \boldsymbol{\psi}^{n},
\quad  \boldsymbol{\psi}^{n} = {\bm y}^{n} - (1-\sigma) \tau {\bm A}) {\bm y}^{n} .
\]
About the original problem (\ref{2.1}), (\ref{2.2}), it is necessary to solve the system of coupled equations
\[
y_\alpha^{n+1} + \sum_{\beta =1}^{p} \sigma \tau A_{\alpha \beta} y_\beta^{n+1} = \psi_\alpha^n,
\quad \alpha = 1, 2, \dots, p.
\]
One or another iterative method can be used for an approximate solution of such a problem, particularly of the block type \cite{Saad2003}.

The second possibility more fully takes into account the specifics of the nonstationary problems under consideration and is associated with the construction of splitting schemes \cite{Marchuk1990,VabishchevichAdditive}, when the transition to a new time level is related to the solution of more straightforward problems.
For problems like (\ref{2.1}), (\ref{2.2}), it is natural to focus on splitting schemes when the transition to a new time level is provided by solving problems like
\[
y_\alpha^{n+1} +
\sigma \tau A_{\alpha \alpha} y_\alpha^{n+1} = \widetilde{\psi}_\alpha^n,
\quad \alpha = 1, 2, \dots, p.
\]
This realization corresponds to the fact that the computational implementation is provided by inverting only the diagonal part of the operator matrix
${\bm I} + \sigma \tau {\bm A}$.

\section{Explicit-implicit schemes}\label{s-3}

When constructing decoupling schemes, we can focus on extracting separate parts of the operator matrix ${\bm A}$.
The simplest variant is related to taking the diagonal part of the operator matrix ${\bm A}$ to the upper level in time.
We use the notation
\[
\bm D = \mathrm{diag} \{A_{11}, A_{22}, \ldots, A_{pp} \} .
\]
Our basic assumption is that there exists a constant $\gamma > 0$ such that
\begin{equation}\label{3.1}
\bm A \leq \gamma \bm D .
\end{equation}
An explicit-implicit time approximation is used to approximate the solution
\begin{equation}\label{3.2}
\frac{{\bm y}^{n+1} - {\bm y}^{n}}{\tau } +
\bm D (\sigma {\bm y}^{n+1} + (1-\sigma) {\bm y}^{n}) + (\bm A - \bm D) {\bm y}^{n} = 0.
\end{equation}

The scheme (\ref{2.8}), (\ref{3.2}) approximates (\ref{2.3}), (\ref{2.4}) with first order $\tau$ at any values of the parameter $\sigma$.
When investigating the stability of the constructed schemes, we rely on the general results of the theory of stability (correctness) of operator-difference schemes.
For two-level approximations, the following result of the theory of operator-difference schemes \cite{Samarskii1989,SamarskiiMatusVabischevich2002} is of crucial importance.

\begin{lem}\label{l-1}
Let ${\bm A}$ in
\begin{equation}\label{3.3}
\bm B \frac{{\bm y}^{n+1} - {\bm y}^{n}}{\tau } +
{\bm A} {\bm y}^{n} = 0 ,
\quad n = 0,1, \ldots, N-1 ,
\end{equation}
be a self-adjoint positive operator, and the operator ${\bm B}$ satisfy the condition
\begin{equation}\label{3.4}
{\bm B} \geq \frac {\tau}{2} {\bm A}  .
\end{equation}
Then the scheme (\ref{2.8}), (\ref{3.3}) is stable in ${\bm H}_{\bm A}$ and the a priori estimate
\begin{equation}\label{3.5}
\|\bm y^{n+1}\|_{\bm A} \leq \|\bm u^0 \|_{\bm A} ,
\quad \ n = 0,1, \ldots, N-1 ,
\end{equation}
is true.
\end{lem}

\pf
For the proof, let us write the scheme (\ref{3.3}) in the form of
\[
\left ( \bm B - \frac {\tau}{2} \bm A \right ) \frac{\bm y^{n+1} - \bm y^{n+1} - \bm y^{n}}{\tau } +
\bm A \frac{\bm y^{n+1} + \bm y^{n}}{2} = 0 .
\]
Let us multiply this equation scalarly in $\bm H$ by $2({\bm y}^{n+1} - {\bm y}^{n})$.
If (\ref{3.4}) is satisfied, we arrive at the inequality
\[
(\bm A \bm y^{n+1}, y^{n+1}) \leq (\bm A \bm y^{n}, y^{n}) ,
\]
from which the estimate (\ref{3.5}) follows.
\qed

Note that the formulated stability condition (\ref{3.4}) is sufficient and necessary.

\begin{thm}\label{t-1}
Explicit-implicit scheme (\ref{2.5}), (\ref{2.8}), (\ref{3.2}) under the inequality (\ref{3.1}) is stable in ${\bm H}_{\bm A}$ for
\begin{equation}\label{3.6}
2 \sigma \geq \gamma ,
\end{equation}
wherein the a priori estimate (\ref{3.5}) is correct.
\end{thm}

\pf
The scheme (\ref{3.2}) is written in the canonical form (\ref{3.3}), where, taking into account (\ref{3.1})
\[
\bm B = \bm I + \sigma \tau \bm D \geq \bm I + \frac{1}{\gamma } \sigma \tau A .
\]
The stability condition (\ref{3.4}) in lemma~\ref{l-1} will be satisfied when we choose the weight parameter $\sigma$ according to (\ref{3.6}).
\qed

For the scheme (\ref{3.3}), the approximate time solution at the new level is obtained from the solution of the problem
\[
\bm B \bm y^{n+1} = \bm B \bm y^{n} - \tau {\bm A} y^{n} .
\]
In the case of (\ref{3.2}), the individual components of the vector sought are found in Eqs.
\[
(I_\alpha + \sigma \tau A_{\alpha \alpha}) y^{n+1}_\alpha = f_\alpha^n ,
\quad \alpha =1,2,\ldots,p .
\]
For known right-hand sides of which we have
\[
f_\alpha^n = (I_\alpha + \sigma \tau A_{\alpha \alpha}) y^{n}_\alpha
- \tau \sum_{\beta=1}^{p}  A_{\alpha \beta} y_\beta^n ,
\quad \alpha =1,2,\ldots,p .
\]
It is possible to organize parallel calculations $y^{n+1}_\alpha, \ \alpha =1,2,\ldots,p$.

We easily solve the linear algebraic equations system with a diagonal and triangular matrix.
In the approximate solution of the Cauchy problem (\ref{2.3}), (\ref{2.4}), we can focus on time approximations with the triangular part of the operator matrix $\bm A$ taken to a new level in time. Let
\begin{equation}\label{3.7}
{\bm L} =
\begin{pmatrix}
0 & 0 & \cdots & 0 & 0 \\
A_{21} & 0 & \cdots & 0 & 0 \\
\cdots & \cdots & \cdots & \cdots & \cdots \\
A_{p1} & A_{p2} & \cdots & A_{p p-1} & 0 \\
\end{pmatrix} .
\end{equation}
Given the self-adjointness of $\bm A$, we have
\[
\bm A = \bm L + \bm D + \bm L^* .
\]
We will determine the approximate solution from
\begin{equation}\label{3.8}
\frac{{\bm y}^{n+1} - {\bm y}^{n}}{\tau } +
(\bm L + \bm D) {\bm y}^{n+1} + \bm L^* {\bm y}^{n} = 0.
\end{equation}

\begin{thm}\label{t-2}
Explicit-implicit scheme (\ref{2.5}), (\ref{2.8}), (\ref{2.8}), (\ref{3.7}), (\ref{3.8}) is stable in ${\bm H}_{\bm A}$ and there is an estimate of (\ref{3.5}) for the solution.
\end{thm}

\pf
Writing (\ref{3.8}) as (\ref{3.3}) we have
\[
\bm B = \bm I + \tau (\bm L + \bm D) .
\]
By
\[
(\bm B \bm u, \bm u) = \left ( \Big ( \bm I + \tau \frac{1}{2} (\bm L + \bm D + \bm L^* \bm L^* \Big ) \bm u, \bm u \right )
\geq \left ( \Big ( \bm I + \tau \frac{1}{2} \bm A \Big ) \bm u, \bm u \right ) ,
\]
the stability condition (\ref{3.4}) in lemma~\ref{l-1} is thus satisfied.
\qed

The individual components of the desired vector at a new level on time $y^{n+1}_\alpha, \ \alpha =1,2,\ldots,p$ are found sequentially from Eqs.
\[
y^{n+1}_\alpha + \tau \sum_{\beta=1}^{\alpha } A_{\alpha \beta} y_\beta^{n+1} = f_\alpha^n ,
\]
where
\[
f_\alpha^n = y^{n}_\alpha
- \tau \sum_{\beta=\alpha + 1}^{p}  A_{\alpha \beta} y_\beta^n ,
\quad \alpha =1,2,\ldots,p .
\]

The scheme (\ref{3.8}) as well as the scheme (\ref{3.2}) has first order $\tau$ accuracy.
Let us note the main possibilities of constructing second-order schemes to extract the diagonal and triangular parts of the operator matrix ${\bm A}$.

In the class of two-level schemes, we can focus on the schemes of the alternating-triangular method \cite{Samarskii1989,VabishchevichAdditive}.
They are based on the triangular splitting of the operator matrix ${\bm A}$:
\begin{equation}\label{3.9}
{\bm A} = {\bm A}_1 + {\bm A}_2,
\quad  {\bm A}_1^* = {\bm A}_2 .
\end{equation}
Additive presentation (\ref{3.9}) corresponds to the choice
\[
{\bm A}_1 =
\begin{pmatrix}
\frac{1}{2} A_{11} & 0 & \cdots & 0 \\
A_{21} & \frac{1}{2} A_{22} & \cdots & 0 \\
\cdots & \cdots & \cdots & \cdots \\
A_{p1} & A_{p2} & \cdots & \frac{1}{2} A_{pp} \\
\end{pmatrix} ,
\quad
{\bm A}_2 =
\begin{pmatrix}
\frac{1}{2} A_{11} & A_{12} & \cdots & A_{1p} \\
0 & \frac{1}{2} A_{22} & \cdots & A_{2p} \\
\cdots & \cdots & \cdots & \cdots \\
0 & 0 & \cdots & \frac{1}{2} A_{pp} \\
\end{pmatrix}  .
\]

Instead of (\ref{3.8}) we will use the scheme (\ref{3.3}), in which
\begin{equation}\label{3.10}
\bm B = \bm I + \sigma \tau \bm A_1 .
\end{equation}
Similarly to Theorem~\ref{t-2}, we prove the stability of the scheme (\ref{2.5}), (\ref{2.8}), (\ref{3.3}), (\ref{3.10}) at $\sigma \geq 1$ in ${\bm H}_{{\bm A}}$.
A more interesting possibility is associated with the choice of the operator $\bm B$, which is represented in the following factorized form
\begin{equation}\label{3.11}
\bm B =
({\bm I} + \sigma \tau {\bm A}_1) ({\bm I} + \sigma \tau {\bm A}_2) .
\end{equation}
When $\sigma = 1/2$, the scheme (\ref{3.9}), (\ref{3.11}) approximates (\ref{2.3}) with second order on $\tau$.

\begin{thm}\label{t-3}
The factorized scheme (\ref{2.5}), (\ref{2.8}), (\ref{3.3}), (\ref{3.9}), (\ref{3.11}) is unconditionally stable in  ${\bm H}_{{\bm A}}$ for $\sigma \geq 1/2$. The difference solution satisfies the estimate (\ref{3.5}).
\end{thm}

\pf
In the case (\ref{3.9}), (\ref{3.11}) we have
\[
\bm B = {\bm I} + \sigma\tau {\bm A} + \sigma^2\tau^2 {\bm A}_1 {\bm A}_2,
\quad \bm B = \bm B^* \geq  {\bm I} + \sigma\tau {\bm A} .
\]
Thus, the inequality (\ref{3.4}) holds for $\sigma \geq 1/2$.
Considering this, the required statement follows from Lemma~\ref{l-1}.
\qed

The construction of explicit-implicit schemes of the second order of accuracy in extracting the diagonal part of the operator matrix $\bm A$ can be done using three-level schemes \cite{Samarskii1989,SamarskiiMatusVabischevich2002}.
The three-level scheme for the approximate solution of the problem (\ref{2.3})--(\ref{2.5}) can be written in the form
\begin{equation}\label{3.12}
\begin{split}
{{\bm B}} \frac{{\bm y}^{n+1} - {\bm y}^{n-1}}{2\tau } & +
\bm R ({\bm y}^{n+1} - 2 {\bm y}^{n} + {\bm y}^{n-1}) +
{\bm A} {\bm y}^{n} = 0 , \\
& \quad n = 1,2, \ldots, N-1,
\end{split}
\end{equation}
for given
\begin{equation}\label{3.13}
\bm y^0 = \bm u^0,
\quad \bm y^1 = \bm v^0 .
\end{equation}
Let us formulate the stability conditions for the scheme (\ref{3.13}) in the form of the following statement \cite{Samarskii1989,SamarskiiMatusVabischevich2002}.

\begin{lem}\label{l-2}
Let in (\ref{3.12}) the operators ${\bm R}$ and ${\bm A}$ be constant self-adjoint operators.
When conditions
\begin{equation}\label{3.14}
{\bm B} \geq 0,
\quad \bm A > 0,
\quad \bm R > \frac {1}{4} {\bm A}
\end{equation} are true,
the scheme (\ref{3.12}), (\ref{3.13}) is stable and the a priori estimate
\begin{equation}\label{3.15}
{\cal E}^{n+1} \leq {\cal E}^1
\end{equation}
holds, in which
\[
\begin{split}
{\cal E}^{n+1} & = \frac 14 (\bm A(\bm y^{n+1}+ \bm y^n), \bm y^{n+1} + \bm y^n) \\
& + (\bm R(\bm y^{n+1}- \bm y^n), \bm y^{n+1}- \bm y^n) -
\frac 14 (\bm A(\bm y^{n+1} - \bm y^n), \bm y^{n+1} - \bm y^n).
\end{split}
\]
\end{lem}

The proof of the estimation is based on the introduction of new quantities
\[
\bm s^n = \frac{1}{2} (\bm y^n + \bm y^{n-1}) ,
\quad \bm r^n = \bm y^n - \bm y^{n-1} .
\]
Given
\[
\bm y^n = \frac{1}{4} (\bm y^{n+1} + 2 \bm y^{n+1} + 2 \bm y^{n-1}) - \frac{1}{4} (\bm y^{n+1} - 2 \bm y^n + \bm y^{n-1}) ,
\]
we can write (\ref{3.12}) as
\[
\bm B \frac{\bm r^{n+1} + \bm r^{n}}{2 \tau } + \left ( \bm R - \frac{1}{4} \bm A \right ) (\bm r^{n+1} - \bm r^{n})
+ \bm A \frac{\bm s^{n+1} + \bm s^n}{2 } = \bm \varphi^{n} ,
\quad \ n = 1,2, \ldots, N-1 .
\]
Multiplying this equation scalarly in $\bm H$ by $2(\bm s^{n+1} - \bm s^{n}) = (\bm r^{n+1} + \bm r^{n})$, we get an inequality
\[
{\cal E}^{n+1} \leq {\cal E}^{n} ,
\]
from which follows the stability estimate from the initial data (\ref{3.15}).

For the approximate solution of (\ref{2.3}), (\ref{2.4}), we can use a three-level scheme instead of the two-level scheme (\ref{3.2})
of the second approximation order
\begin{equation}\label{3.16}
\frac{{\bm y}^{n+1} - {\bm y}^{n-1}}{2\tau } +
\bm D (\sigma {\bm y}^{n+1} + (1-2\sigma) {\bm y}^{n} + \sigma {\bm y}^{n-1}) + (\bm A - \bm D) {\bm y}^{n} = 0 .
\end{equation}

\begin{thm}\label{t-4}
Explicit-implicit three-level scheme (\ref{2.5}), (\ref{3.13}), (\ref{3.16}) with constraints (\ref{3.1}) is stable for
\begin{equation}\label{3.17}
4 \sigma > \gamma ,
\end{equation}
and the a priori estimate (\ref{3.15}) is true with
\[
\bm B = \bm I,
\quad \bm R = \sigma \bm D .
\]
\end{thm}

\pf
The three-level scheme (\ref{3.16}) is written as (\ref{3.16}).
Taking into account (\ref{3.1}), the last inequality (\ref{3.14}) will be satisfied under the constraints (\ref{3.17}) on the weight $\sigma$.
Thus, the conditions of the lemma~\ref{l-2} are satisfied.
\qed

The explicit-implicit approximations under consideration can be interprete as  regularized schemes \cite{VabishchevichAdditive}.
As a primary one, we take a weighted scheme (\ref{2.7}), which satisfies us both in accuracy (first order for all $\sigma \neq 1/2$ and second order --- for $\sigma = 1/2$) and stability (it is enough to take $\sigma \geq 1/2$).
Regularization by perturbations of the difference scheme operators aims to obtain a good decoupling system of equations for finding an approximate solution at a new time level.

In particular, the scheme (\ref{3.9}), (\ref{3.11}) can be written as follows:
\begin{equation}\label{3.18}
(\bm I + \sigma \tau \bm A) \frac{{\bm y}^{n+1} - {\bm y}^{n}}{\tau } +
\sigma^2 \tau^2 \bm A_1 \bm A_2 \frac{{\bm y}^{n+1} - {\bm y}^{n}}{\tau } + +
{\bm A} {\bm y}^{n} = 0 .
\end{equation}
The regularization is provided by the summand $\sigma^2 \tau \bm A_1 \bm A_2 ({\bm y}^{n+1} - {\bm y}^{n} ) = \mathcal{O}(\tau^2)$.
In \cite{vabishchevich2014PTM}, three-level variants of the alternating-triangle method schemes are proposed.
Instead of (\ref{3.18}) we will use the scheme
\begin{equation}\label{3.19}
(\bm I + \sigma \tau \bm A) \frac{{\bm y}^{n+1} - {\bm y}^{n}}{\tau } +
\sigma^2 \tau \bm A_1 \bm A_2 ({\bm y}^{n+1} - 2 {\bm y}^{n} + {\bm y}^{n-1}) +
{\bm A} {\bm y}^{n} = 0 ,
\end{equation}
which is related to the choice of the regularizing summand $\sigma^2 \tau \bm A_1 \bm A_2 ({\bm y}^{n+1} - 2 {\bm y}^{n} + {\bm y}^{n-1}) = \mathcal{O}(\tau^3)$.

Taking into account that
\[
\frac{{\bm y}^{n+1} - {\bm y}^{n}}{\tau } = \frac{{\bm y}^{n+1} - {\bm y}^{n-1}}{2\tau }
+ \frac{{\bm y}^{n+1} - 2 {\bm y}^{n} + {\bm y}^{n-1}}{2\tau } ,
\]
we write the scheme (\ref{3.19}) in the form (\ref{3.12}) for
\[
\bm B = \bm I + \sigma \tau \bm A ,
\quad \bm R = \frac{1}{2 \tau }  (\bm I + \sigma \tau \bm A) + \sigma^2 \tau \bm A_1 \bm A_2 .
\]
For $\sigma \geq 1/2$, conditions (\ref{3.14}) are satisfied.
Lemma~\ref{l-2} allows us to state the following statement.

\begin{thm}\label{t-5}
For $\sigma \geq 1/2$, the scheme (\ref{2.5}), (\ref{3.9}), (\ref{3.13}), (\ref{3.19}) is unconditionally stable, and the a priori estimate (\ref{3.5}) is valid for the approximate solution.
\end{thm}

The main disadvantages of the three-level schemes in the approximate solution of the Cauchy problem for first-order evolution equations are well known:
\begin{enumerate}
\item calculation of the first time step;
\item using non-uniform time grids.
\end{enumerate}
We can use a two-level time scheme to set the second initial condition (\ref{3.13}).
In the simplest version, $\bm v^0$ is determined from Eq.
\[
(\bm I + \sigma \tau \bm A) \frac{{\bm v}^0 - {\bm y}^{0}}{\tau } + \bm A y^0 = 0 .
\]
The problem of non-uniform meshes for three-level schemes is discussed in the literature (see, for example, the article \cite{vabishchevich2023nonunif} and the bibliography given therein), but the direct use of the obtained results for the considered problems is complex.

\section{Multicomponent decomposition-composition schemes}\label{s-4}

Above, we have noted the possibilities of constructing splitting schemes for two-component splitting of the operator matrix $\bm A$ by extracting its diagonal part or splitting it into two triangular matrices.
In the theory and practice \cite{Marchuk1990,VabishchevichAdditive} of splitting schemes, one often focuses on the multicomponent additive representation of the problem operator.

We use a general computational technique of decomposition and composition \cite{DC} for the approximate solution of the Cauchy problem for nonstationary problems.
At the decomposition (analysis) stage, the problem operator $\bm A$ is represented as a sum of operators $\bm A_\alpha, \ \alpha = 1,2, \ldots, p$:
\begin{equation}\label{4.1}
\bm A = \sum_{\alpha = 1}^{p} \bm A_\alpha ,
\quad \alpha=1,2,\ldots,p .
\end{equation}
The problem for a single-operator summand is more convenient for computational realization.
At the composition stage (synthesis), the approximate solution at the new time level is constructed as a linear superposition of the solutions of the auxiliary problems.

The additive decomposition of the problem operator is based on the representation of the unit operator $\bm I$ in $\bm H$ in the form of
\begin{equation}\label{4.2}
\bm I = \sum_{\alpha = 1}^{p} \bm R_\alpha ,
\quad \bm R_\alpha = \bm R_\alpha^* \geq 0,
\quad \alpha=1,2,\ldots,p .
\end{equation}
To derive (\ref{4.1}) from the decomposition of the unit operator (\ref{4.2}), we have two main possibilities.
In the first possibility, we have
\begin{equation}\label{4.3}
\bm A \longrightarrow \sum_{\alpha = 1}^{p} \bm R_\alpha \bm A,
\quad \bm A_\alpha = \bm R_\alpha \bm A,
\quad \alpha=1,2,\ldots, p .
\end{equation}
In the second option, we use a slightly different construction:
\begin{equation}\label{4.4}
\bm A \longrightarrow \sum_{\alpha = 1}^{p} \bm A \bm R_\alpha ,
\quad \bm A_\alpha = \bm A \bm R_\alpha ,
\quad \alpha=1,2,\ldots, p .
\end{equation}
This approach was used earlier in constructing domain decomposition schemes for the numerical solutions of parabolic problems \cite{vabishchevich2023subdomain}, and as a general methodological technique, was formulated in  \cite{DC}.

In the decomposition (\ref{4.2}), we will put
\begin{equation}\label{4.5}
\bm R_\alpha = \operatorname{diag} (0, \ldots, 0, I_\alpha, 0, \ldots, 0),
\quad \alpha=1,2,\ldots, p ,
\end{equation}
where $I_\alpha$ is a unit operator in $H_\alpha, \ \alpha = 1,2, \ldots, p$.
Under this construction, the operators $\bm R_\alpha, \ \alpha = 1,2, \ldots, p$ are orthogonal projectors in the corresponding subspaces of $\bm H$.

Using (\ref{4.5}), the first decomposition variant (\ref{4.3}) leads to an additive representation (\ref{4.1}) in which the operator summands are associated with the selection of individual rows of the operator matrix $\bm A$:
\begin{equation}\label{4.6}
\bm A_\alpha = \begin{pmatrix}
0 & 0 & \cdots & 0 \\
\cdots & \cdots & \cdots & \cdots \\
0 & 0 & \cdots & 0 \\
A_{\alpha 1} & A_{\alpha 2} & \cdots & A_{\alpha p} \\
0 & 0 & \cdots & 0 \\
\cdots & \cdots & \cdots & \cdots \\
0 & 0 & \cdots & 0 \\
\end{pmatrix} ,
\quad \alpha=1,2,\ldots, p .
\end{equation}
The second variant of decomposition (\ref{4.4}) leads to a split of $\bm A$ into separate columns:
\begin{equation}\label{4.7}
\bm A_\alpha = \begin{pmatrix}
0 & \cdots & 0 & A_{1 \alpha} & 0 &  \cdots & 0 \\
\cdots & \cdots & \cdots & \cdots& \cdots & \cdots & \cdots \\
0 & \cdots & 0 & A_{\alpha \alpha} & 0 &  \cdots & 0 \\
\cdots & \cdots & \cdots & \cdots & \cdots & \cdots & \cdots \\
0 & \cdots & 0 & A_{p \alpha} & 0 & \cdots & 0 \\
\end{pmatrix} ,
\quad \alpha=1,2,\ldots, p .
\end{equation}

At the composition stage, the approximate solution of the problem (\ref{2.3}), (\ref{2.4}) is constructed based on the use of additive splitting (\ref{4.1}).
We first consider the time approximation problem for the decomposition of (\ref{4.2}), (\ref{4.3}).
From (\ref{2.3}), (\ref{4.1}), (\ref{4.2}), (\ref{4.4}), we obtain Eq.
\begin{equation}\label{4.8}
\frac{d \bm u}{d t} + \sum_{\alpha = 1}^{p} \bm R_\alpha \bm A \bm u = 0,
\quad 0 < t \leq T .
\end{equation}
We cannot simply use the results of the theory of additive operator-difference schemes in constructing splitting schemes for the problem (\ref{2.4}), (\ref{4.8}) because the operators $\bm R_\alpha \bm A, \ \alpha =1,2, \ldots, p$ are not non-negative and are not self-adjoint.

Let us symmetrize  \cite{SamarskiiMatusVabischevich2002,SamVabGulin} the operator of the equation (\ref{4.8}).
Instead of $\bm u$ we introduce a new quantity $\bm v = \bm A^{1/2} \bm u$.
From equation (\ref{4.8}) we obtain
\begin{equation}\label{4.9}
\frac{d \bm v}{d t} + \sum_{\alpha = 1}^{p} \bm A_\alpha \bm v = 0,
\quad 0 < t \leq T .
\end{equation}
For the operator summands $\bm A_\alpha, \ \alpha =1,2, \ldots, p$ we have
\begin{equation}\label{4.10}
\bm A_\alpha = \bm A^{1/2} \bm R_\alpha \bm A^{1/2} = \bm A_\alpha^* \geq 0,
\quad \alpha =1,2, \ldots, p .
\end{equation}
We can now construct composition schemes for equation (\ref{4.9}), which is augmented by the initial condition
\begin{equation}\label{4.11}
\bm v(0) = \bm v^0 .
\end{equation}
For the problem (\ref{2.4}), (\ref{4.8}) we have $\bm v^0 = \bm A^{1/2} \bm u^0$.

At the second variant of decomposition (\ref{4.2}), (\ref{4.4}) we come to the equation
\begin{equation}\label{4.12}
\frac{d \bm u}{d t} + \sum_{\alpha = 1}^{p} \bm A \bm R_\alpha \bm u = 0,
\quad 0 < t \leq T .
\end{equation}
From the problem (\ref{2.4}), (\ref{4.12}) we get to the problem (\ref{4.9})--(\ref{4.11}) for $\bm v = \bm A^{-1/2} \bm u$.

For the problem (\ref{4.9})--(\ref{4.11}), the standard stability estimate on the initial data is as follows
\begin{equation}\label{4.13}
\|\bm v(t)\| \leq \|\bm v^0\|,
\quad 0 < t \leq T .
\end{equation}
In the first decomposition of (\ref{2.4}), (\ref{4.8}) the estimate of (\ref{4.13}) corresponds to the estimate of (\ref{2.6}) at $\bm D = \bm A$, and in the second decomposition variant (\ref{2.4}), (\ref{4.12}) from (\ref{4.13}) we arrive at the estimate (\ref{2.6}) at $\bm D = \bm A^{-1}$.

Various additive operator-difference schemes \cite{VabishchevichAdditive} can be used for approximate solutions of the problem (\ref{4.9})--(\ref{4.11}).
We will distinguish factorized schemes under two-component splitting ($p=2$), acting as operator analogies of classical alternating-direction schemes. Such composition schemes under more interesting decomposition variants have been discussed above.
In computational practice, component-wise splitting schemes of first and second-order time approximation for $p \ge 2$ are the most widely used.
A complicating factor is that such composition schemes are sum approximation schemes since the individual equations do not approximate the original equations.
Additive-averaging schemes can be viewed as parallel variants of regularized component-wise splitting schemes.
Regularized additive full approximation schemes are most easily constructed based on multiplicative perturbation of individual operator terms.
Let us also note vector additive schemes, the construction of which is based on the transition to a system of equivalent equations.
For the problems  (\ref{2.4}), (\ref{4.8}) and (\ref{2.4}), (\ref{4.12}), which we consider, they are of no interest since we initially have a system of equations.

We begin our consideration of composition schemes for the problem (\ref{4.9})--(\ref{4.11}) with the usual component-wise splitting schemes.
We formulate an auxiliary result on the stability of two-level schemes with weights.

\begin{lem}\label{l-3}
Let the constant operators
\begin{equation}\label{4.14}
 \bm B = \bm B^* > 0,
 \quad \bm A = \bm A^* \geq 0 
\end{equation} 
in the scheme 
\begin{equation}\label{4.15}
  \bm B \frac{{\bm y}^{n+1} - {\bm y}^{n}}{\tau } +
  \bm A (\sigma \bm y^{n+1} +(1-\sigma) \bm y^{n} ) = 0 ,
  \quad n = 0,1, \ldots, N-1 .
\end{equation} 
Then the scheme (\ref{4.14}), (\ref{4.15}) is stable in ${\bm H}_{\bm B}$ when $\sigma = \operatorname{const} \geq 1/2$ and the estimate 
\[
\|\bm y^{n+1}\|_{\bm B} \leq \|\bm u^0 \|_{\bm B} ,
\quad \ n = 0,1, \ldots, N-1 ,
\]
holds.
\end{lem}

\pf
For $\bm y^{n+\sigma} = \sigma \bm y^{n+1} +(1-\sigma) \bm y^{n}$ we have
\[
 \bm y^{n+\sigma} = \left (\sigma - \frac{1}{2} \right ) (\bm y^{n+1} - \bm y^{n}) + \frac{1}{2} (\bm y^{n+1} + \bm y^{n}) .
\] 
Let us multiply (\ref{4.15}) scalarly by $\bm y^{n+\sigma}$. Under the constraints $\sigma \geq 1/2$ and taking into account (\ref{4.14}), we obtain the inequalities 
\[ 
 (\bm B \bm y^{n+1}, \bm y^{n+1}) \leq (\bm B \bm y^{n}, \bm y^{n}),
 \quad n = 0,1, \ldots, N-1 , 
\] 
from which the stability estimate follows.
\qed

The approximate solution of the problem (\ref{4.9})--(\ref{4.11}) at a new level in time is determined from the sequential solution of $p$ problems \cite{Marchuk1990,VabishchevichAdditive}:
\begin{equation}\label{4.16}
\begin{split}
\frac{\bm w^{n+\alpha/p} - \bm w^{n+(\alpha-1)/p}}{\tau }
+ & \bm A_\alpha (\sigma \bm w^{n+\alpha/p} + (1-\sigma) \bm w^{n+(\alpha-1)/p}) = 0, \\
& \alpha = 1,2, \ldots, p ,
\quad n = 0,1, \ldots, N-1 ,
\end{split}
\end{equation}
\begin{equation}\label{4.17}
\bm w^0 = \bm v^0 .
\end{equation}
The main result regarding the stability of initial data is formulated as follows.

\begin{thm}\label{t-6}
The additive operator-difference scheme (\ref{4.16}), (\ref{4.17}) for the approximate solution of the problem (\ref{4.9}), (\ref{4.11}) is unconditionally stable under the splitting (\ref{4.10}) and $\sigma \geq 1/2$, and the estimation
\begin{equation}\label{4.18}
\|{\bm w}^{n+1}\| \leq \|{\bm v}^{0}\| ,
\quad n = 0,1, \ldots, N-1 ,
\end{equation}
holds.
\end{thm}

\pf
If we fulfill the conditions of the theorem from lemma~\ref{l-3} for the auxiliary quantities, we obtain the inequality
\[
\|\bm w^{n+\alpha/p}\| \leq \|\bm w^{n+(\alpha-1)/p}\|,
\quad \alpha = 1,2, \ldots, p .
\]
Thus, when we move to a new level in time, we get the inequality
\[
\|\bm w^{n+1}\| \leq \|\bm w^{n}\| .
\]
This inequality ensures that the estimate (\ref{4.18}) is satisfied.
\qed

In problem (\ref{2.4}), (\ref{4.8}), the component-wise splitting scheme (\ref{4.16}), (\ref{4.17}) corresponds to the definition of an approximate solution from Eqs.
\begin{equation}\label{4.19}
\begin{split}
\frac{\bm y^{n+\alpha/p} - \bm y^{n+(\alpha-1)/p}}{\tau }
+ & \bm R_\alpha \bm A  (\sigma \bm y^{n+\alpha/p} + (1-\sigma) \bm y^{n+(\alpha-1)/p}) = 0, \\
& \alpha = 1,2, \ldots, p ,
\quad n = 0,1, \ldots, N-1 ,
\end{split}
\end{equation}
when the initial condition (\ref{2.8}) is given.
Theorem~\ref{t-6} provides stability in $\bm H_{\bm A}$ --- estimate (\ref{3.5}).

For the problem (\ref{2.4}), (\ref{4.12}) we use the scheme
\begin{equation}\label{4.20}
\begin{split}
\frac{\bm y^{n+\alpha/p} - \bm y^{n+(\alpha-1)/p}}{\tau }
+ & \bm A \bm R_\alpha  (\sigma \bm y^{n+\alpha/p} + (1-\sigma) \bm y^{n+(\alpha-1)/p}) = 0, \\
& \alpha = 1,2, \ldots, p ,
\quad n = 0,1, \ldots, N-1 .
\end{split}
\end{equation}
The stability of the approximate solution holds in $\bm H_{\bm A^{-1}}$:
\begin{equation}\label{4.21}
\|\bm y^{n+1}\|_{\bm A^{-1}} \leq \|\bm u^0 \|_{\bm A^{-1}} ,
\quad \ n = 0,1, \ldots, N-1 .
\end{equation}

Let us note the key points of the computational realization of the scheme (\ref{2.8}), (\ref{4.19}).
To approximate the solution at the new time level $\bm y^{n+\alpha/p} = \bm y^{n+(\alpha-1)/p} + \bm z$, we solve $p$  equations
\begin{equation}\label{4.22}
\bm B \bm z = \bm f ,
\end{equation}
in which
\[
\bm B = \bm I+\sigma \tau \bm R_\alpha \bm A ,
\quad \bm f =  \tau \bm R_\alpha \bm A \bm y^{n+(\alpha-1)/p} ,
\quad \alpha = 1,2, \ldots, p .
\]
The non-zero elements of the matrix $\bm B$ lie on the main diagonal and the row numbered $\alpha$:
\[
B_{\alpha \alpha} = I_\alpha + \sigma \tau A_{\alpha \alpha},
\quad  B_{\beta \beta  } = I_\beta ,
\quad B_{\alpha \beta } = \sigma \tau A_{\alpha \beta},
\quad \beta \neq \alpha ,
\quad \beta =1,2,\ldots,p .
\]
For non-zero components of the right-hand side, we have
\[
f_\alpha  = \tau \sum_{\beta = 1}^{p} A_{\alpha \beta} y_\beta^{n+(\alpha-1)/p} ,
\quad \alpha  =1,2,\ldots,p .
\]
An algorithm for finding an approximate solution on a new level by time:
\begin{enumerate}
\item Find $z_\alpha$ from Eq.
\[
(I_\alpha + \sigma \tau A_{\alpha \alpha}) z_\alpha = \tau \sum_{\beta = 1}^{p} A_{\alpha \beta} y_\beta^n ;
\]
\item All other components of $\bm z$ are zero:
\[
z_\beta = 0,
\quad \beta \neq \alpha ,
\quad \beta =1,2,\ldots,p .
\]
\end{enumerate}
The computational work to find $z_\alpha$ is associated with solving an equation that involves the operator $A_{\alpha \alpha}$ at each $\alpha=1,2,\ldots,p$.
Algorithmically, the implementation of the decomposition-composition scheme (\ref{2.8}), (\ref{4.19}) is close to the implementation of the scheme (\ref{2.8}), (\ref{3.8}).
The differences can be traced when considering the most straightforward problem with $p=2$ (see section~\ref{s-5} with numerical experiments).

Using the composition scheme (\ref{2.8}), (\ref{4.20}) we solve equation (\ref{4.22}), where now
\[
\bm B = \bm I+\sigma \tau \bm A \bm R_\alpha ,
\quad \bm f =  \tau \bm A \bm R_\alpha \bm y^{n+(\alpha-1)/p} ,
\quad \alpha = 1,2, \ldots, p .
\]
The non-zero elements of the matrix $\bm B$ lie on the main diagonal, and the column numbered $\alpha$, and
\[
B_{\alpha \alpha} = I_\alpha + \sigma \tau A_{\alpha \alpha},
\quad  B_{\beta  \beta } = I_\beta ,
\quad B_{\beta \alpha} = \sigma \tau A_{\beta  \alpha},
\quad \beta \neq \alpha ,
\quad \beta =1,2,\ldots,p .
\]
For the components of the right-hand side of the equation (\ref{4.22}) we have
\[
f_\beta = \tau A_{\beta \alpha} y_\alpha^{n+(\alpha-1)/p},
\quad \beta =1,2,\ldots,p .
\]
The calculations are performed in the following order:
\begin{enumerate}
\item Find $z_\alpha$ from Eq.
\[
(I_\alpha + \sigma \tau A_{\alpha \alpha}) z_\alpha = \tau A_{\alpha \alpha} y_\alpha^{n+(\alpha-1)/p} ;
\]
\item For the other components of $\bm z$, we have
\[
z_\beta = \tau A_{\beta \alpha} (y_\alpha^{n+(\alpha-1)/p} - \sigma z_\alpha),
\quad \beta \neq \alpha ,
\quad \beta =1,2,\ldots,p .
\]
\end{enumerate}
The main computational work to find $\bm y^{n+\alpha/p}$ is generated by solving an equation that involves the operator $A_{\alpha \alpha}$ at each $\alpha=1,2,\ldots,p$.
The other elements of the operator matrix $A_{\beta \alpha}, \ \beta \neq \alpha$ are associated with explicit computations.

Let us consider the second class of composition schemes for the approximate solution of (\ref{4.9})--(\ref{4.11}).
We will apply regularized additive schemes \cite{VabishchevichAdditive,SamarskiiVabischevich1998} when the approximate solution at a new level in time is determined from Eq.
\begin{equation}\label{4.23}
\frac{\bm w^{n+1} - \bm w^n}{\tau} + \sum_{\alpha = 1}^{p} (\bm I+\sigma \tau \bm A_\alpha )^{-1}  \bm A_\alpha \bm w^n = 0 .
\end{equation}
Here $\sigma$ is a numerical parameter (weight).
The scheme (\ref{4.17}), (\ref{4.23}) is a full approximation scheme and approximates (\ref{4.9}), (\ref{4.11}) with first order $\tau$.

\begin{thm}\label{t-7}
The scheme (\ref{4.10}), (\ref{4.17}), (\ref{4.23}) is unconditionally stable at $\sigma \geq p/2$ in $\bm H$.
For the approximate solution, there is a stability estimate on the initial data (\ref{4.18}).
\end{thm}

\pf
We use a special additive representation of the solution at the new level:
\begin{equation}\label{4.24}
\bm w^{n+1} = \frac{1}{p} \sum_{\alpha = 1}^{p} \bm w_\alpha^{n+1} .
\end{equation}
The auxiliary quantities $\bm w_\alpha^{n+1}$ are determined from the system of equations
\begin{equation}\label{4.25}
(\bm I+\sigma \tau \bm A_\alpha ) \frac{\bm w_\alpha^{n+1} - \bm w^n}{p \tau} + \bm A_\alpha \bm w^n = 0 ,
\quad \alpha=1,2,\ldots,p .
\end{equation}
Taking into account (\ref{4.10}) and applying lemma~\ref{l-3}, for $\bm w_\alpha^{n+1}$ we obtain
\[
\|\bm w_\alpha^{n+1}\| \leq \|\bm w^n\| ,
\quad \alpha=1,2,\ldots,p ,
\]
at $\sigma \geq p/2$.
From the representation (\ref{4.24}), the estimation of the transition to a new level in time is as follows
\[
\|\bm w^{n+1}\| \leq \frac{1}{p} \sum_{\alpha = 1}^{p} \|\bm w_\alpha^{n+1} \| \leq \|\bm w^n \| .
\]
This estimation allows us to establish a stability estimate (\ref{4.18}).
\qed

By introducing the auxiliary quantities $\bm w_\alpha^{n+1}, \ \alpha=1,2,\ldots,p$ according to (\ref{4.24}), we have the additive-averaged scheme (\ref{4.17}), (\ref{4.24}), (\ref{4.25}) for the approximate solution of the problem (\ref{4.9})--(\ref{4.11}).
Thus, the regularized scheme (\ref{4.10}), (\ref{4.17}), (\ref{4.23}) is a more general composition scheme.

The organization of the computation is similar to that of the component-wise splitting schemes (\ref{4.16}), (\ref{4.17}).
For example, under the decomposition (\ref{4.2}), (\ref{4.3}), the solution at the new level in time is defined by
\[
\frac{\bm y^{n+1} - \bm y^n}{\tau} + \sum_{\alpha = 1}^{p} (\bm I+\sigma \tau \bm R_\alpha \bm A )^{-1} \bm R_\alpha \bm A \bm y^n = 0 .
\]
First, similarly to (\ref{4.22}) find the auxiliary quantities $\bm z_\alpha, \ \alpha=1,2,\ldots,p$, from Eqs.
\[
\bm B_\alpha \bm z_\alpha = \bm f_\alpha ,
\]
in which
\[
\bm B_\alpha = \bm I+\sigma \tau \bm R_\alpha \bm A ,
\quad \bm f_\alpha = \tau \bm R_\alpha \bm A \bm y^{n} ,
\quad \alpha = 1,2, \ldots, p .
\]
For the approximate solution at the new level, we have
\[
\bm y^{n+1} = \bm y^n - \sum_{\alpha = 1}^{p}\bm z_\alpha .
\]
The most important aspect of regularized additive schemes is the possibility of parallel (independent) calculation of these quantities.

Accuracy improvement up to the second order of the component-wise splitting scheme (\ref{4.16}), (\ref{4.17}) is carried out based on symmetrization when organizing calculations according to the Fryazinov-Strang regulation \cite{Fryazinov1968,Strang1968}
\[
\bm A_1 \rightarrow \bm A_2 \rightarrow \cdots \rightarrow \bm A_p \rightarrow \bm A_p \rightarrow \bm A_{p-1} \rightarrow \cdots \rightarrow \bm A_1.
\]
and choosing $\sigma = 1/2$.
Instead of (\ref{4.9}), (\ref{4.10}) we will use the equation
\[
\frac{d \bm v}{d t} + \sum_{\alpha = 1}^{2p} \bm A_\alpha \bm v = 0,
\quad 0 < t \leq T ,
\]
with operator summands
\[
\bm A_\alpha = \bm A_\alpha^* \geq 0 ,
\quad \bm A_\alpha = \frac{1}{2} \left \{  \begin{array}{ll}
\bm A^{1/2} \bm R_\alpha \bm A^{1/2} ,   & \alpha =1,2, \ldots, p , \\
\bm A^{1/2} \bm R_{2p+1-\alpha}\bm A^{1/2} ,   & \alpha =p+1,p+2, \ldots, 2p . \\
\end{array}
\right .
\]
Similar to (\ref{4.16}), for the approximate solution, we use the scheme
\begin{equation}\label{4.26}
\begin{split}
\frac{\bm w^{n+\alpha/2p} - \bm w^{n+(\alpha-1)/2p}}{\tau }
+ & \bm A_\alpha \frac{\bm w^{n+\alpha/p} + \bm w^{n+(\alpha-1)/p}}{2} = 0, \\
& \alpha = 1,2, \ldots, 2p ,
\quad n = 0,1, \ldots, N-1 ,
\end{split}
\end{equation}
which is unconditionally stable in $\bm H$ and has second-order accuracy on $\tau$.

We have some possibilities for constructing regularized additive schemes of second-order accuracy schemes related to using three-level time approximations, noted in \cite{VabishchevichAdditive,SamarskiiVabischevich1998}.
For our problems (\ref{2.1}), (\ref{2.2}) with operator matrices $\bm A$, the known variants seem to be too heavy.

\section{Numerical experiments}\label{s-5}

We illustrate the possibilities of the proposed decomposition-composition schemes on the model diffusion problem.
In the unit square $\Omega= \{\bm{x}\mid \bm{x}=(x_1,x_2),\ 0<x_1, x_2 <1\}$ we find the solution $w_\alpha(\bm x, t), \ \alpha =1,2$ of the system of equations
\begin{equation}\label{5.1}
\begin{split}
\frac{\partial w_{1}}{\partial t} & - \operatorname{div} (d_{1 1} ({\bm x}) \operatorname{grad} \, w_{1})
- \operatorname{div} (d_{1 2} ({\bm x}) \operatorname{grad} \, w_{2}) = 0 , \\
\frac{\partial w_{2}}{\partial t} & - \operatorname{div} (d_{2 1} ({\bm x}) \operatorname{grad} \, w_{1})
- \operatorname{div} (d_{2 2} ({\bm x}) \operatorname{grad} \, w_{2}) = 0 ,
\quad  {\bm x} \in \Omega,
\quad  0 < t \leq T .
\end{split}
\end{equation}
The boundary and initial conditions are
\begin{equation}\label{5.2}
\frac{\partial w_{1}}{\partial \nu} (\bm x,t)  = 0,
\quad  \frac{\partial w_{2}}{\partial \nu} (\bm x,t) = 0,
\quad  {\bm x} \in \partial \Omega,
\quad  0 < t \leq T ,
\end{equation}
\begin{equation}\label{5.3}
w_1(\bm x,0) = w_1^0(\bm x),
\quad  w_2(\bm x,0) = w_2^0(\bm x) ,
\quad  {\bm x} \in \Omega ,
\end{equation}
where $\nu$  is the normal to $\partial \Omega$.

We use a finite element approximation over the space \cite{brenner2008mathematical,ciarlet2002finite}. 
In $\Omega$, we perform triangulation while uniformly partitioning into $m$ segments on each individual variable $0 \leq x_1, x_2 \leq 1$.
The approximate solution $u_\alpha(\bm x, t), \ \alpha =1,2$ of the problem (\ref{5.1})--(\ref{5.3}) is based on piecewise linear finite elements: $u_\alpha(\bm x, t) \in V_h \subset H^1(\Omega), \ \alpha =1,2$.
We need to find $u_\alpha(\bm x, t) \in V_h, \ \alpha =1,2$ from the following conditions
\begin{equation}\label{5.4}
\begin{split}
 \left (\frac{\partial u_{1}}{\partial t}, v_1 \right ) & + a_{11} (u_1, v_1) + a_{12} (u_2, v_1) = 0 , \\
 \left (\frac{\partial u_{2}}{\partial t}, v_2 \right ) & + a_{21} (u_1, v_2) + a_{22} (u_2, v_2) = 0 ,
\end{split}
\end{equation} 
\begin{equation}\label{5.5}
 (u_1(\bm x,0), v_1) = (w_1^0(\bm x), v_1) ,
 \quad (u_2(\bm x,0), v_2) = (w_2^0(\bm x), v_2) ,
\end{equation} 
which take place for all $v_\alpha \in V_h, \ \alpha =1,2$.
In (\ref{5.4}), the bilinear forms are
\[
 a_{\alpha, \beta} (u,v) = \int_{\Omega}  d_{\alpha, \beta} (\bm x) \operatorname{grad} u(\bm x) \operatorname{grad} v(\bm x) \, d \bm x ,
 \quad \alpha, \beta =1,2.   
\] 

Let us define the operators $A_{\alpha \beta}: V_h \rightarrow V_h$ by the rule
\[
 a_{\alpha, \beta} (u,v) = (A_{\alpha, \beta} \, u,v) ,
 \quad \alpha, \beta =1,2.   
\]
We arrive from (\ref{5.4}), (\ref{5.5}) to a system of equations (\ref{2.1}) in which $p = 2, \ H_1 = H_2 = V_h$.
The initial conditions are given in the form (\ref{2.2}) when $u_\alpha^0(\bm x) = P w_\alpha^0(\bm x), \ \alpha =1,2$, and $P$ is the projection operator on the finite element space $V_h$. 
The property (\ref{2.5}) will be satisfied when 
\begin{equation}\label{5.6}
 d_{11}(\bm x) > 0, 
 \quad d_{22}(\bm x) > 0, 
 \quad 4 d_{11}(\bm x) \, d_{22}(\bm x)  > \big (d_{12}(\bm x) + d_{21}(\bm x) \big )^2 ,
 \quad  {\bm x} \in \Omega . 
\end{equation} 

We want to illustrate the possibilities of decoupling methods for solving the boundary value problem for the system of diffusion equations (\ref{5.1})--(\ref{5.3}) under rather general assumptions on diffusion coefficients.
For this purpose, we consider the test problem in which
\[
 d_{11}(\bm x), \ d_{12}(\bm x), \ d_{22}(\bm x) =  
\left \{\begin{array}{rrrrr}
 5, & -2, & 1, & x_2 > 0.5 ,\\
 1, & 2, & 5, & x_2 \leq 0.5,\\
\end{array}
\right . 
 \quad {\bm x} \in \Omega .  
\] 
The initial conditions are taken as
\[
 w_1^0(\bm x) = \frac{1}{2} + \frac{1}{4}  \big (\cos(2 \pi x_1) +\cos(2 \pi x_2) \big ),
 \quad w_2^0(\bm x) = 16 x_1^2 (1-x_1)^2 .
\] 
Thus $0 \leq w_\alpha^0(\bm x) \leq 1 , \ \alpha =1,2$, and the stationary solution is $w_1(\bm x, \infty) = 1/2, \ w_2(\bm x, \infty) = 8/15$.

\begin{figure}[htbp]
\centering
\includegraphics[width=0.4\linewidth]{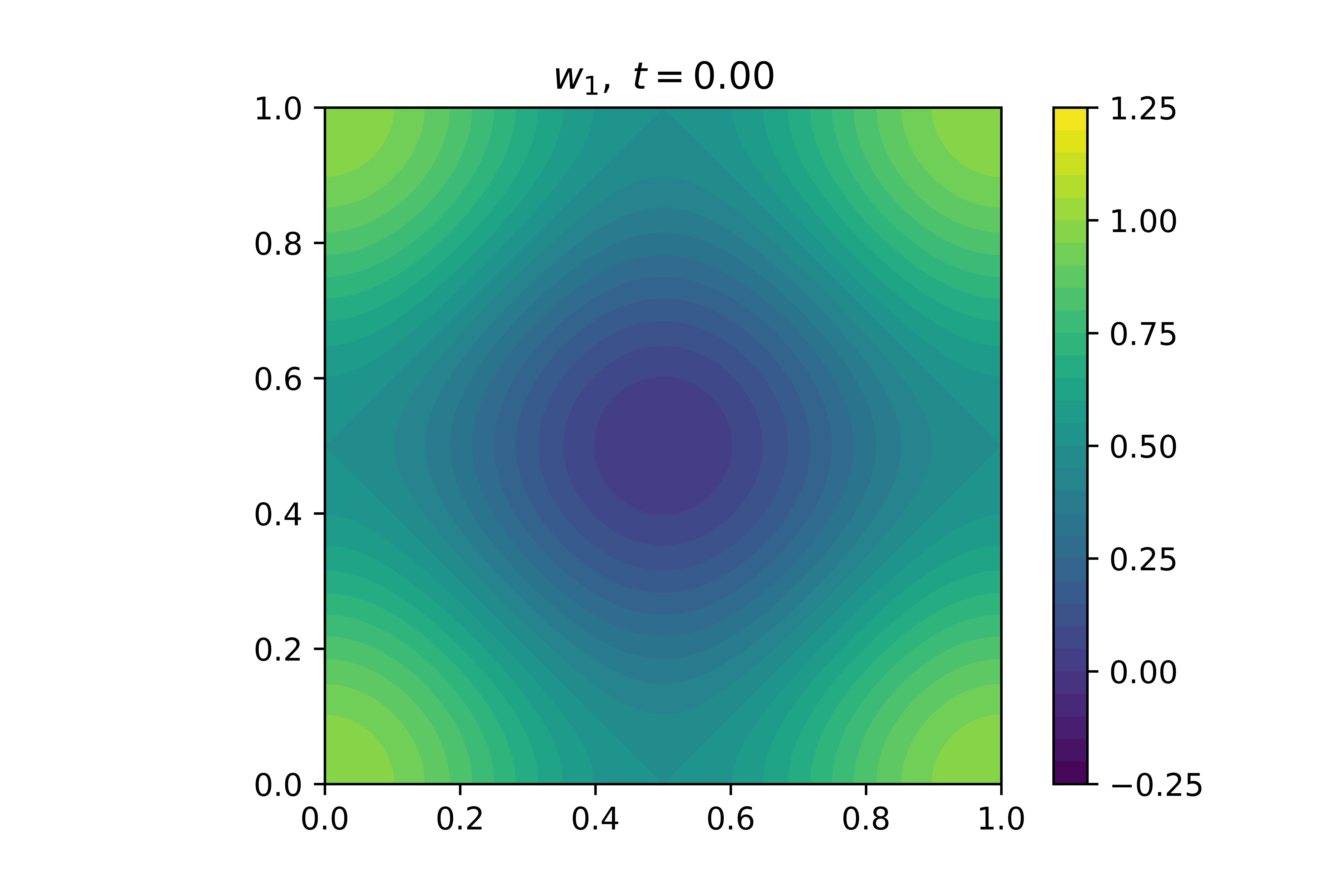} \includegraphics[width=0.4\linewidth]{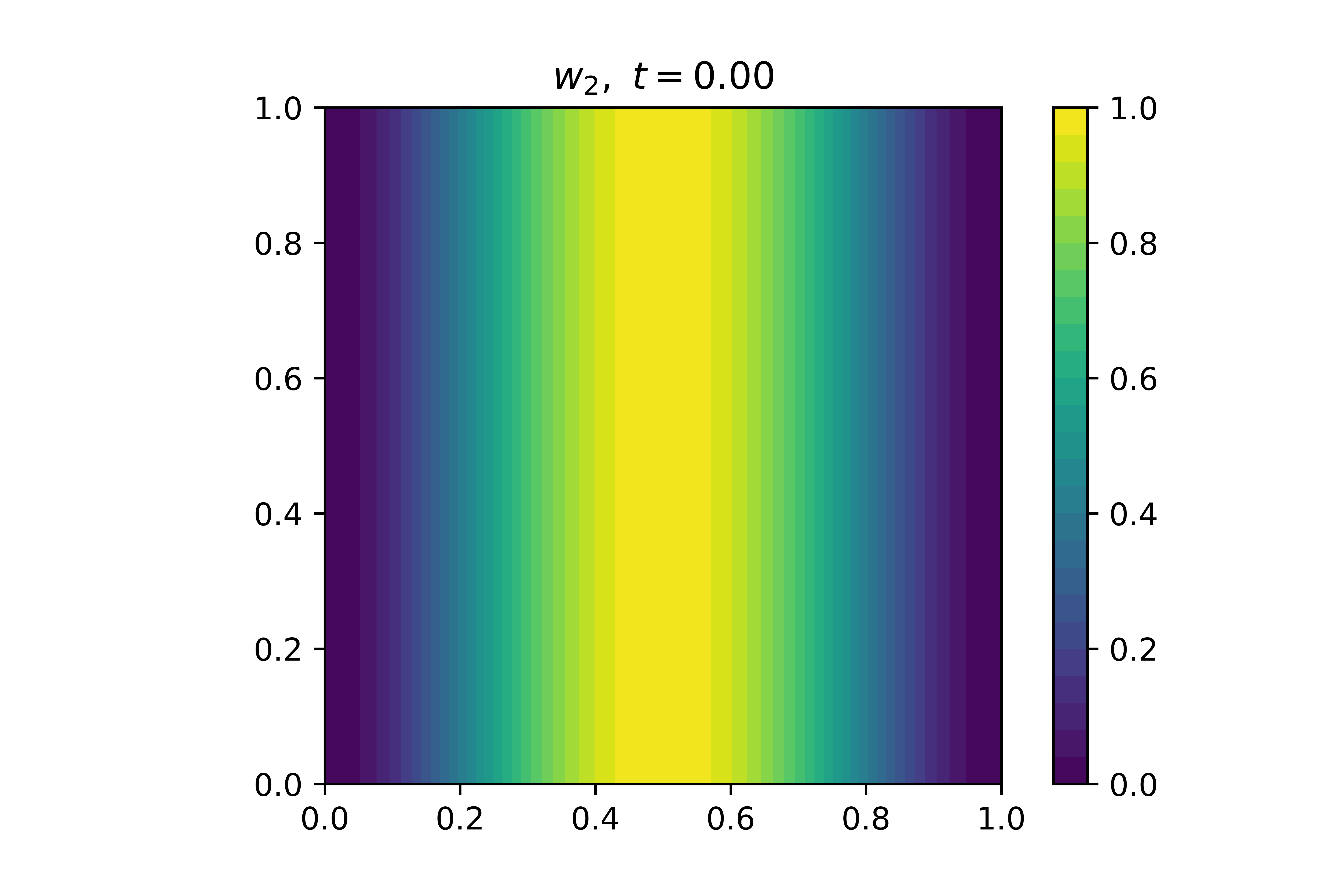} \\
\includegraphics[width=0.4\linewidth]{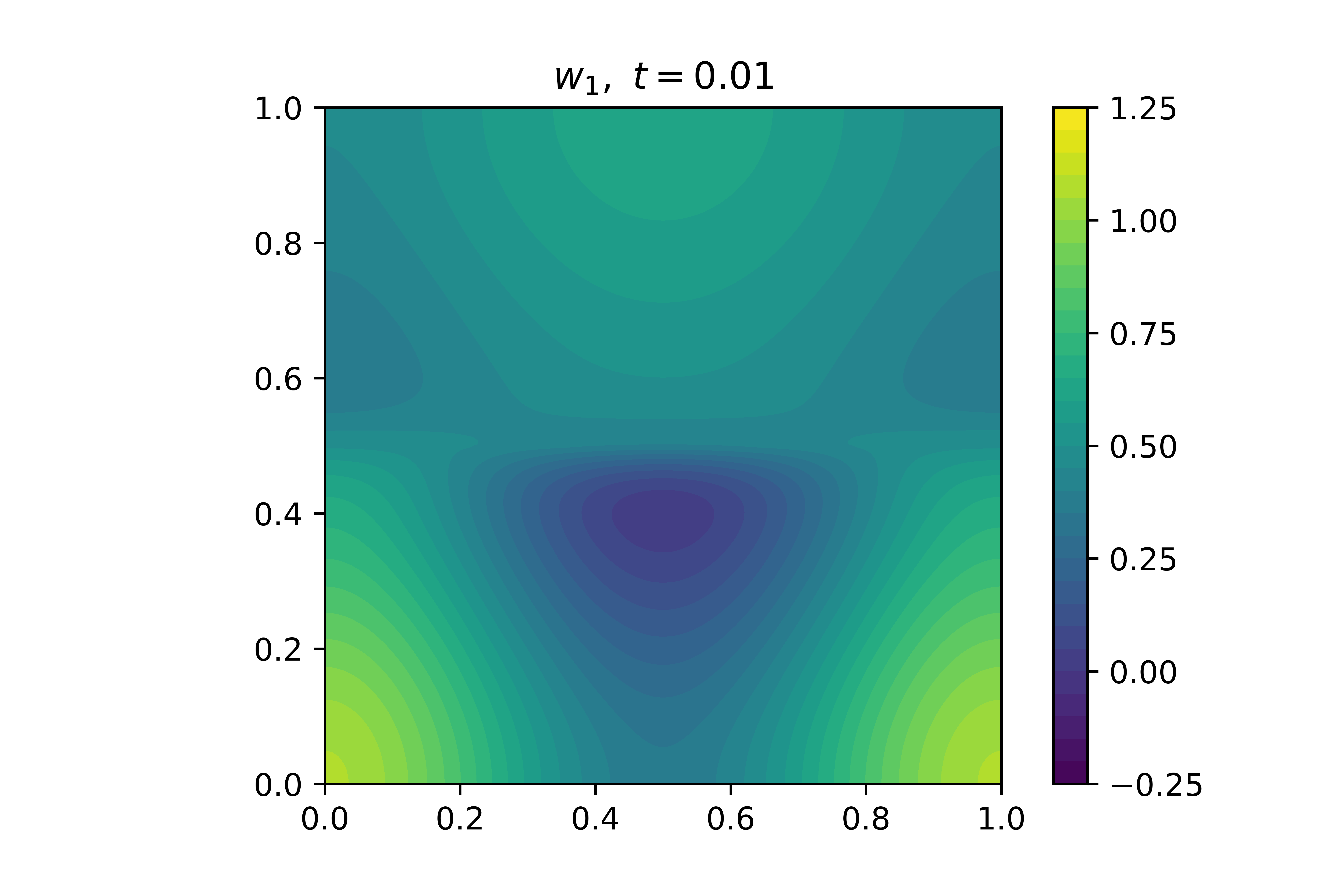} \includegraphics[width=0.4\linewidth]{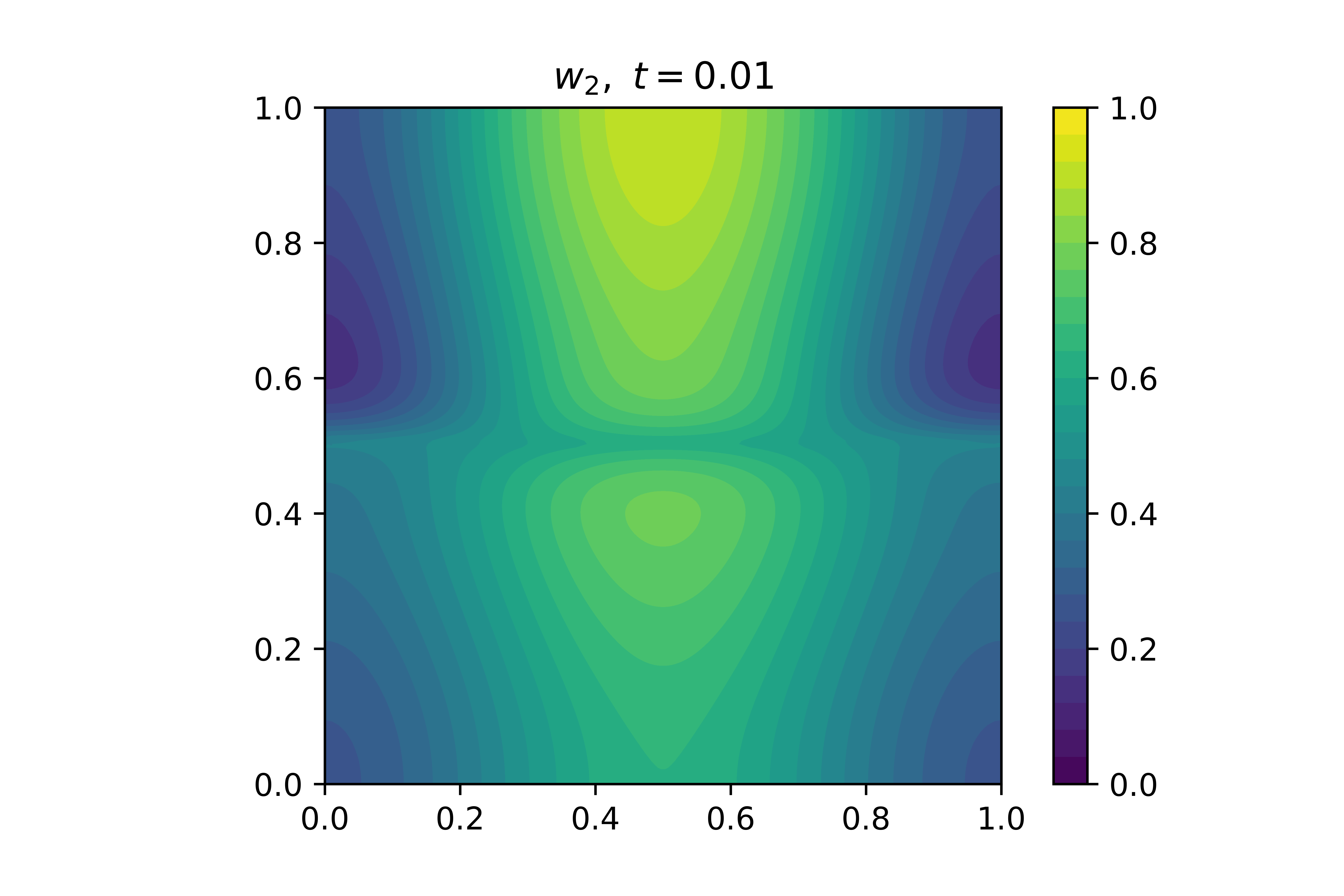} \\
\includegraphics[width=0.4\linewidth]{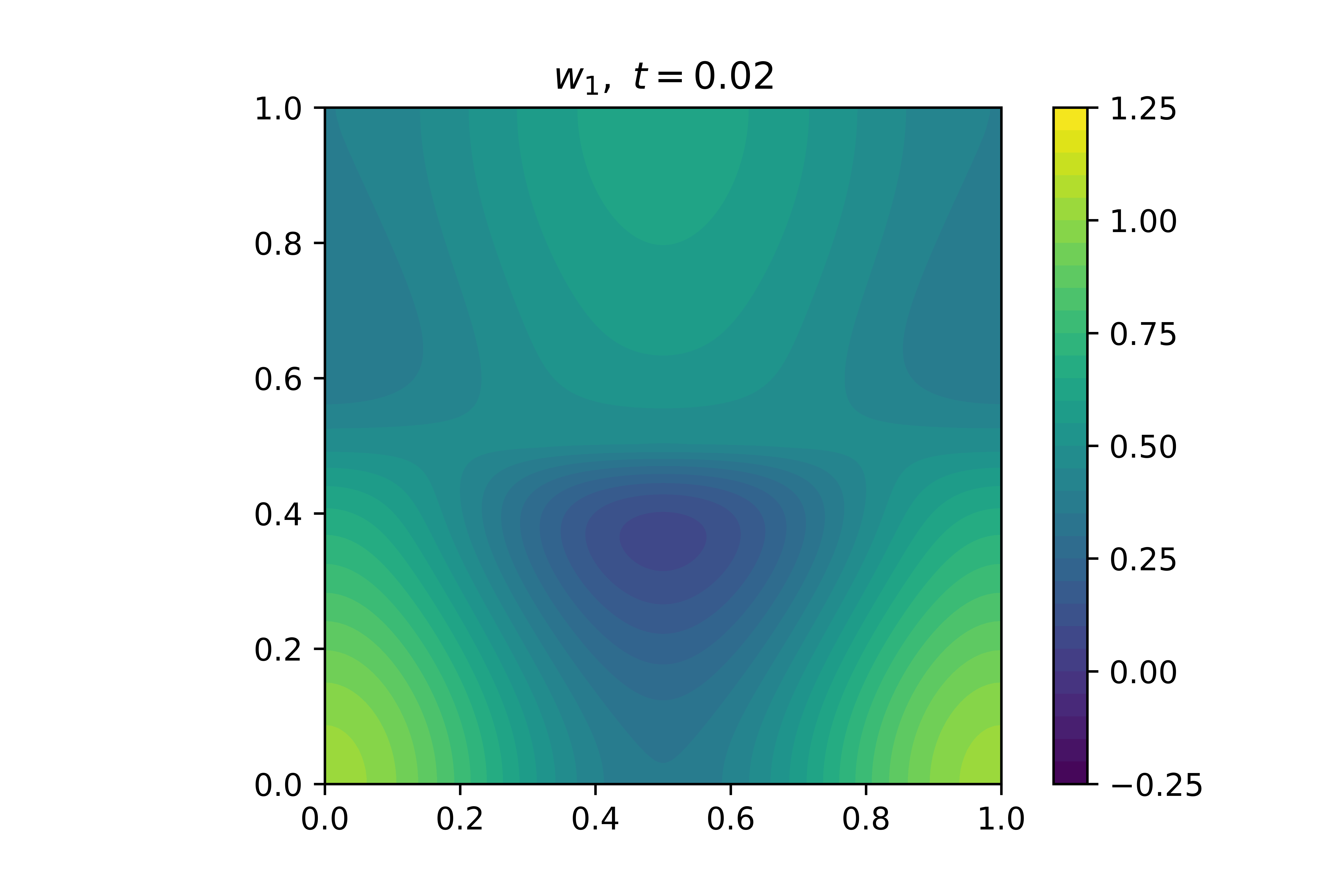} \includegraphics[width=0.4\linewidth]{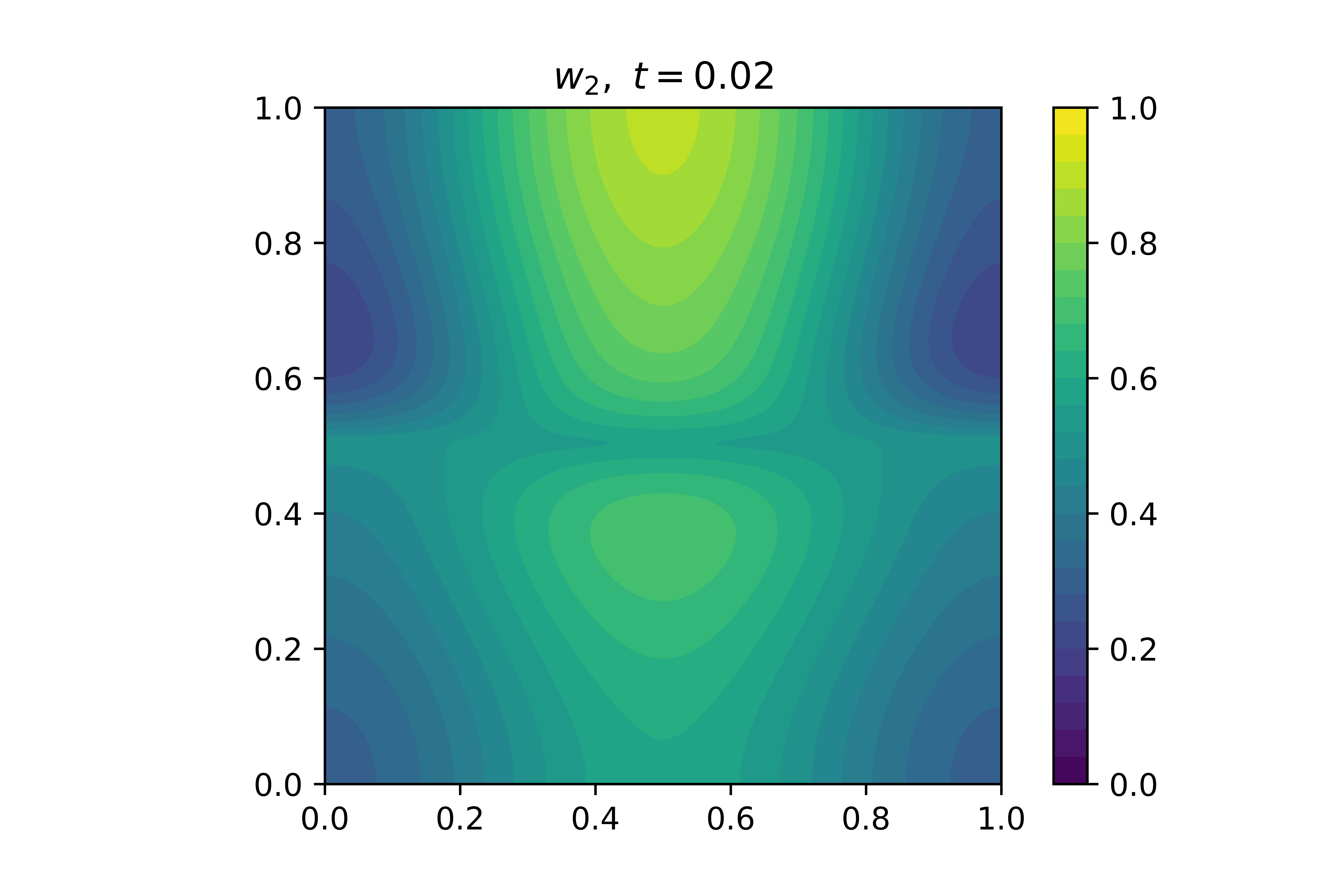} \\
\includegraphics[width=0.4\linewidth]{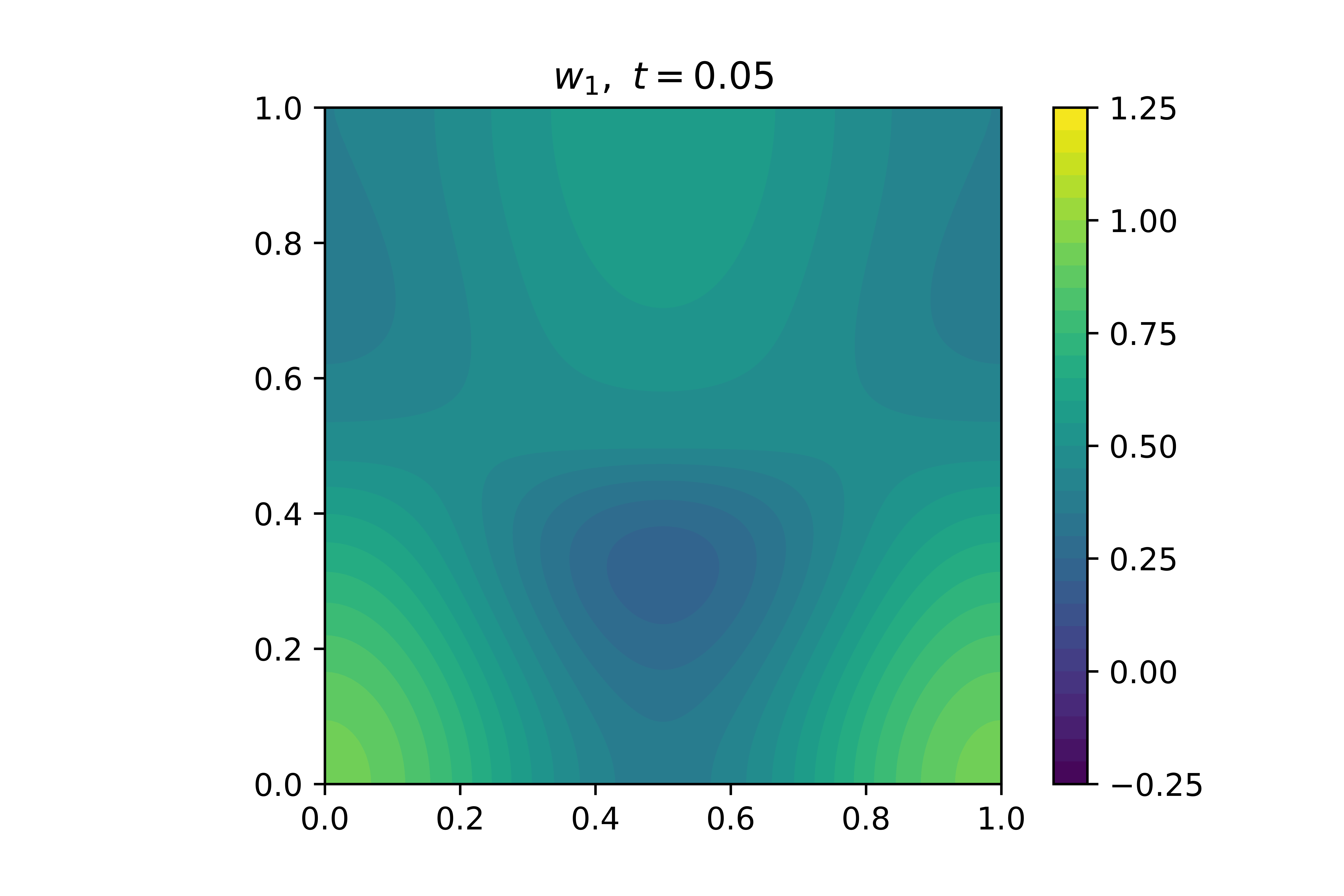} \includegraphics[width=0.4\linewidth]{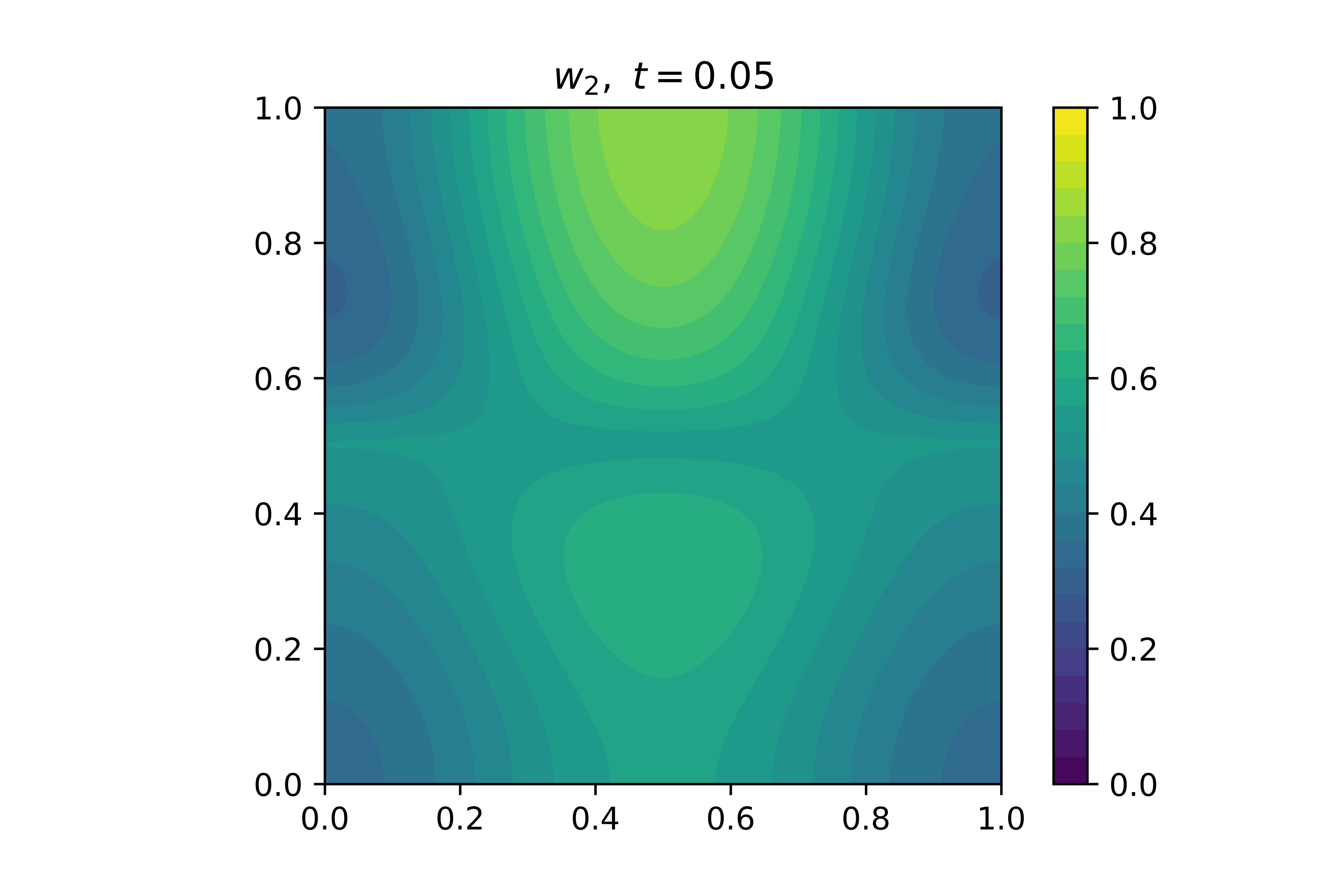} \\
\includegraphics[width=0.4\linewidth]{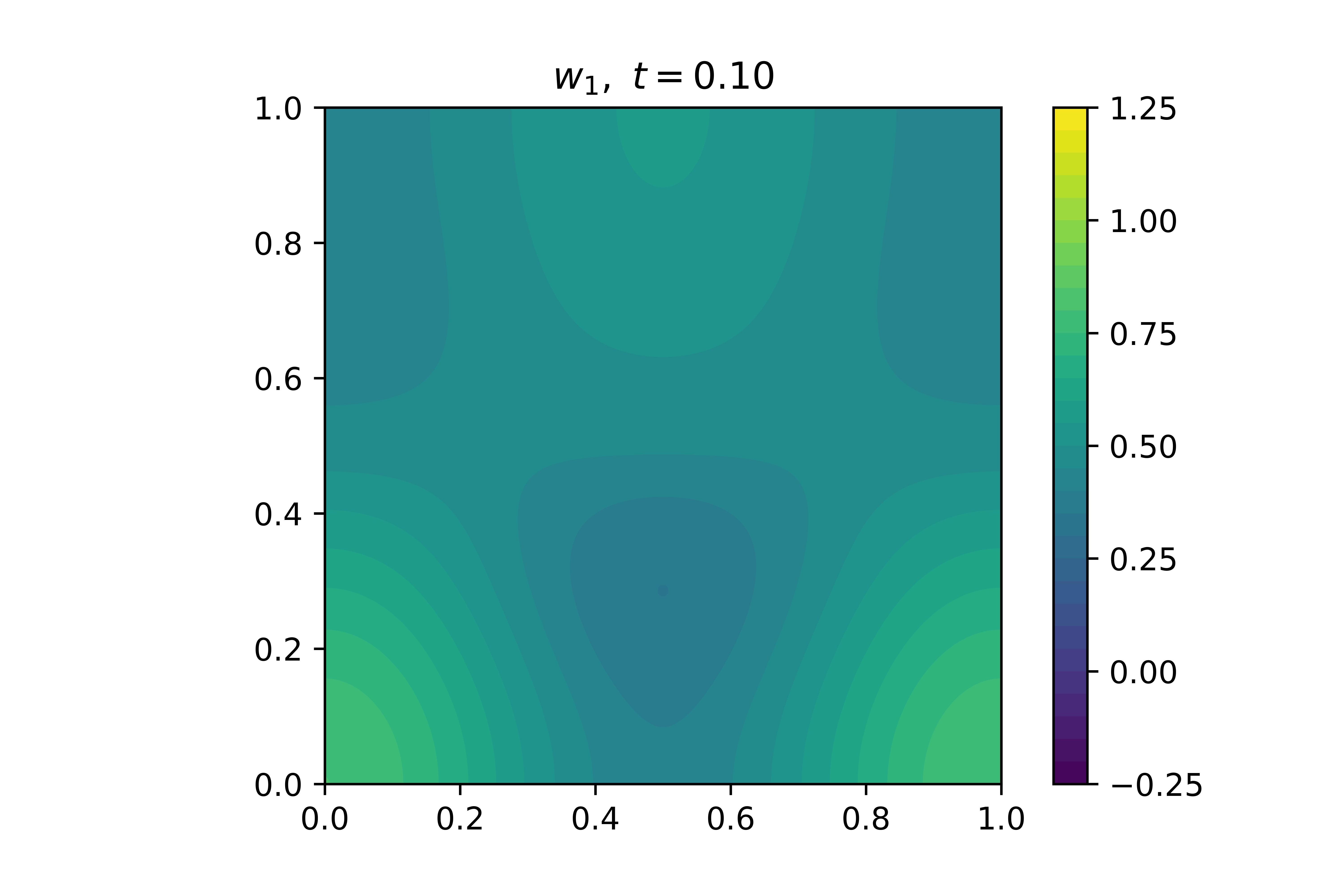} \includegraphics[width=0.4\linewidth]{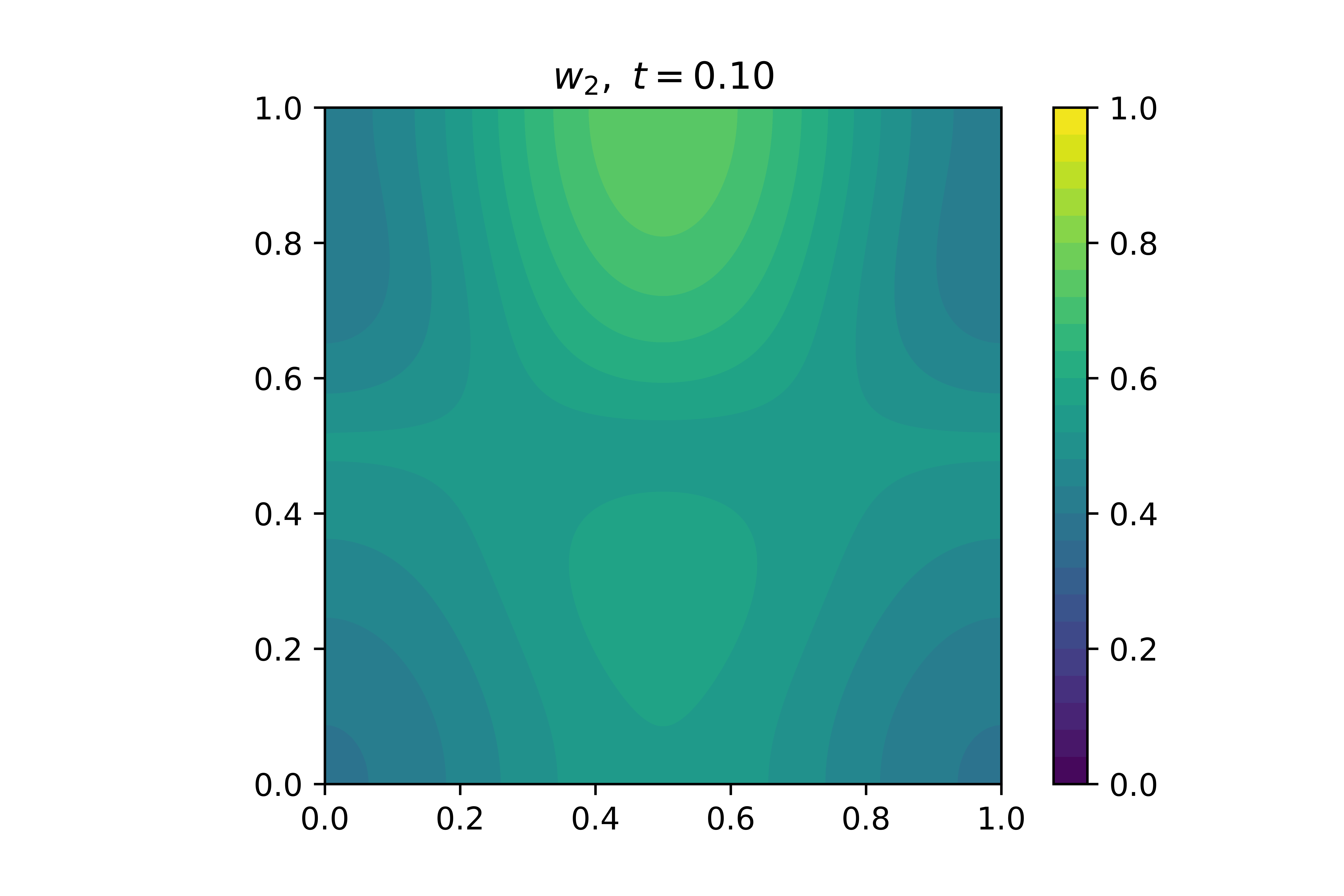} 
\caption{Solution of the initial boundary value problem for different time instants.}
\label{f-1}
\end{figure}

The approximate solution for individual time moments is shown in Fig.~\ref{f-1}.
The space approximation was performed at $m=100$, time step $\tau = 10^{-3} \ (N=100)$.
We used a weighted scheme (see (\ref{2.7}))
\begin{equation}\label{5.7}
\begin{split}
 \frac{\overline{y}_1^{n+1} - \overline{y}_1^{n}}{\tau } + A_{11} (\sigma \overline{y}_1^{n+1} + (1-\sigma)\overline{y}_1^{n}) 
 + A_{12} (\sigma \overline{y}_2^{n+1} + (1-\sigma)\overline{y}_2^{n} ) = 0, \\
 \frac{\overline{y}_2^{n+1} - \overline{y}_2^{n}}{\tau } + A_{21} (\sigma \overline{y}_1^{n+1} + (1-\sigma)\overline{y}_1^{n} ) 
 + A_{22} (\sigma \overline{y}_2^{n+1} + (1-\sigma)\overline{y}_2^{n} ) = 0 \\
\end{split}
\end{equation} 
for $\sigma = 1/2$.  

A significant restructuring of the initial state is observed with the stationary transition.
Large gradients at the boundary of subdomains $x_2 = 0.5$ ($\Omega = \Omega_1 \cup \Omega_2$) are associated with discontinuities of diffusion and cross-diffusion coefficients.

In our test problem, the exact solution is unknown.
Therefore, we evaluate the accuracy of the decomposition-composition schemes by comparing the obtained solutions with the approximate solution, which is obtained using standard weighted schemes (\ref{5.7}).
In this case, we evaluate the error due to the splitting of the problem rather than the overall error due to time and space approximations.
We will compare first-order precision splitting schemes to a purely implicit scheme ($\sigma = 1$) and second-order precision schemes --- to a symmetric scheme ($\sigma = 1/2$).
Let $y_\alpha^n(\bm x), \ \alpha =1,2$ be the components of the approximate solution of the splitting scheme, then let us determine the error rate at a separate time step in two norms:
\[
\varepsilon_\alpha (t^n) = \|y_\alpha^n(\bm x) - \overline{y}_\alpha^n(\bm x)\|,
\quad \delta_\alpha (t^n) = \max_{\bm x \in \Omega} |y_\alpha^n(\bm x) - \overline{y}_\alpha^n(\bm x)| ,
\quad n = 0, \ldots, N ,
\quad \alpha = 1,2 .
\]

The error introduced by decomposition-composition can be analyzed by extracting additional terms in the equation.
For first-order accuracy schemes, the reference solution is determined from Eq.
\[
\frac{\overline{\bm y}^{n+1} - \overline{\bm y}^{n}}{\tau } +
{\bm A} \overline{\bm y}^{n+1} = 0 .
\]
For $\bm z^n = \overline{\bm y}^n - \bm y^n $ we have Eq.
\[
\frac{{\bm z}^{n+1} - {\bm z}^{n}}{\tau } +
{\bm A} \bm z^{n+1}  = \bm \psi^n .
\]

For the explicit-implicit scheme (\ref{3.2}) with the diagonal part of the operator matrix $\bm A$ taken to a new level in time, the error of the splitting scheme is estimated by the value of
\[
\bm \psi^n = - \tau (\bm A - \sigma \bm D) \frac{\bm y^{n+1} - \bm y^{n}}{\tau } .
\]
In this case, for the individual components of the solution of the system of equations (\ref{5.1}), we obtain
\begin{equation}\label{5.8}
\begin{split}
(1+ \sigma \tau A_{11}) \frac{y_1^{n+1} - y_1^{n}}{\tau } + A_{11} y_1^{n} + A_{12} y_2^{n} = 0, \\
(1+ \sigma \tau A_{22}) \frac{y_2^{n+1} - y_2^{n}}{\tau } + A_{21} y_1^{n} + A_{22} y_2^{n} = 0 .
\end{split}
\end{equation}

\begin{figure}[htbp]
\centering
\includegraphics[width=0.4\linewidth]{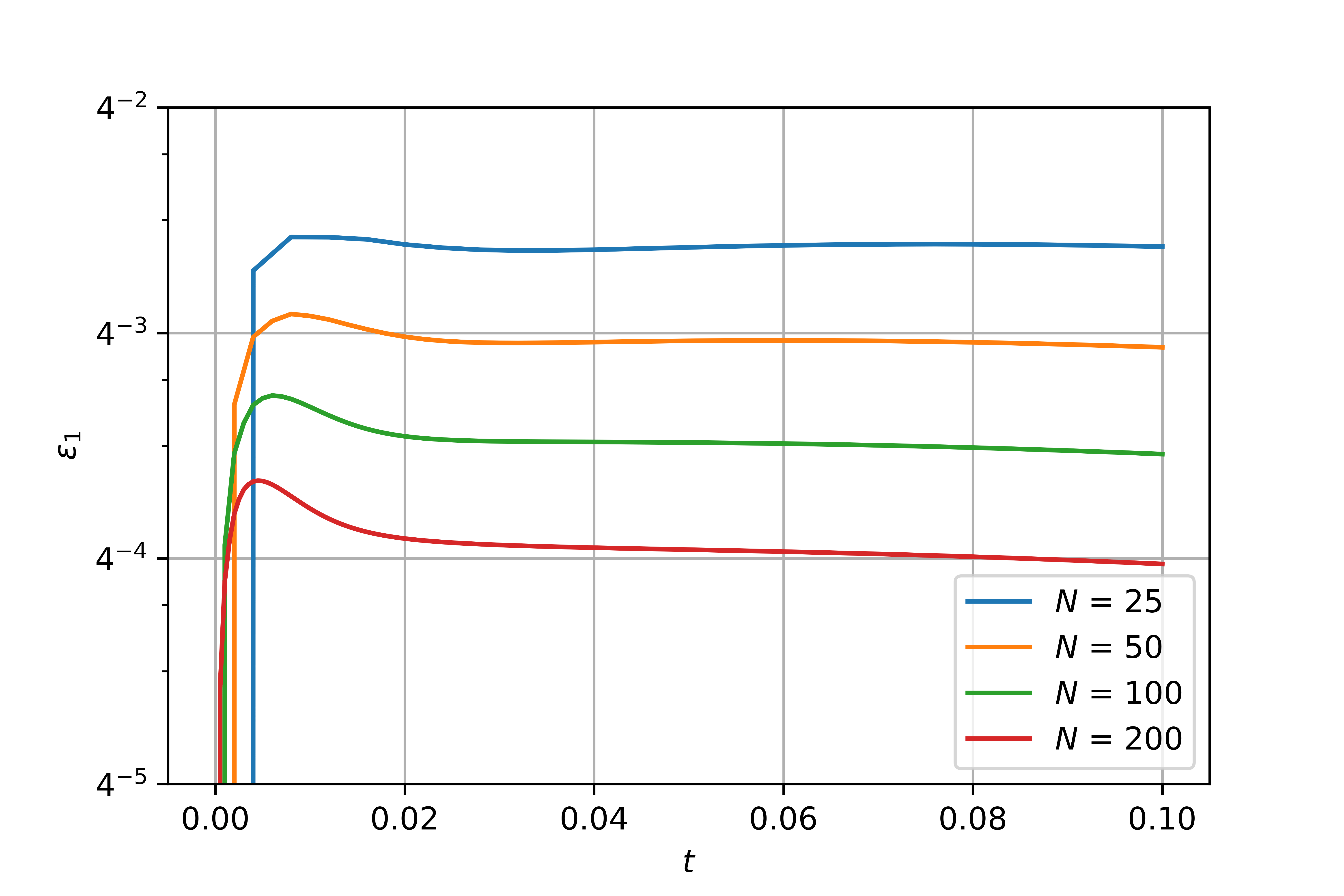} \includegraphics[width=0.4\linewidth]{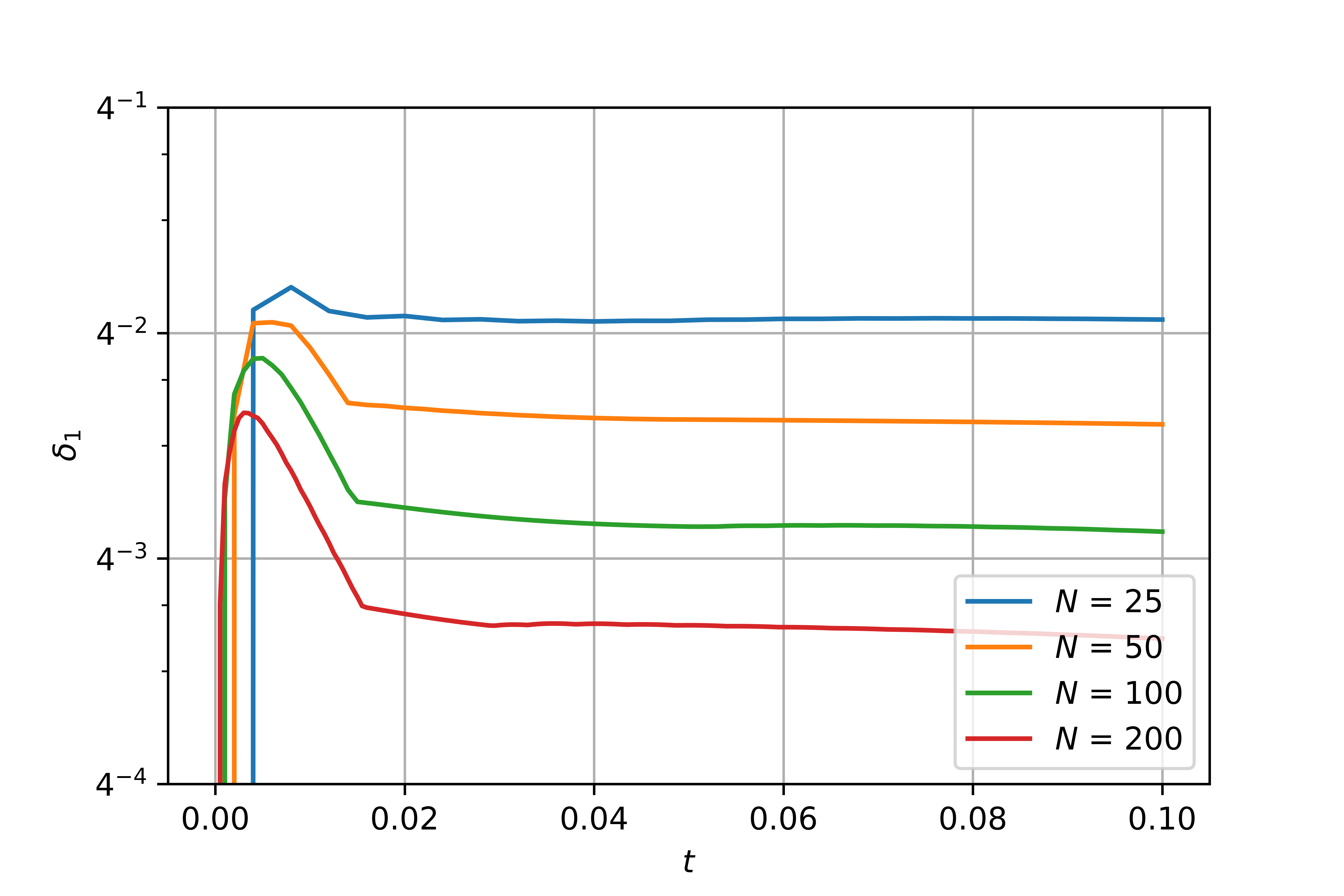} \\
\includegraphics[width=0.4\linewidth]{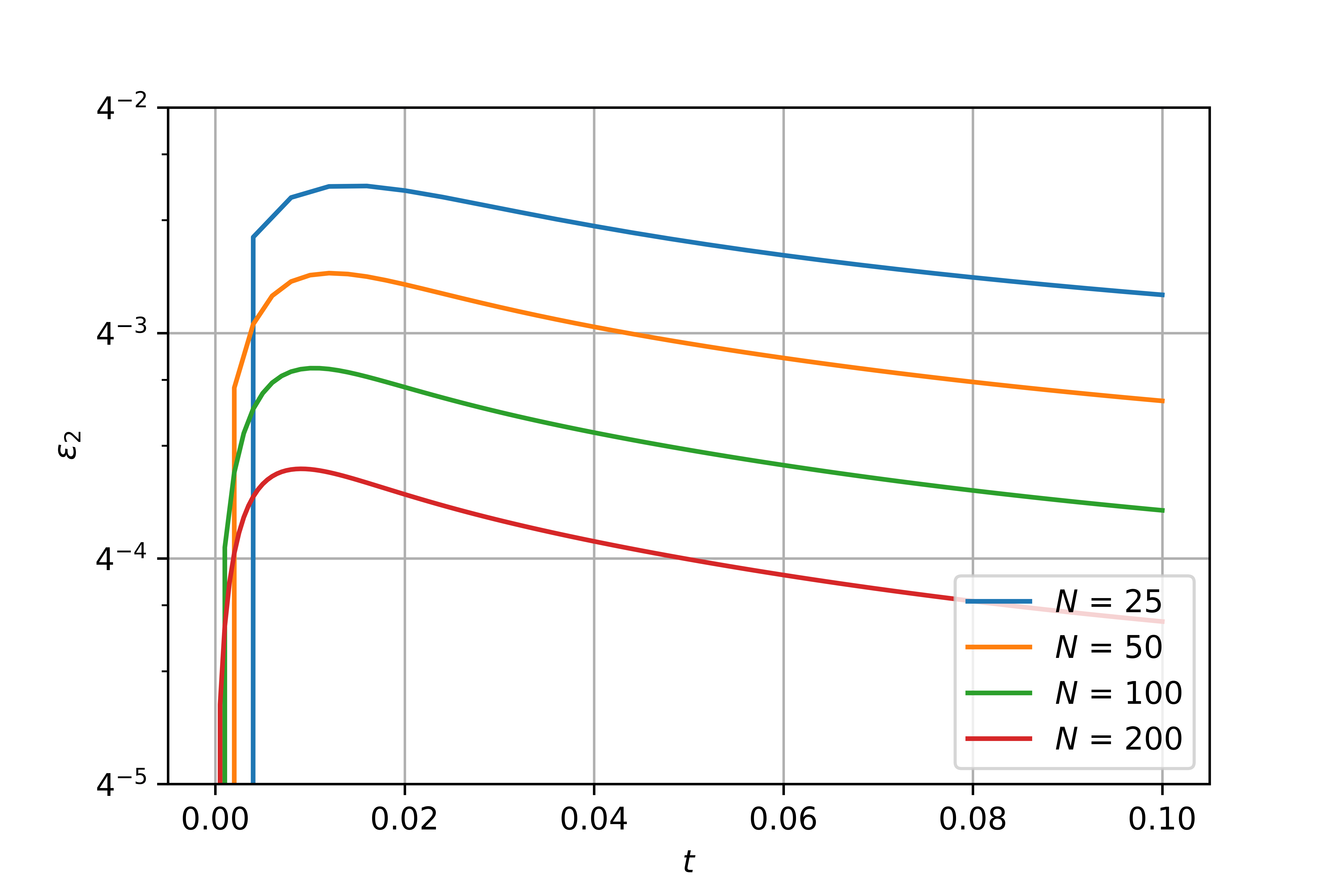} \includegraphics[width=0.4\linewidth]{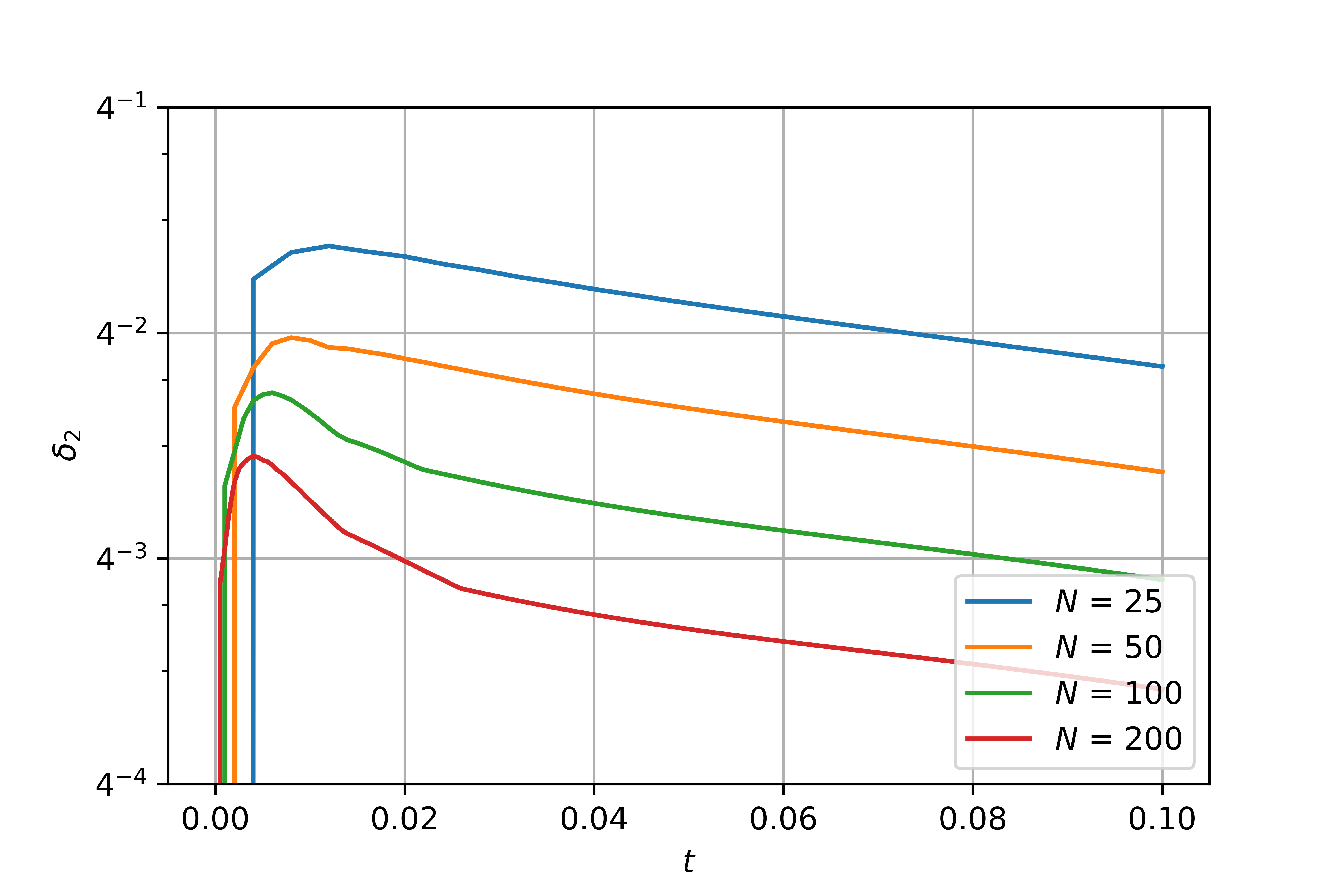}
\caption{Decomposition-composition scheme (\ref{5.8}).}
\label{f-2}
\end{figure}

The results of comparing the decomposition-composition scheme (\ref{5.8}), for which $\gamma =2$ in (\ref{3.1}), and the standard weighted scheme (\ref{5.7}) at $\sigma = 1$ are shown in Fig.~\ref{f-2}. The closeness of the approximate solutions is estimated by the value $\mathcal{O}(\tau)$.

The decomposition-composition scheme (\ref{3.8}) for the system of equations (\ref{5.1}) takes the following form
\begin{equation}\label{5.9}
\begin{split}
\frac{y_1^{n+1} - y_1^{n}}{\tau } + A_{11} y_1^{n} + A_{12} y_2^{n} = 0, \\
\frac{y_2^{n+1} - y_2^{n}}{\tau } + A_{21} y_1^{n+1} + A_{22} y_2^{n+1} = 0 .
\end{split}
\end{equation}
The value 
\[
\bm \psi^n = \tau \bm L^* \frac{\bm y^{n+1} - \bm y^{n}}{\tau } 
\]
estimates the error of this splitting scheme.

\begin{figure}[htbp]
\centering
\includegraphics[width=0.4\linewidth]{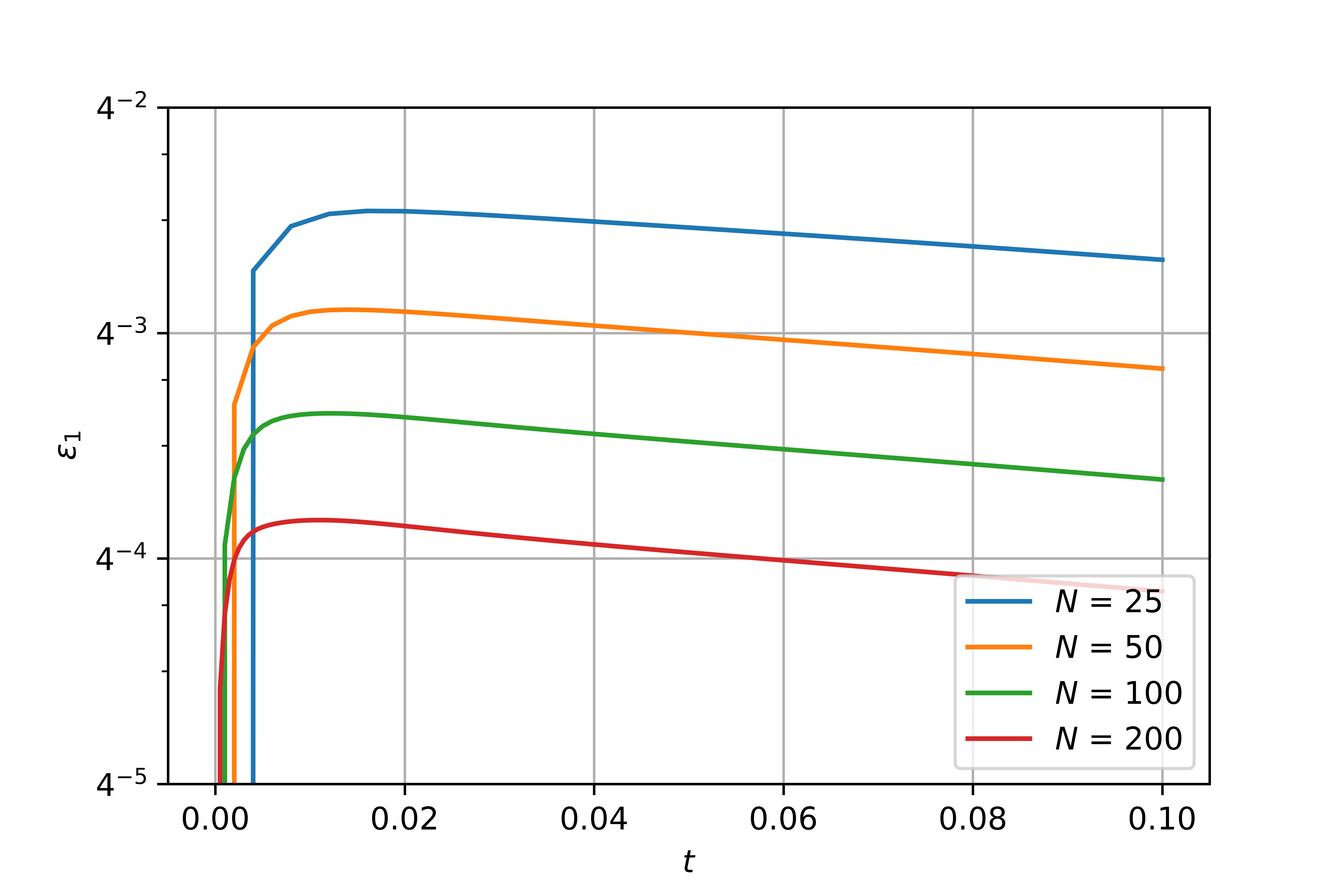} \includegraphics[width=0.4\linewidth]{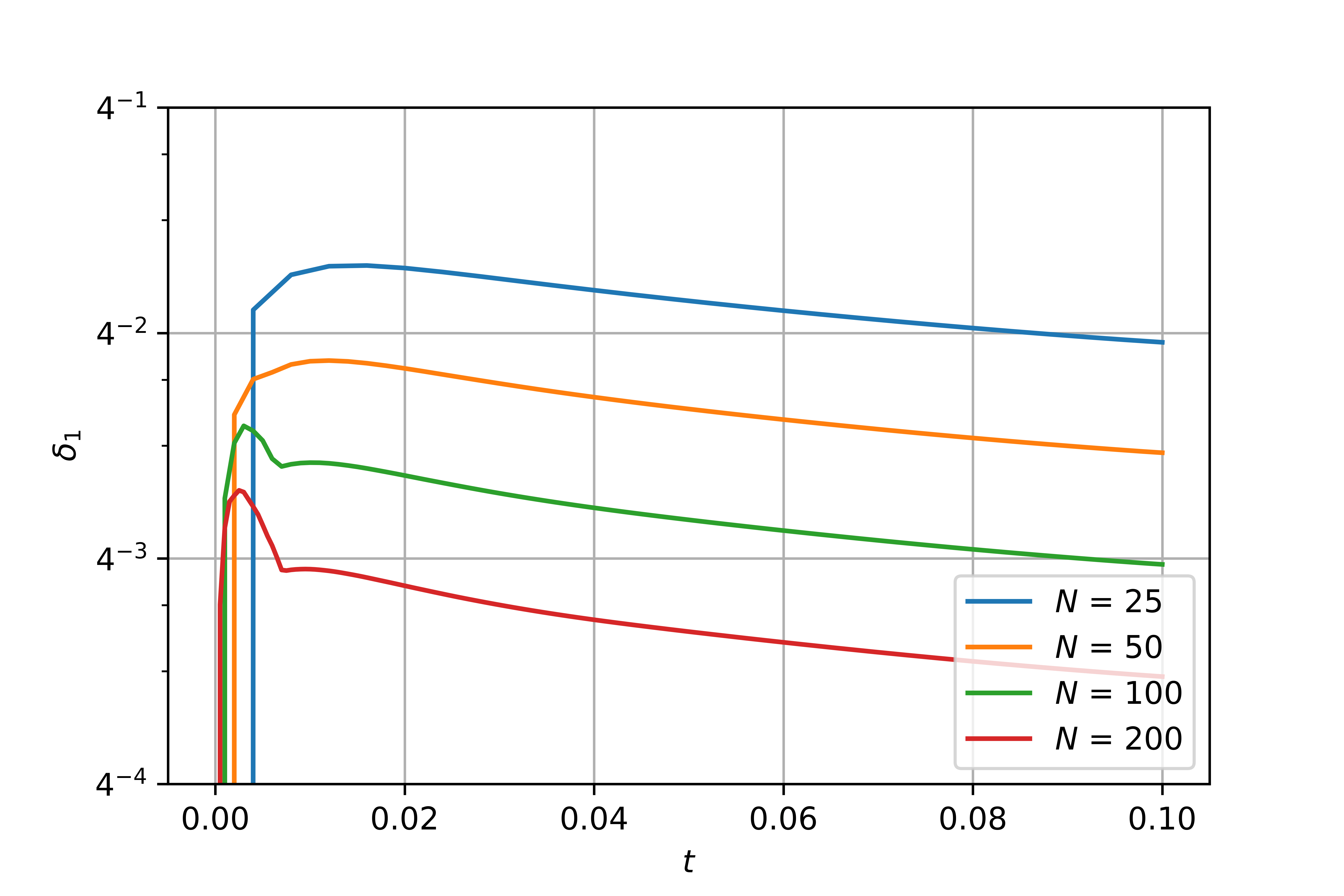} \\
\includegraphics[width=0.4\linewidth]{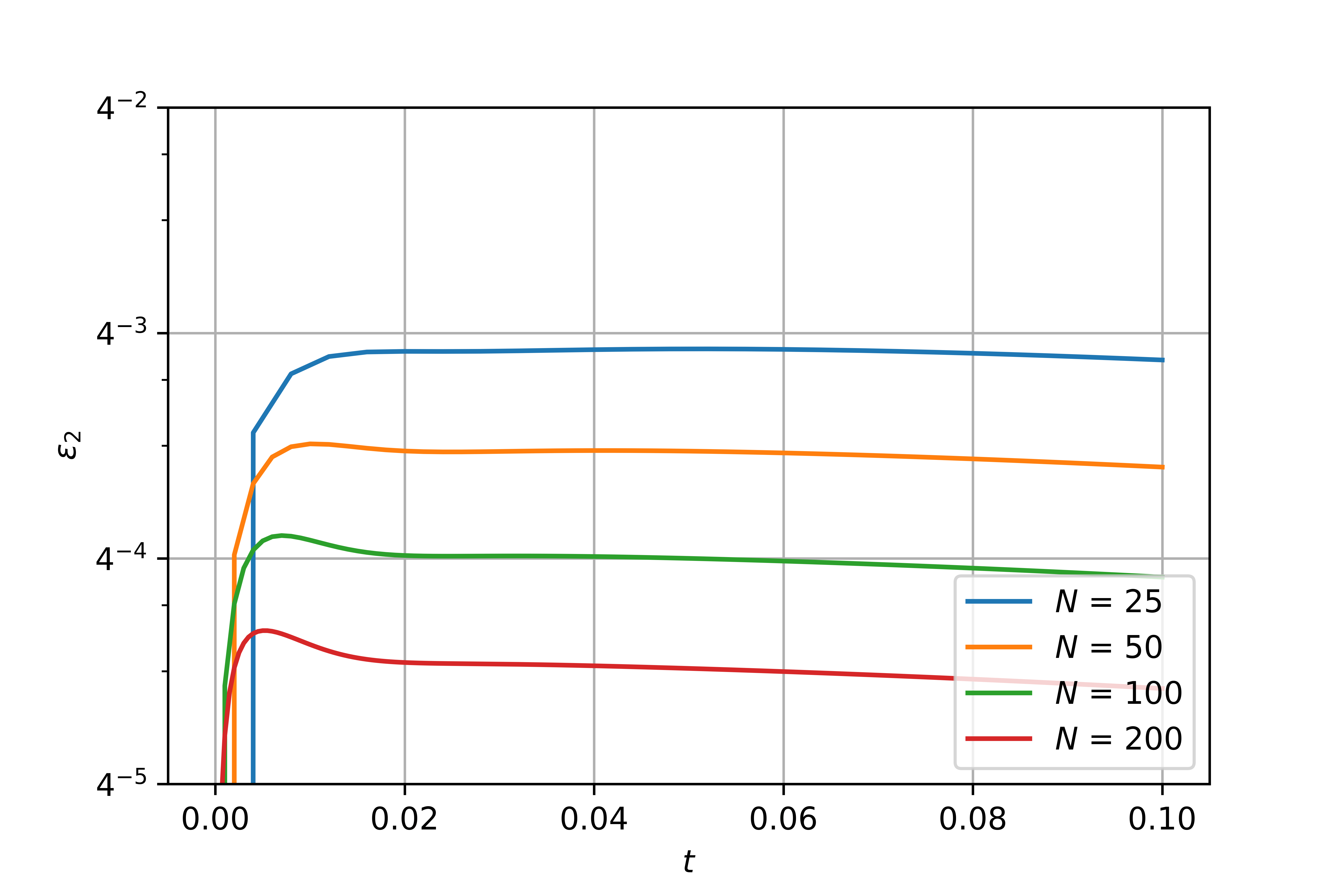} \includegraphics[width=0.4\linewidth]{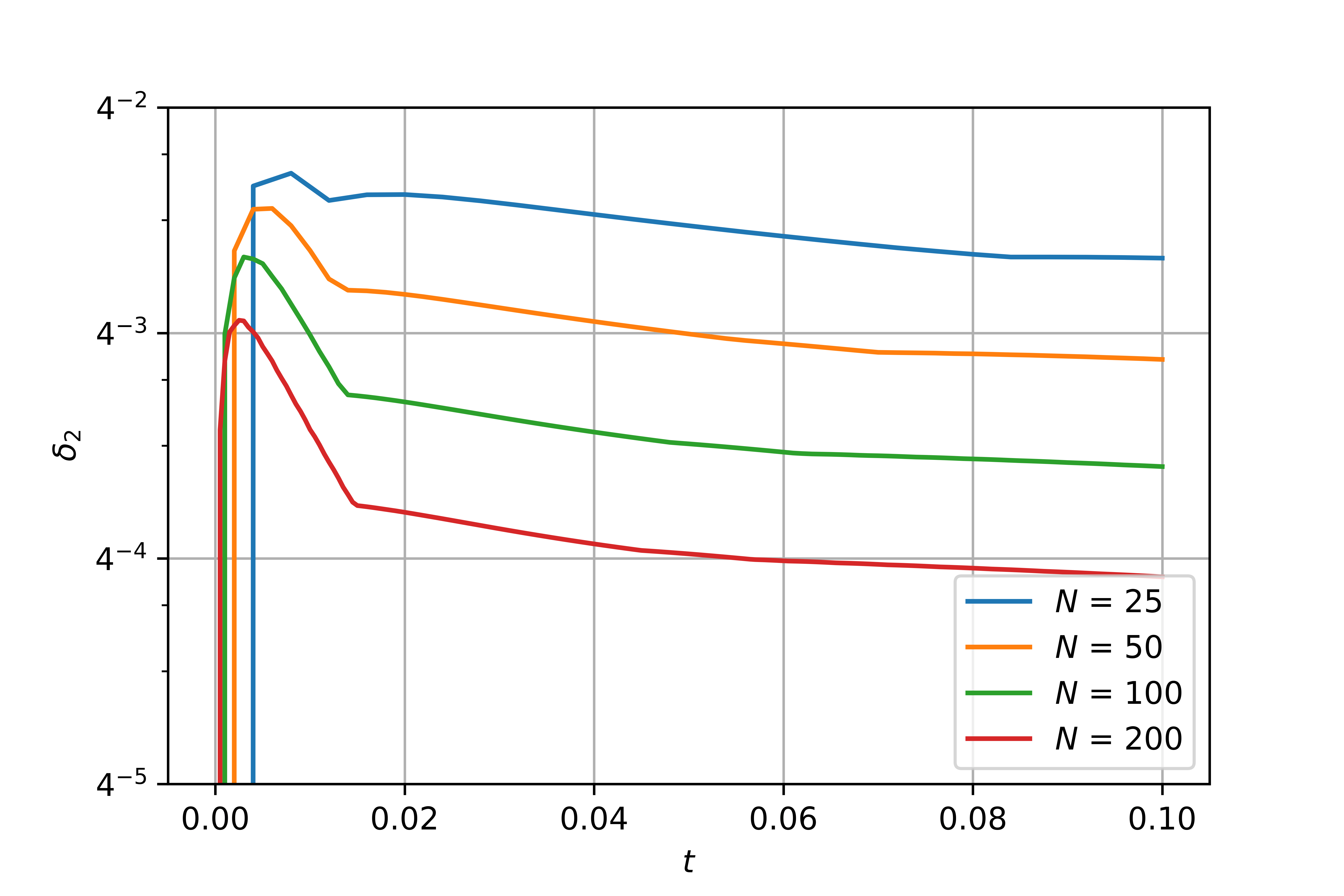}
\caption{Decomposition-composition scheme (\ref{5.9}).}
\label{f-3}
\end{figure}

The numerical results for this scheme and the purely implicit scheme are shown in Fig.~\ref{f-3}.
We can conclude that there is no significant accuracy improvement in the scheme (\ref{5.9}) compared to the scheme (\ref{5.8}).
Under these conditions, we can give preference to the scheme (\ref{5.8}), taking into account the possibility of parallel organization of computations for calculating individual components of the approximate solution.

We consider the alternating-triangle method schemes (\ref{3.3}), (\ref{3.9}), (\ref{3.11}) in two basic variants.
At $\sigma = 1$ we have the operator analog of the Douglas-Rachford scheme \cite{DouglasRachford1956} when
\[
\begin{split}
\frac{{\bm y}^{n+1/2} - {\bm y}^{n}}{\tau } + \bm A_1 \bm y^{n+1/2} + \bm A_2 \bm y^{n} = 0 , \\
\frac{{\bm y}^{n+1} - {\bm y}^{n+1/2}}{\tau } + \bm A_2 (\bm y^{n+1} - y^{n}) = 0 .
\end{split}
\]
For the splitting error, we have
\[
\bm \psi^n = - \tau^2 \bm A_1 \bm A_2 \frac{\bm y^{n+1} - \bm y^{n}}{\tau } .
\]

\begin{figure}[htbp]
\centering
\includegraphics[width=0.4\linewidth]{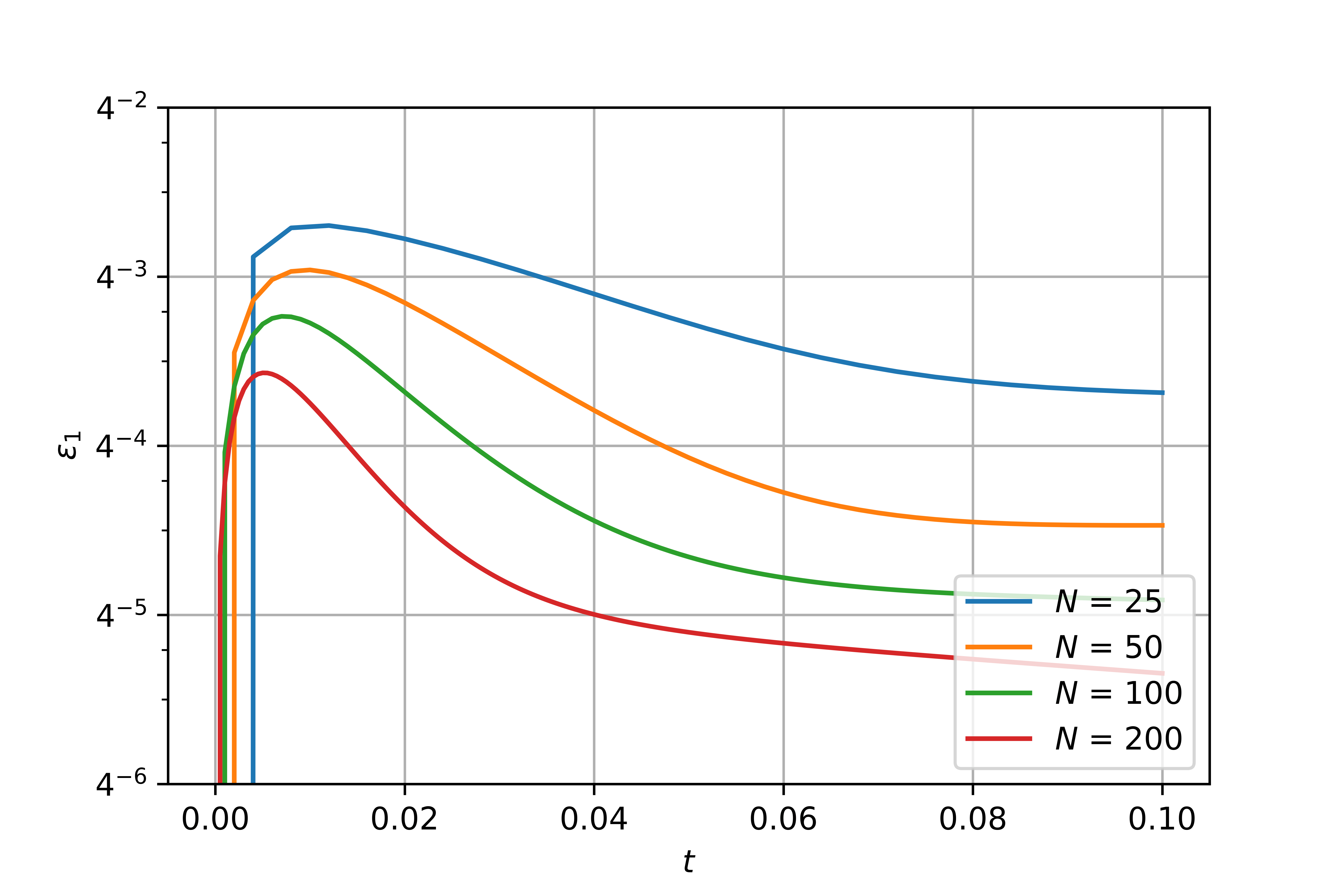} \includegraphics[width=0.4\linewidth]{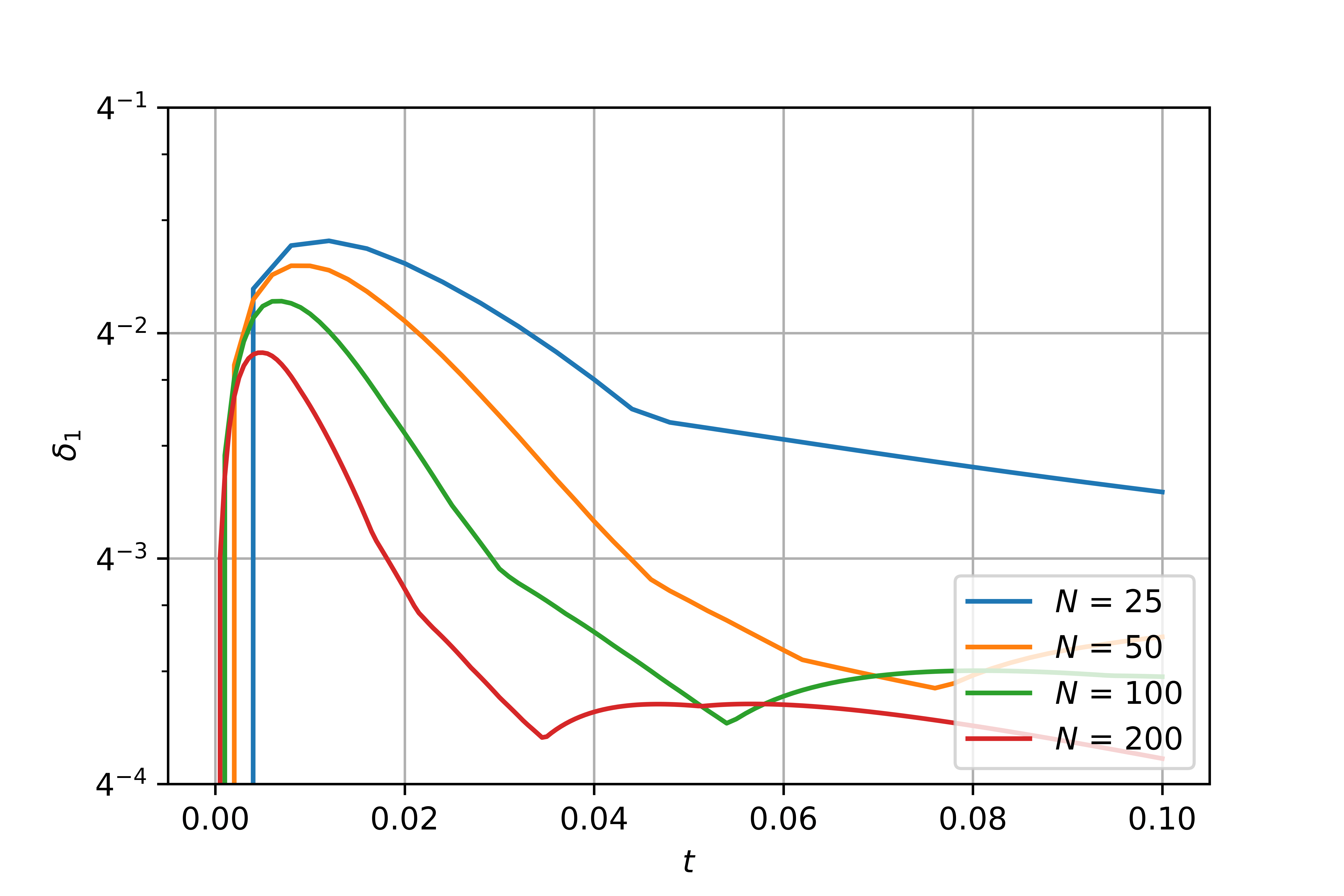} \\
\includegraphics[width=0.4\linewidth]{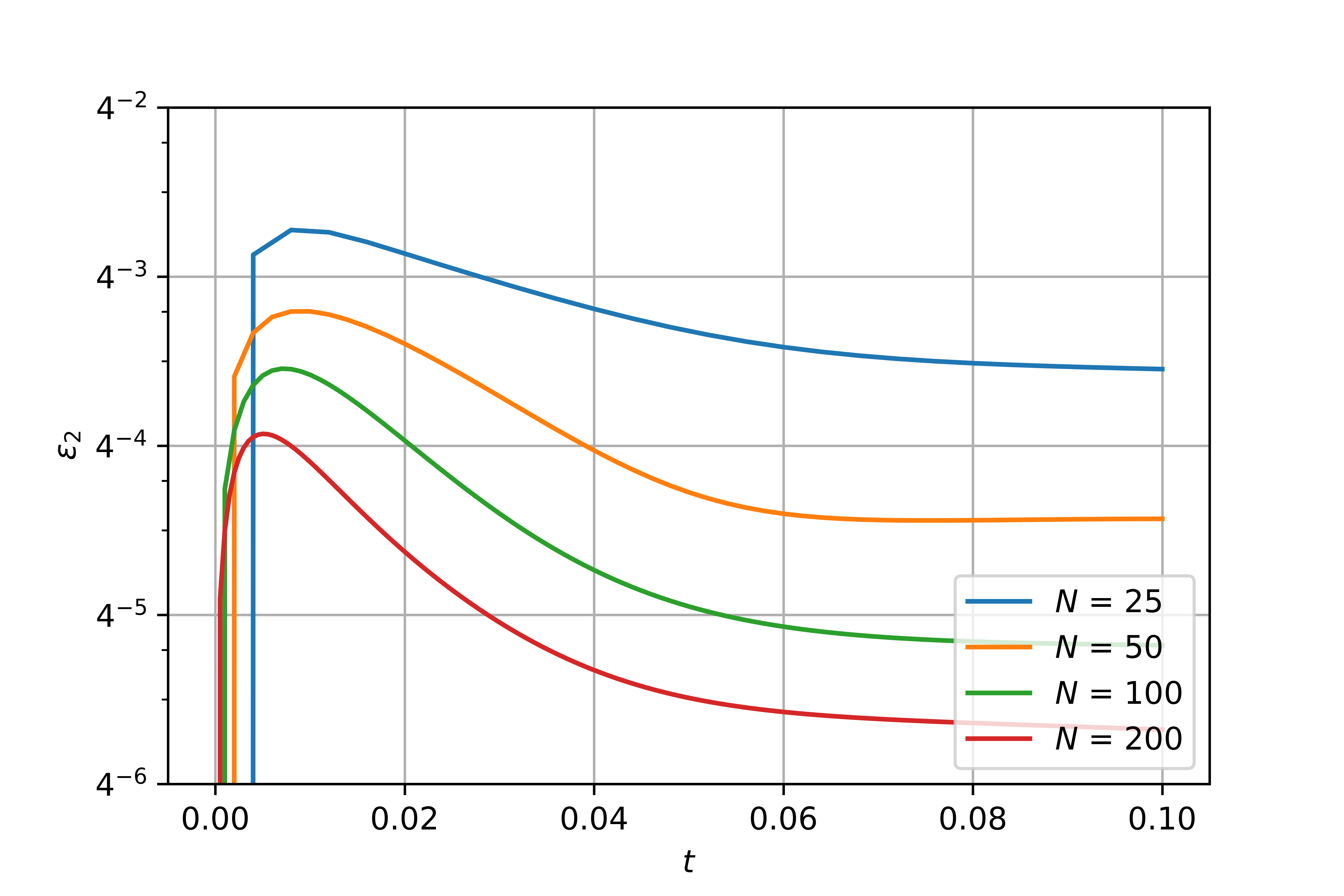} \includegraphics[width=0.4\linewidth]{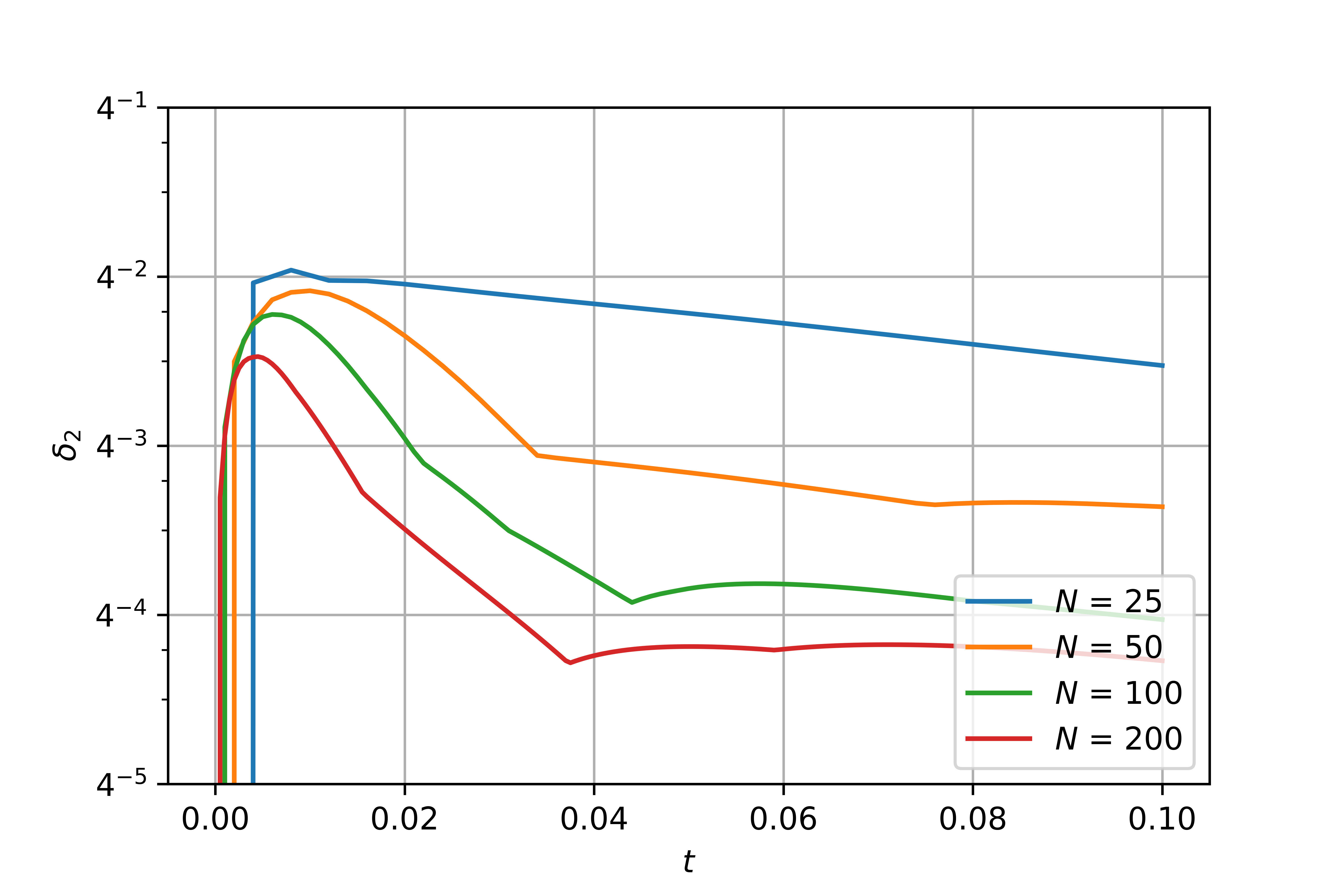}
\caption{Decomposition-composition scheme (\ref{5.10}).}
\label{f-4}
\end{figure}

For the test problem under consideration, this leads to the following scheme
\begin{equation}\label{5.10}
\begin{split}
& \left (1+ \frac{1}{2} \tau A_{11} \right ) \frac{y_1^{n+1/2} - y_1^{n}}{\tau } + A_{11} y_1^{n} + A_{12} y_2^{n} = 0, \\
& \left (1+ \frac{1}{2} \tau A_{22} \right ) \frac{y_2^{n+1/2} - y_2^{n}}{\tau } + A_{21} y_1^{n+1/2} + A_{22} y_2^{n} = 0 , \\
& \frac{y_2^{n+1} - y_2^{n+1/2}}{\tau } + \frac{1}{2} A_{22} (y_2^{n+1} - y_2^{n} ) = 0, \\
& \frac{y_1^{n+1} - y_1^{n+1/2}}{\tau } + \frac{1}{2} A_{11} (y_1^{n+1} - y_1^{n} ) + A_{12} (y_2^{n+1} - y_2^{n} ) = 0 .
\end{split}
\end{equation}
The data, shown in Fig.~\ref{f-4}, demonstrate an unprincipled improvement in accuracy for the alternating-triangle method scheme.

Of particular interest are the results of the approximate solution of the test problem using second-order accuracy decomposition-composition schemes.
The comparison is made with the symmetric scheme ($\sigma = 1/2$ in (\ref{2.7})), when the reference solution is determined from Eqs.
\[
\frac{\overline{\bm y}^{n+1} - \overline{\bm y}^{n}}{\tau } +
{\bm A} \frac{\overline{\bm y}^{n+1} + \overline{\bm y}^{n}}{2 } = 0 .
\]
For $\bm z^n = \overline{\bm y}^n - \bm y^n $ we have the equation
\[
\frac{{\bm z}^{n+1} - {\bm z}^{n}}{\tau } +
{\bm A} \frac{{\bm z}^{n+1} + {\bm z}^{n}}{\tau }  = \bm \psi^n .
\]

We start with the alternating-triangle scheme of (\ref{3.3}), (\ref{3.9}), (\ref{3.11}) when $\sigma = 1/2$.
We write such an operator analog of the Peaceman-Rachford scheme \cite{PeacemanRachford1955}  as
\[
\begin{split}
\left (\bm I + \frac{1}{2} \tau \bm A_1 \right ) \bm y^{n+1/2} - \left (\bm I - \frac{1}{2} \tau \bm A_2 \right ) \bm y^{n} = 0 , \\\bm
\left (\bm I + \frac{1}{2} \tau \bm A_2 \right ) \bm y^{n+1} - \left (\bm I - \frac{1}{2} \tau \bm A_1 \right ) \bm y^{n+1/2} = 0 .
\end{split}
\]
In this case, the splitting error is estimated as
\[
\bm \psi^n = - \frac{1}{4} \tau^2 \bm A_1 \bm A_2 \frac{\bm y^{n+1} - \bm y^{n}}{\tau } .
\]
For the individual components of the solution, we have
\begin{equation}\label{5.11}
\begin{split}
& \left (1+\frac{1}{4} \tau A_{11} \right ) y_1^{n+1/2}
- \left (1-\frac{1}{4} \tau A_{11} \right ) y_1^{n} + \frac{1}{2} \tau A_{12} y_2^{n} = 0, \\
& \left (1+\frac{1}{4} \tau A_{22} \right ) y_2^{n+1/2}
- \left (1-\frac{1}{4} \tau A_{22} \right ) y_2^{n} + \frac{1}{2} \tau A_{21} y_1^{n+1/2} = 0, \\
& \left (1+\frac{1}{4} \tau A_{22} \right ) y_2^{n+1}
- \left (1-\frac{1}{4} \tau A_{22} \right ) y_2^{n+1/2} + \frac{1}{2} \tau A_{21} y_1^{n+1/2} = 0, \\
& \left (1+\frac{1}{4} \tau A_{11} \right ) y_1^{n+1}
- \left (1-\frac{1}{4} \tau A_{11} \right ) y_1^{n+1/2} + \frac{1}{2} \tau A_{12} y_2^{n+1} = 0 .
\end{split}
\end{equation}

\begin{figure}[htbp]
\centering
\includegraphics[width=0.4\linewidth]{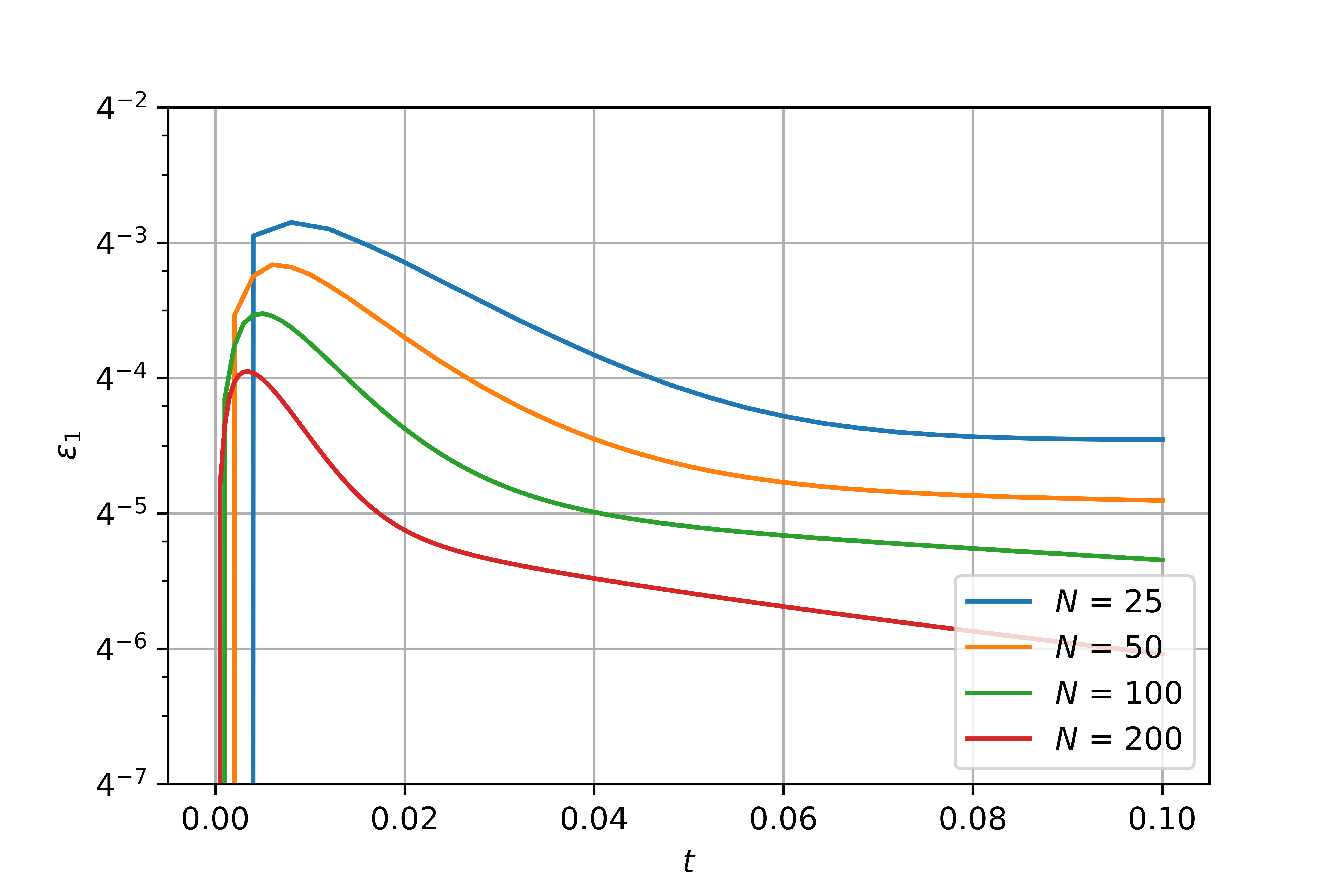} \includegraphics[width=0.4\linewidth]{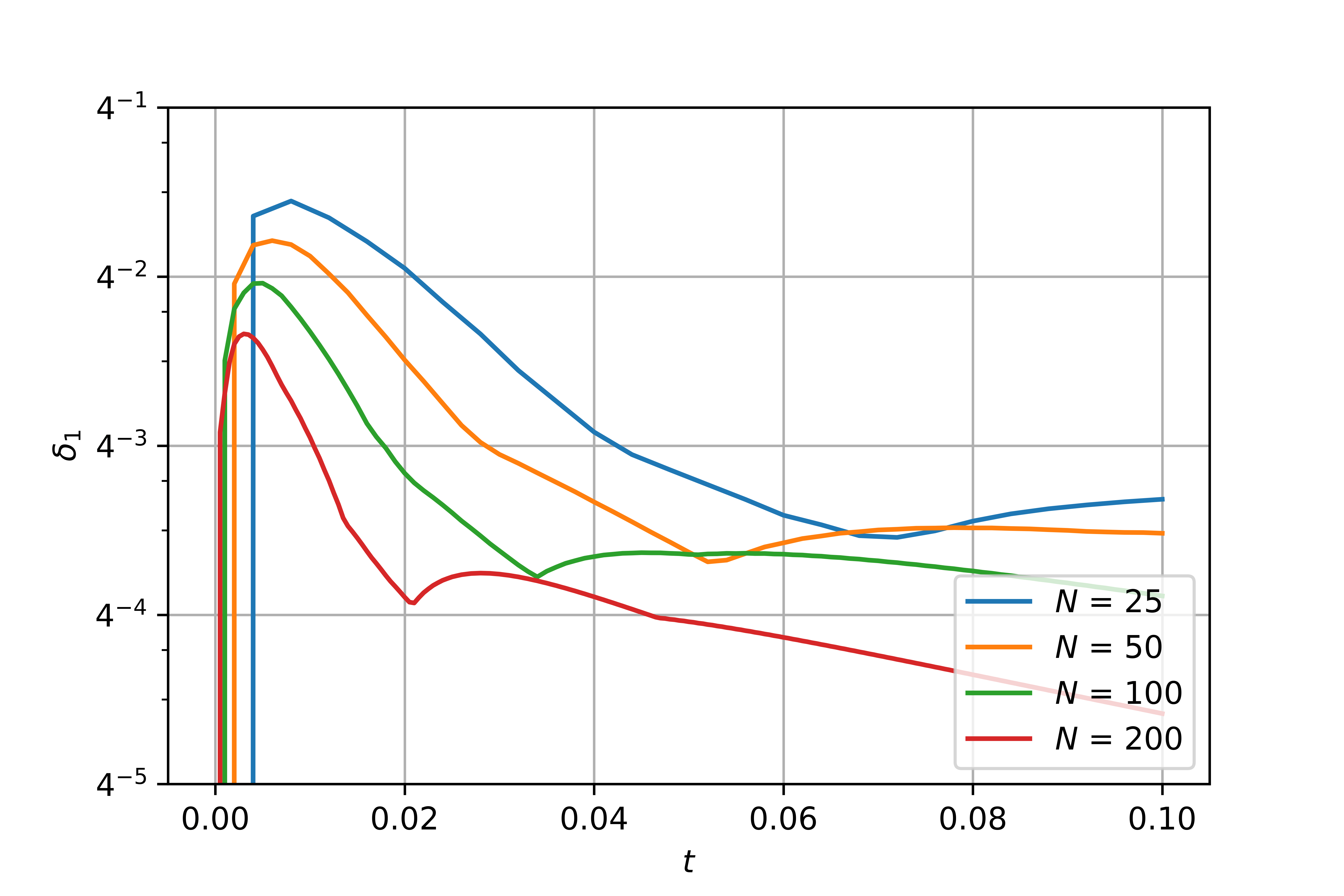} \\
\includegraphics[width=0.4\linewidth]{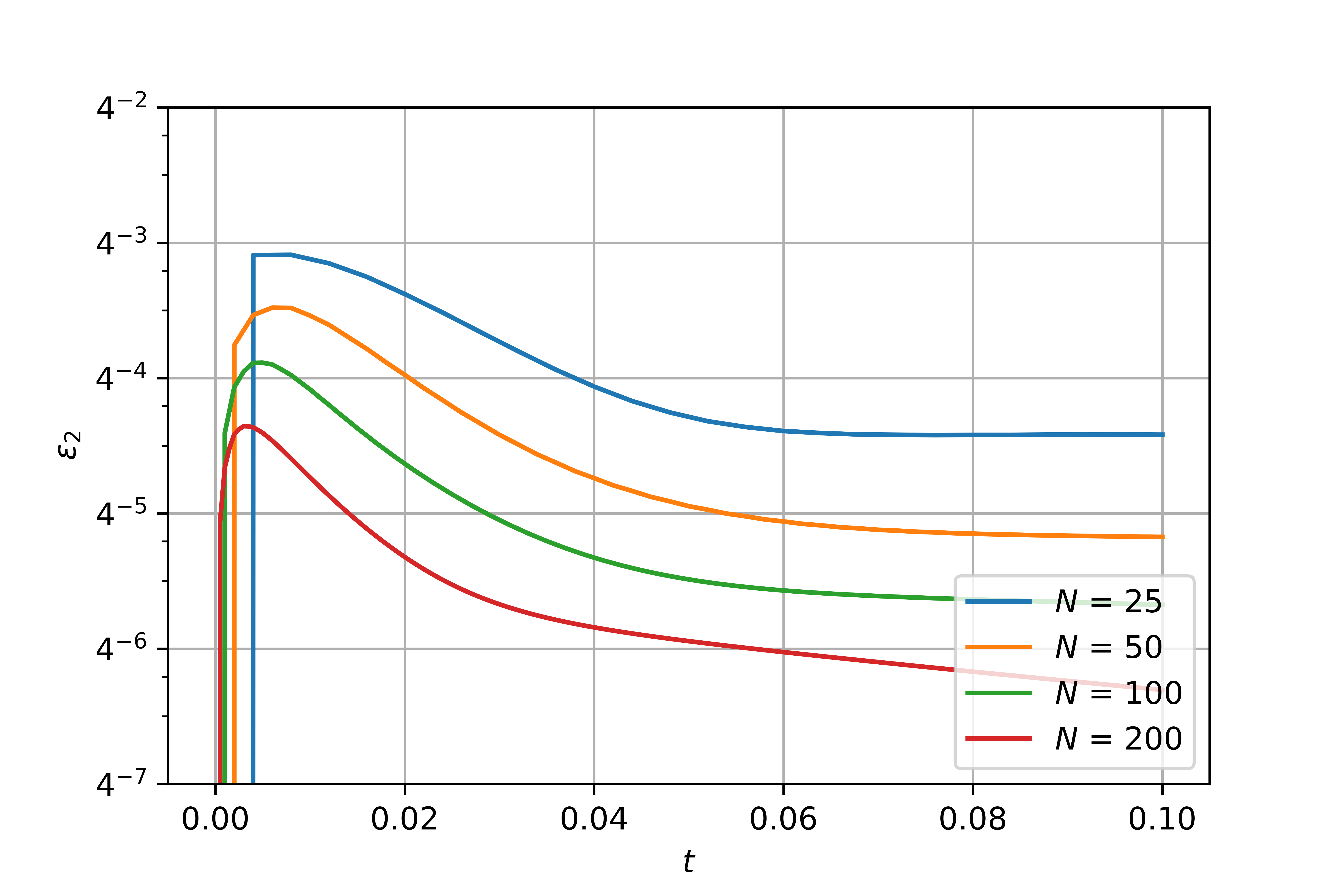} \includegraphics[width=0.4\linewidth]{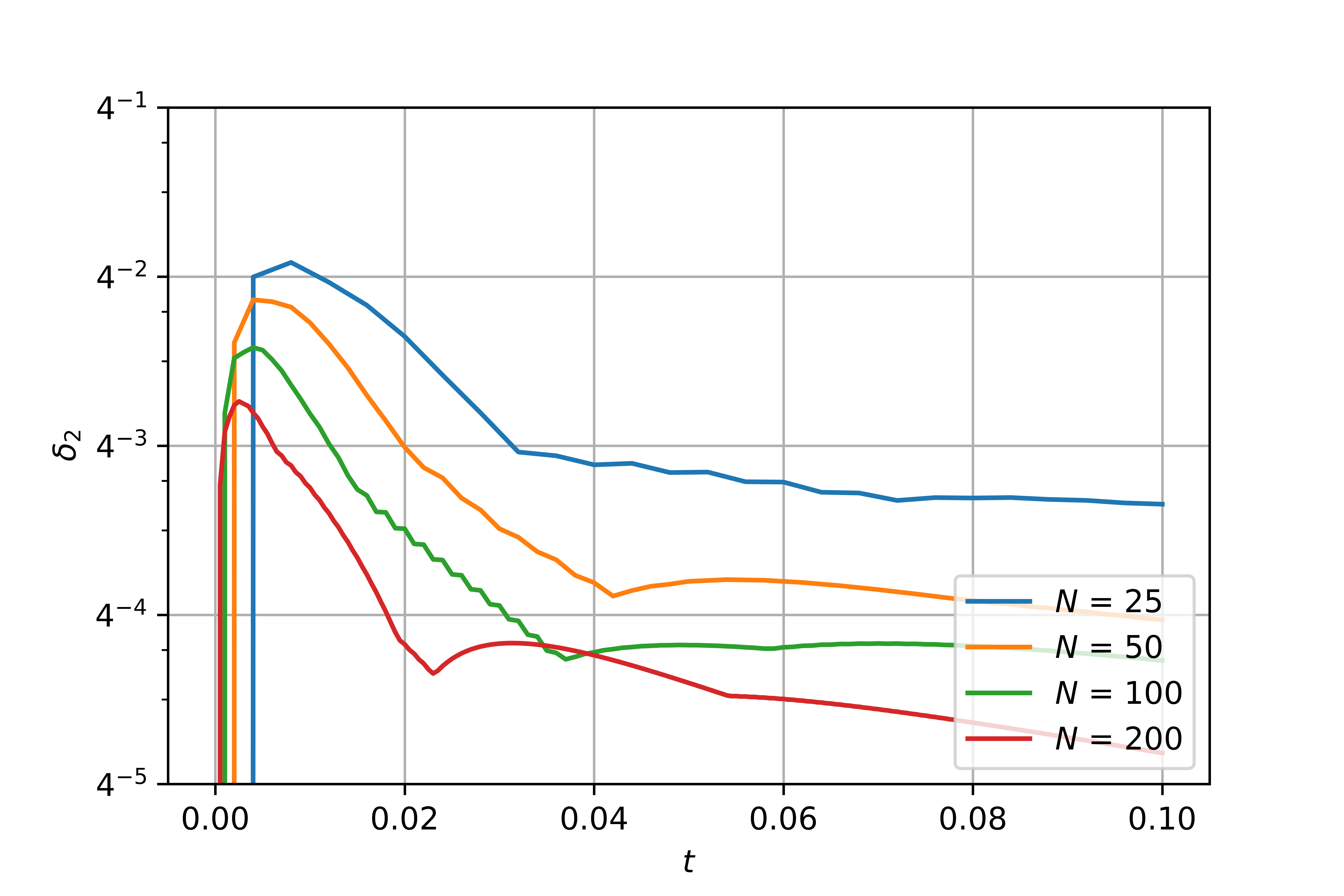}
\caption{Decomposition-composition scheme (\ref{5.11}).}
\label{f-5}
\end{figure}

The results of comparing the decomposition-composition scheme (\ref{5.11}) and the standard symmetric scheme are shown in Fig.~\ref{f-5}.
Compared to the scheme (\ref{5.11}), we have two times smaller $\sigma$ value and could expect to improve the closeness to the reference solution by a factor of four.

We will first illustrate the possibilities of using three-level decomposition-composition schemes for the alternating-triangle method scheme (\ref{3.19}) when $\sigma = 1/2$.
The realization is based on the two-level Peaceman-Rachford scheme with given right-hand sides:
\[
\begin{split}
\left (\bm I + \frac{1}{2} \tau \bm A_1 \right ) \bm y^{n+1/2} - \left (\bm I - \frac{1}{2} \tau \bm A_2 \right ) \bm y^{n} = \bm f_1^n, \\
\left (\bm I + \frac{1}{2} \tau \bm A_2 \right ) \bm y^{n+1} - \left (\bm I - \frac{1}{2} \tau \bm A_1 \right ) \bm y^{n+1/2} =  \bm f_2^n  .
\end{split}
\]
If we eliminate $\bm y^{n+1/2}$, we get
\[
\left (\bm I + \frac{1}{2} \tau \bm A_1 \right ) \left (\bm I + \frac{1}{2} \tau \bm A_2 \right ) \frac{{\bm y}^{n+1} - {\bm y}^{n}}{\tau }
+ \bm A {\bm y}^{n}
= \bm f^n ,
\]
where
\[
\tau \bm f^n = \left (\bm I - \frac{1}{2} \tau \bm A_1 \right ) \bm f_1^n + \left (\bm I + \frac{1}{2} \tau \bm A_1 \right ) \bm f_2^n .
\]
Taking into account (\ref{3.18}), (\ref{3.19}), we have
\[
\bm f^n = \frac{1}{4} \tau \bm A_1 \bm A_2 (\bm y^{n} - \bm y^{n-1}) .
\]
For this reason, let us assume
\[
\bm f_2^n = - \bm f_1^n = \frac{1}{4} \tau \bm A_2 (\bm y^{n} - \bm y^{n-1}) .
\]
For the splitting error, we obtain
\[
\bm \psi^n = - \frac{1}{4} \tau^3 \bm A_1 \bm A_2 \frac{\bm y^{n+1} - 2 \bm y^{n} + \bm y^{n-1}}{\tau^2} .
\]

For the considered system from the two equations, we obtain
\begin{equation}\label{5.12}
\begin{split}
& \left (1+\frac{1}{4} \tau A_{11} \right ) y_1^{n+1/2}
- \left (1-\frac{1}{4} \tau A_{11} \right ) y_1^{n} + \frac{1}{2} \tau A_{12} y_2^{n} = - \varphi_1^n, \\
& \left (1+\frac{1}{4} \tau A_{22} \right ) y_2^{n+1/2}
- \left (1-\frac{1}{4} \tau A_{22} \right ) y_2^{n} + \frac{1}{2} \tau A_{21} y_1^{n+1/2} = - \varphi_2^n, \\
& \left (1+\frac{1}{4} \tau A_{22} \right ) y_2^{n+1}
- \left (1-\frac{1}{4} \tau A_{22} \right ) y_2^{n+1/2} + \frac{1}{2} \tau A_{21} y_1^{n+1/2} = \varphi_2^n, \\
& \left (1+\frac{1}{4} \tau A_{11} \right ) y_1^{n+1}
- \left (1-\frac{1}{4} \tau A_{11} \right ) y_1^{n+1/2} + \frac{1}{2} \tau A_{12} y_2^{n+1} = \varphi_1^n ,
\end{split}
\end{equation}
where
\[
\begin{split}
& \varphi_1^n = \frac{1}{8} \tau \big (A_{11} (y_1^{n} - y_1^{n-1}) + 2A_{12} (y_2^{n} - y_2^{n-1}) \big ) , \\
& \varphi_2^n = \frac{1}{8} \tau A_{22} (y_2^{n} - y_2^{n-1}) .
\end{split}
\]

\begin{figure}[htbp]
\centering
\includegraphics[width=0.4\linewidth]{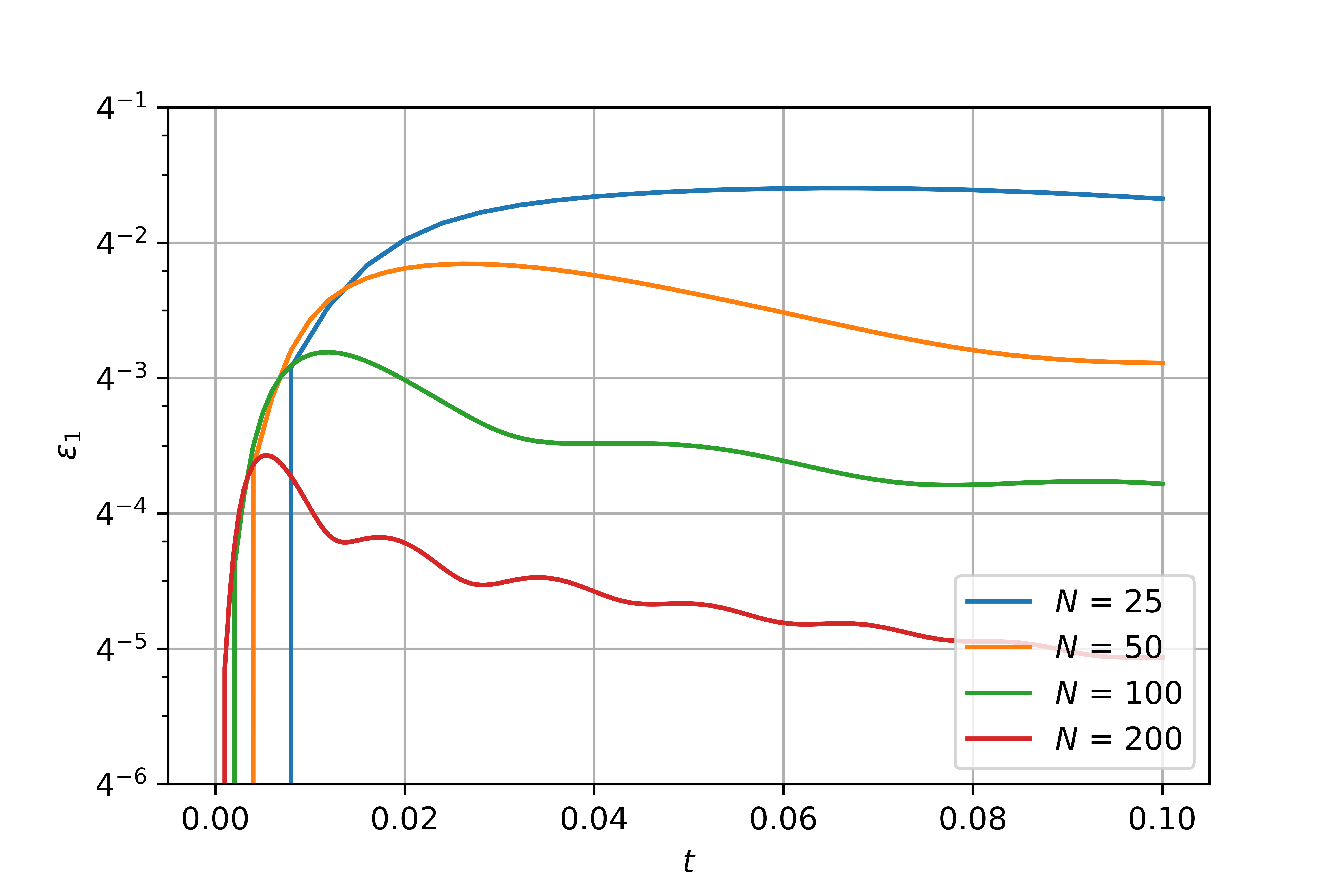} \includegraphics[width=0.4\linewidth]{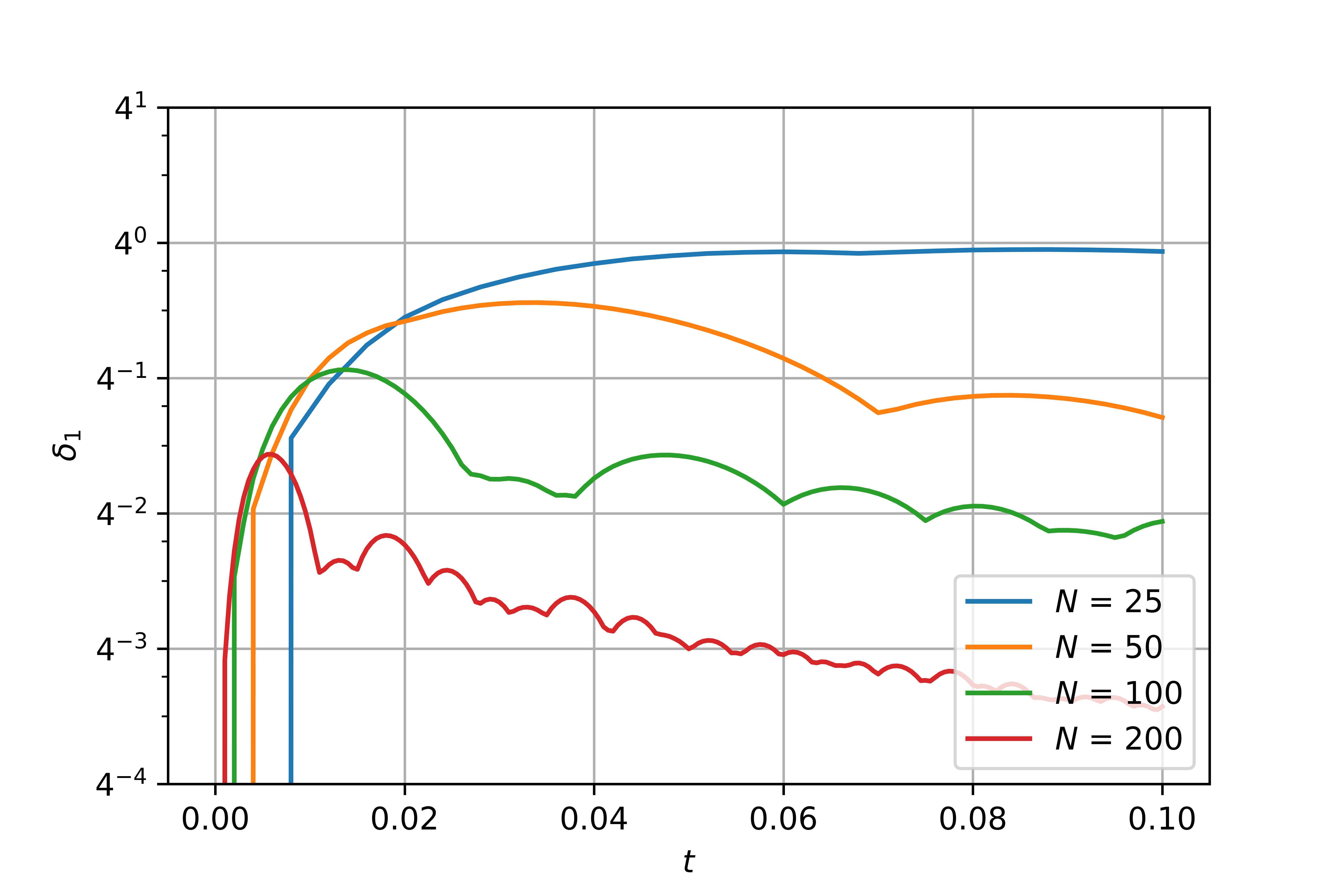} \\
\includegraphics[width=0.4\linewidth]{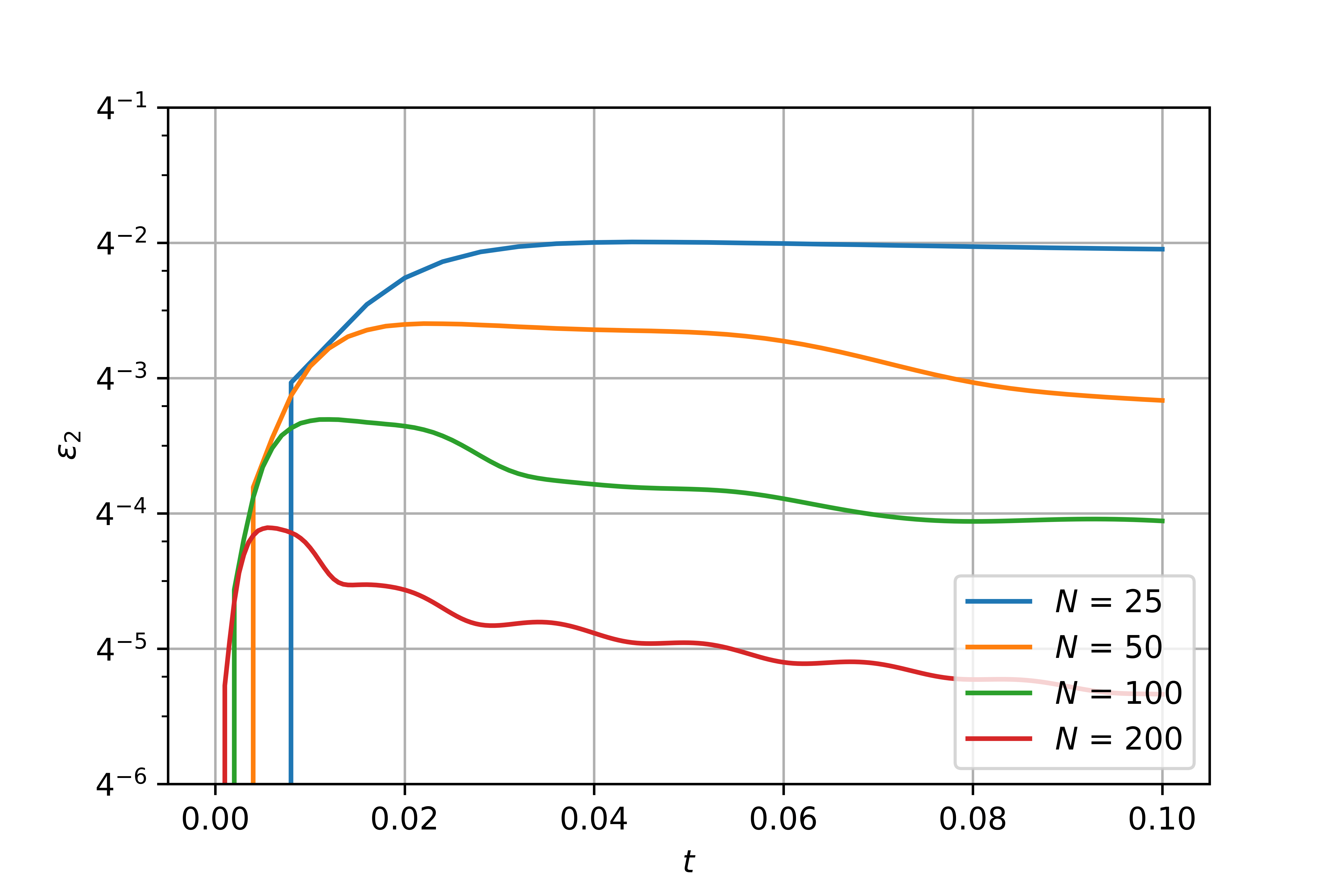} \includegraphics[width=0.4\linewidth]{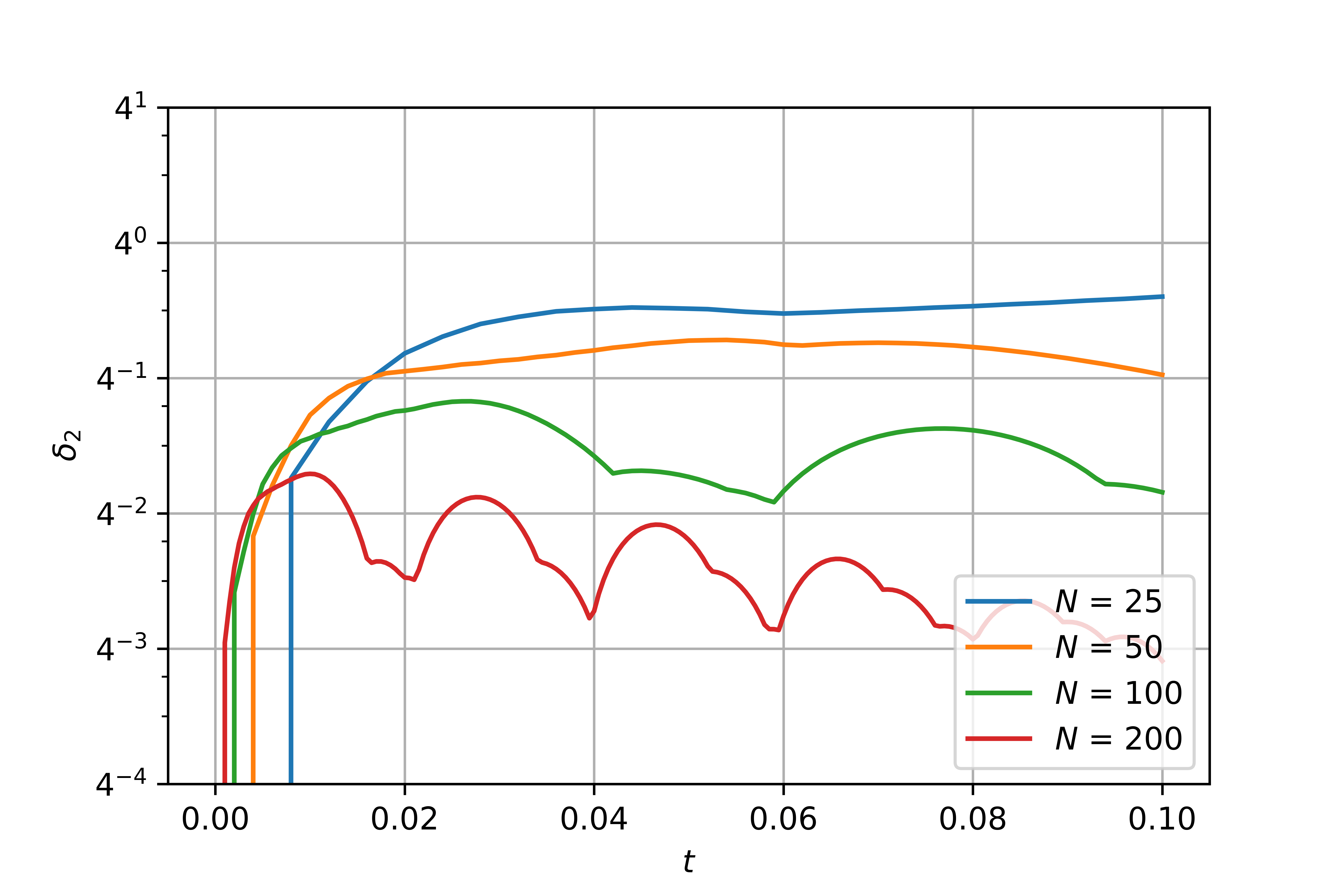}
\caption{Decomposition-composition scheme (\ref{5.12}).}
\label{f-6}
\end{figure}

The calculated data for the scheme (\ref{5.12}) are presented in Fig.~\ref{f-6}.
The main error is observed at the initial stage.
At the developed stage of the process, the result for the symmetric scheme without splitting becomes more close to the original.

When extracting the diagonal part of the problem operator, we obtained the most accurate results using a three-level decomposition-composition scheme.
We associate the reference solution with the three-level scheme of the second order of accuracy
\[
\frac{\overline{\bm y}^{n+1} - \overline{\bm y}^{n-1}}{\tau } +
\frac{1}{4} {\bm A} (\overline{\bm y}^{n+1} - 2 \overline{\bm y}^{n} + \overline{\bm y}^{n-1} ) = 0 .
\]
The splitting error is
\[
\bm \psi^n = \tau \left (\frac{1}{4} \bm A - \sigma \bm D \right ) (\bm y^{n+1} - 2 \bm y^{n} + \bm y^{n-1} ) .
\]

Defining
\[
\overline{\bm v}^{n+1} = \frac{1}{2} (\overline{\bm y}^{n+1} + \overline{\bm y}^{n}) ,
\]
we have the standard symmetric scheme
\[
\frac{\overline{\bm v}^{n+1} - \overline{\bm v}^{n}}{\tau } +
{\bm A} \frac{\overline{\bm v}^{n+1} +\overline{\bm v}^{n}}{2} = 0 .
\]
Because of this, we can compare the solution of the three-level decomposition-composition scheme (\ref{3.16}) with the solution of the usual two-level symmetric scheme, as well as other decomposition-composition schemes of the second order of accuracy.

For the test problem we consider, we use the decomposition-composition scheme (\ref{3.16}) when $\sigma = 1/2$.
In this case, for the individual components of the approximate solution, we have
\begin{equation}\label{5.13}
\begin{split}
\frac{y_1^{n+1} - y_1^{n}}{2\tau } + A_{11} (y_1^{n+1} - 2y_1^{n} + y_1^{n-1} ) + A_{11} y_1^{n} + A_{12} y_2^{n} = 0, \\
\frac{y_2^{n+1} - y_2^{n}}{2\tau } + A_{22} (y_2^{n+1} - 2y_2^{n} + y_2^{n-1} ) + A_{21} y_1^{n} + A_{22} y_2^{n} = 0 .
\end{split}
\end{equation}

\begin{figure}[htbp]
\centering
\includegraphics[width=0.4\linewidth]{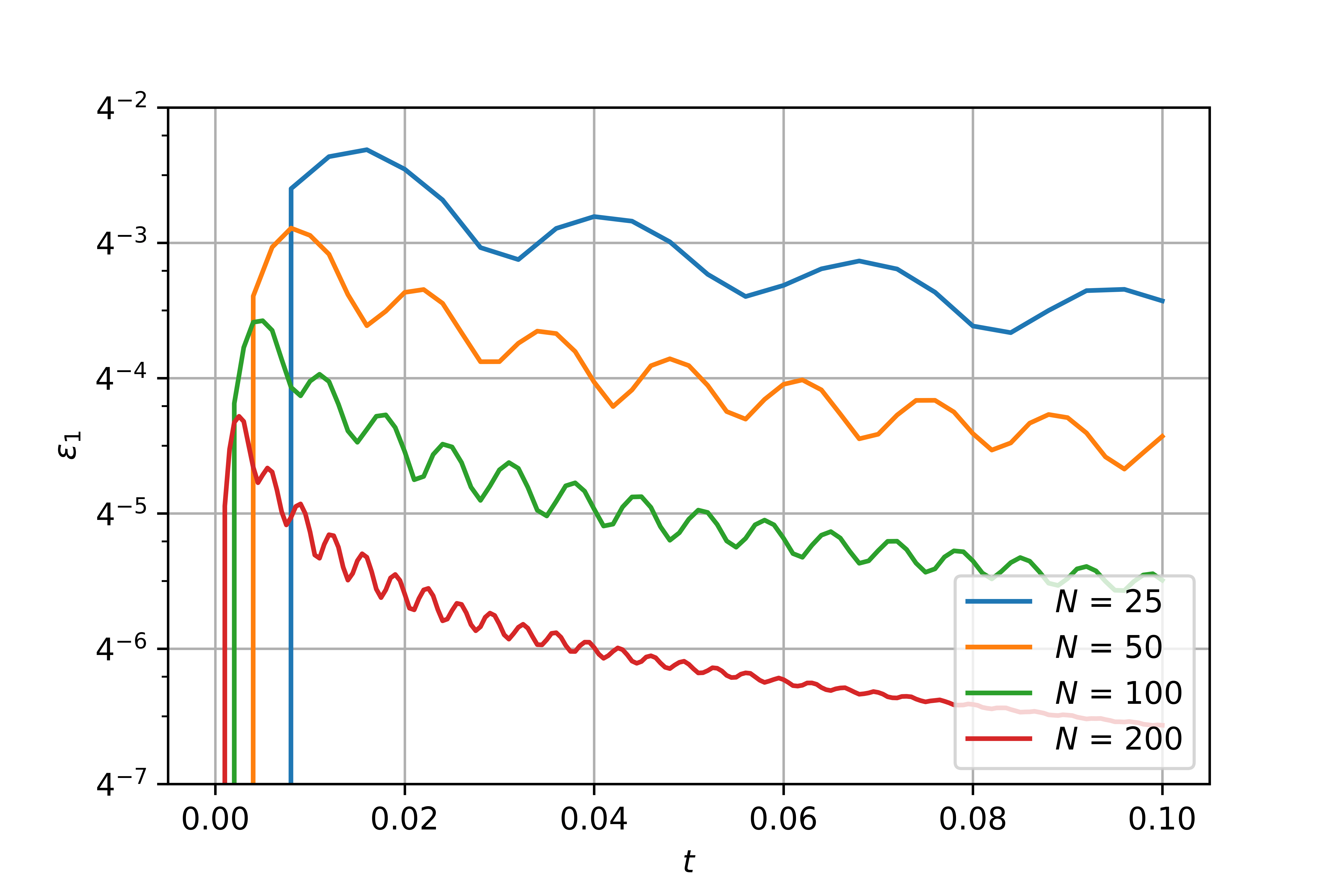} \includegraphics[width=0.4\linewidth]{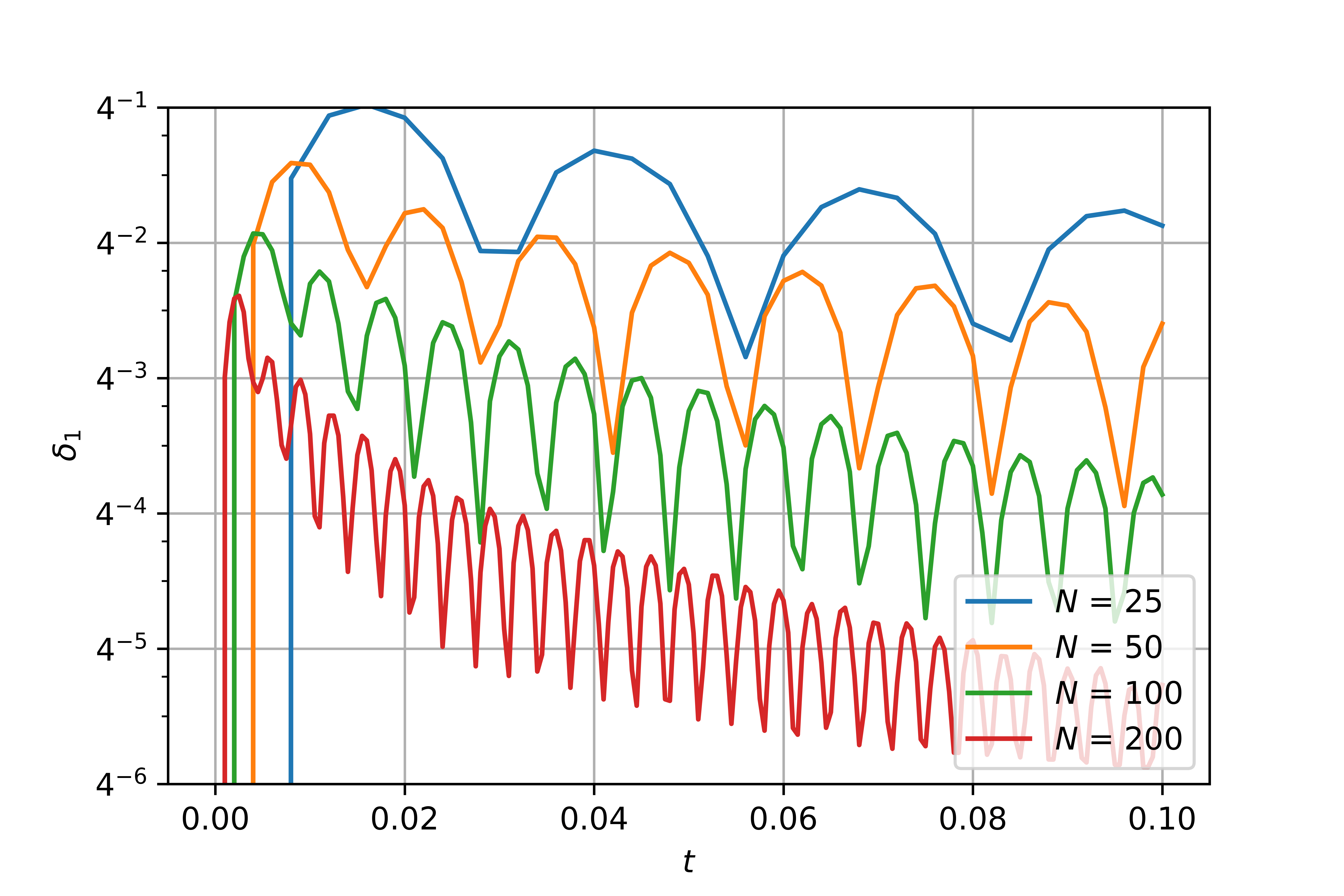} \\
\includegraphics[width=0.4\linewidth]{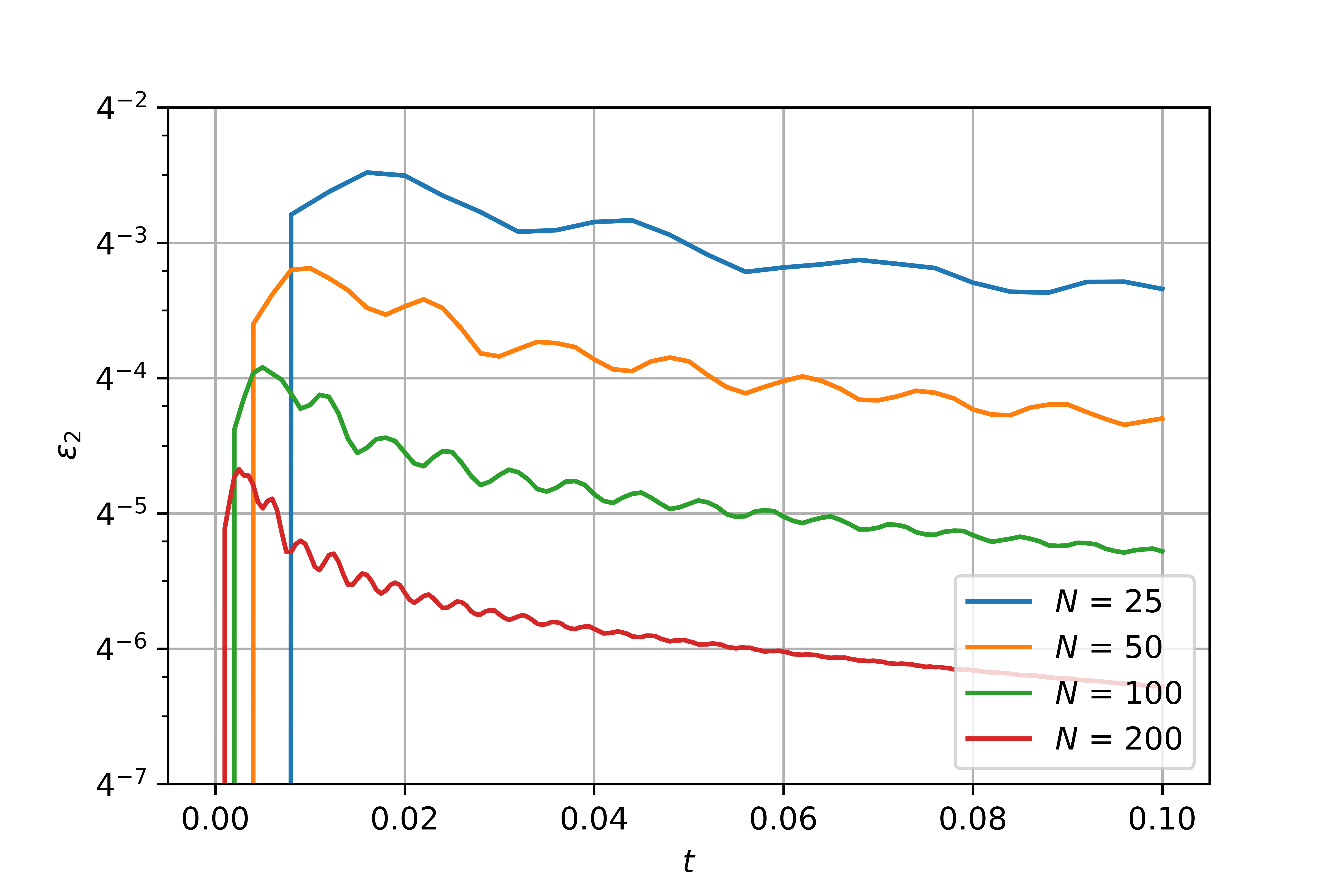} \includegraphics[width=0.4\linewidth]{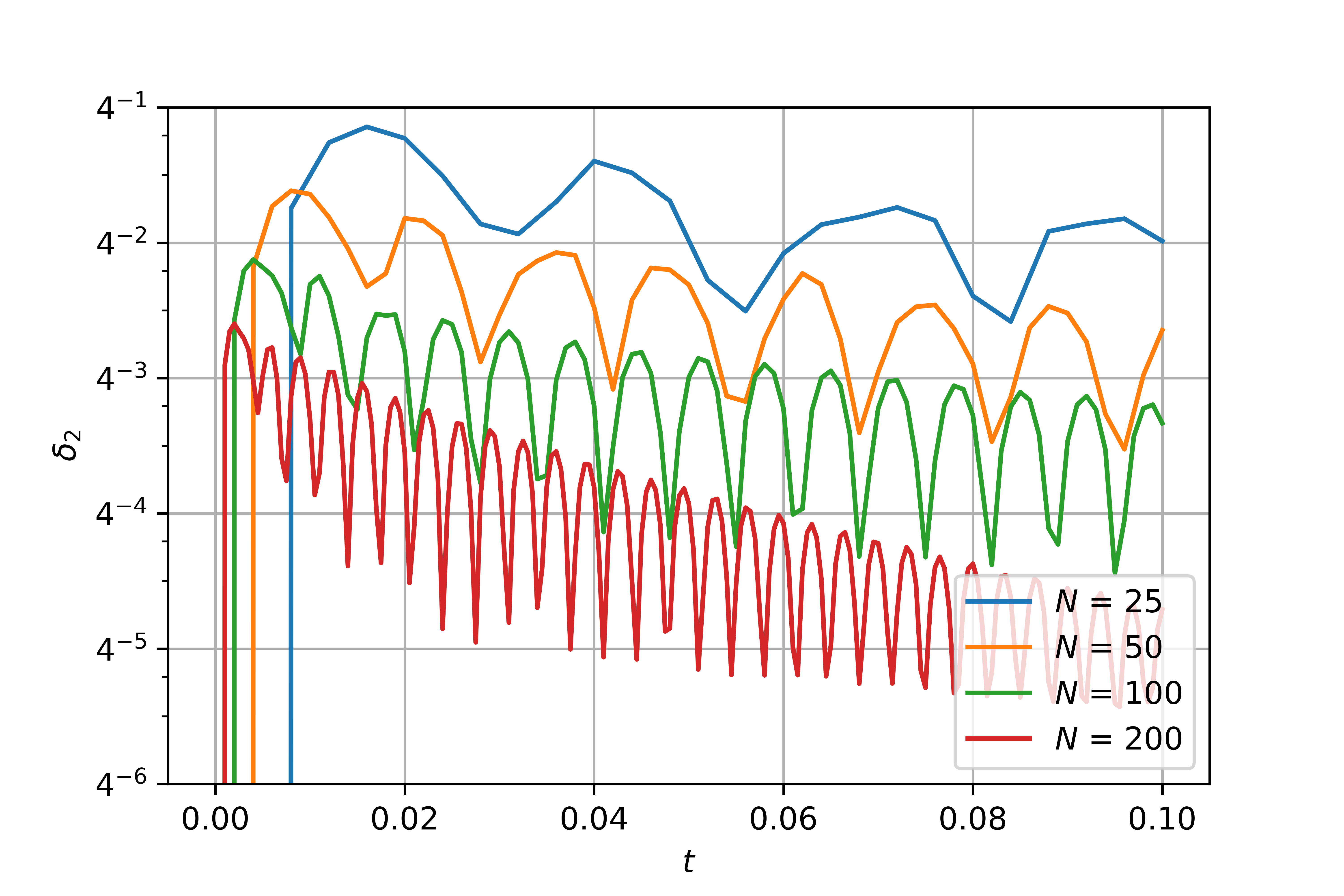}
\caption{Decomposition-composition scheme (\ref{5.13}).}
\label{f-7}
\end{figure}

A comparison of the scheme (\ref{5.13}) and the standard symmetric scheme is shown in Fig.~\ref{f-7}.
The accuracy of this scheme is noticeably higher than that of other decomposition-composition schemes of the second order of accuracy.
This is due to the peculiarities of the problem we consider.
For the alternating-triangular method scheme, the splitting error is proportional to the product of the splitting operators: $\bm A_1 \bm A_2 \bm w$,
where $\bm w$ --- the first or second-time derivative of the solution.
Our test problem is characterized by discontinuous diffusion coefficients, contributing to the approximate solution's main error.
For the scheme (\ref{5.13}) with separation of the diagonal part of the problem operator, the splitting error is due only to the operator terms and not to their product.

We also present results for decomposition-composition schemes based on the additive representation of the unit operator (\ref{4.2}).
For our test problem (\ref{5.2})--(\ref{5.3}) we have a two-component split
\begin{equation}\label{5.14}
\bm A = \bm A_1 + \bm A_2 .
\end{equation}
When decomposing by strings (see (\ref{4.3}) and (\ref{4.6})) we get
\begin{equation}\label{5.15}
\bm A_1 = \begin{pmatrix}
A_{11} & A_{12} \\
0 & 0
\end{pmatrix}  ,
\quad \bm A_2 = \begin{pmatrix}
0 & 0 \\
A_{21} & A_{22}
\end{pmatrix} .
\end{equation}
When decomposing by columns (see (\ref{4.4}) and (\ref{4.7})) we have
\begin{equation}\label{5.16}
\bm A_1 = \begin{pmatrix}
A_{11} & 0 \\
A_{21} & 0
\end{pmatrix}  ,
\quad \bm A_2 = \begin{pmatrix}
0 & A_{12} \\
0 & A_{22}
\end{pmatrix} .
\end{equation}

The scheme of component-wise splitting at decomposition (\ref{5.14}) takes the following form
\begin{equation}\label{5.17}
\begin{split}
& (\bm I + \sigma \tau \bm A_1 ) \frac{\bm y^{n+1/2} -  \bm y^{n}}{\tau} + \bm A_1  \bm y^{n} = 0 , \\
& (\bm I + \sigma \tau \bm A_2 ) \frac{\bm y^{n+1} -  \bm y^{n+1/2}}{\tau} + \bm A_2\bm y^{n+1/2} = 0  .
\end{split}
\end{equation}
When decomposing (\ref{5.15}) for individual components from (\ref{5.17}), we obtain
\begin{equation}\label{5.18}
\begin{split}
& (1+ \sigma \tau A_{11} ) \frac{y_1^{n+1} - y_1^{n}}{\tau } + A_{11} y_1^{n} + A_{12} y_2^{n} = 0, \\
& (1+ \sigma \tau A_{22} ) \frac{y_2^{n+1} - y_2^{n}}{\tau } + A_{21} y_1^{n+1} + A_{22} y_2^{n} = 0 ,
\end{split}
\end{equation}

\begin{figure}[htbp]
\centering
\includegraphics[width=0.4\linewidth]{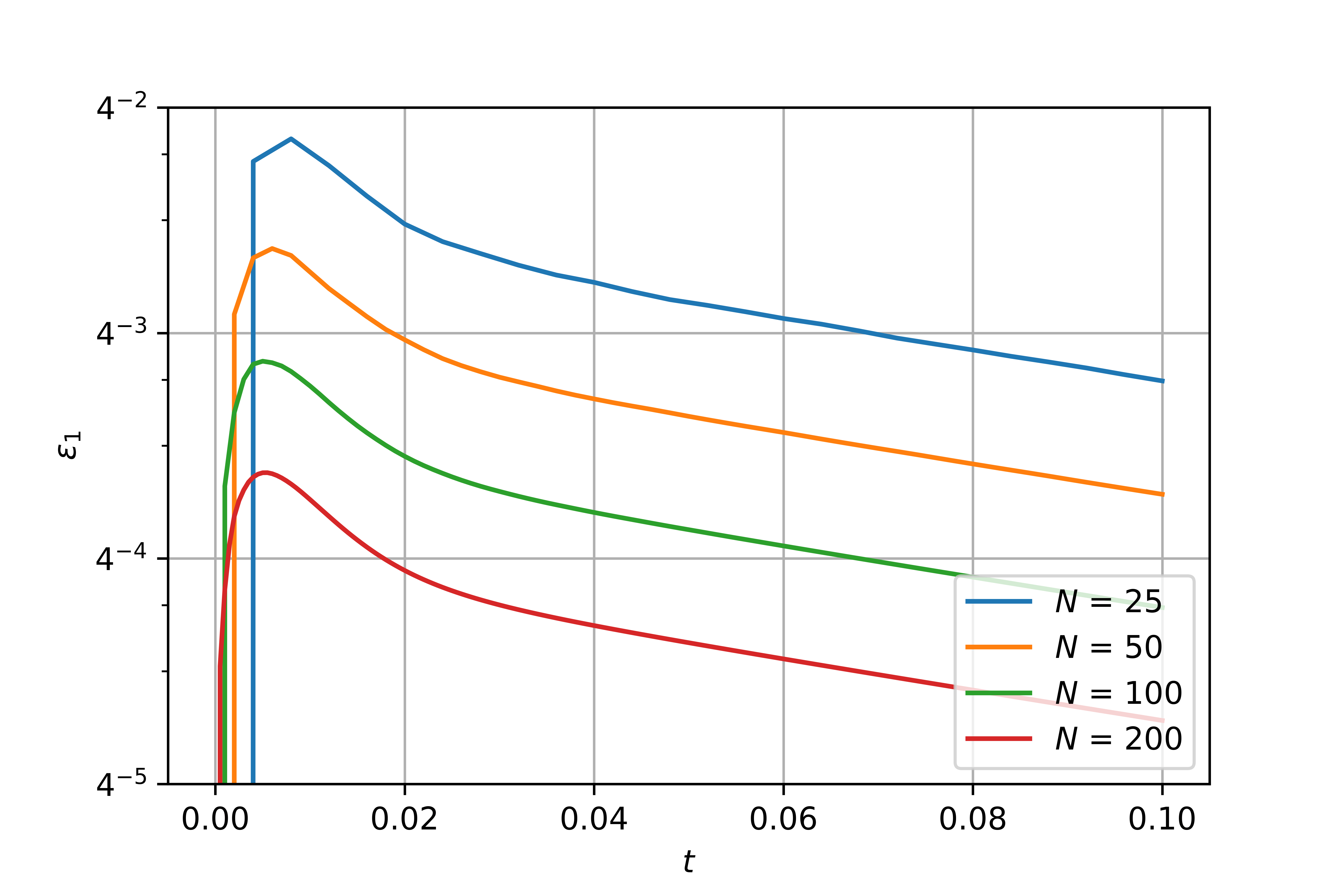} \includegraphics[width=0.4\linewidth]{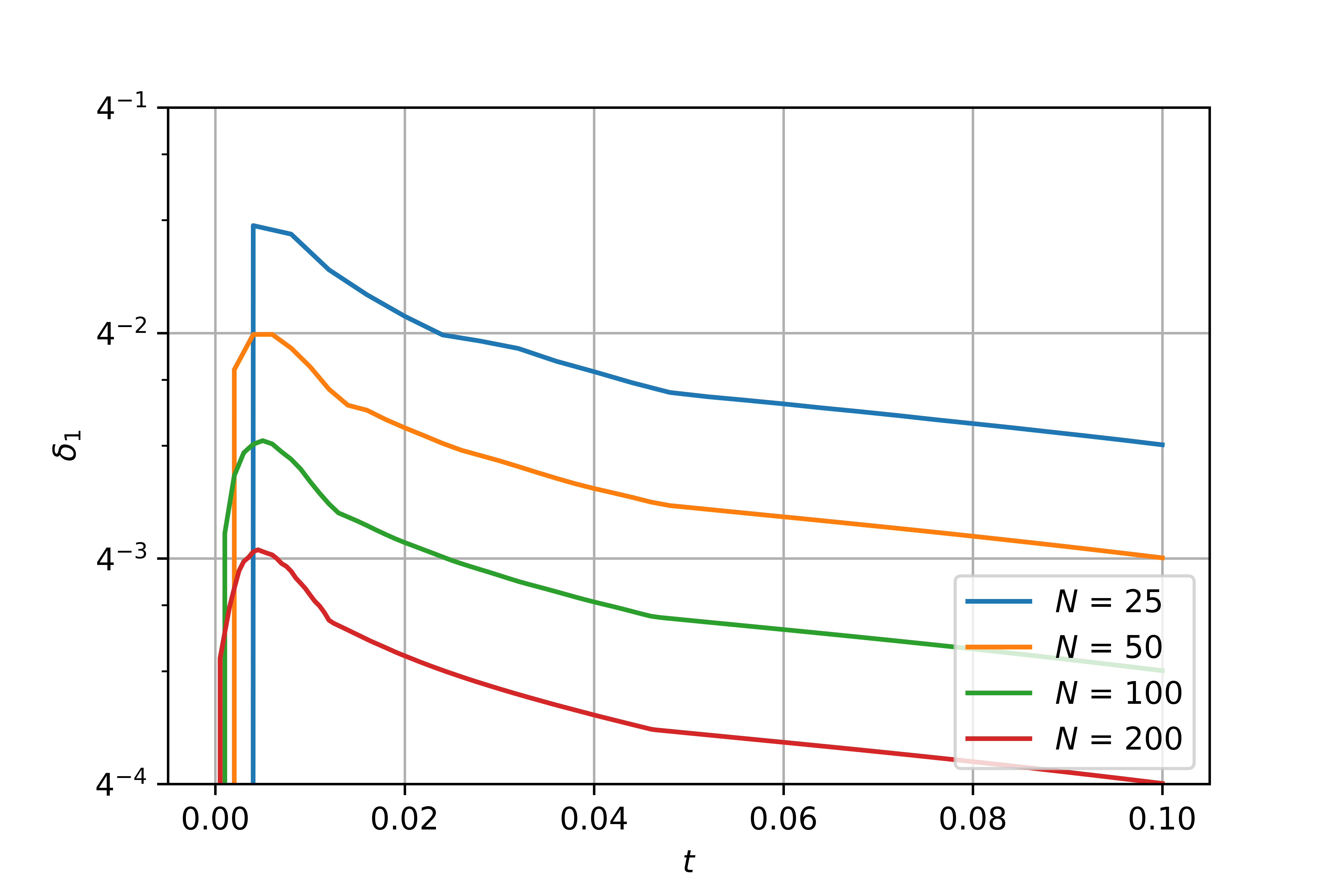} \\
\includegraphics[width=0.4\linewidth]{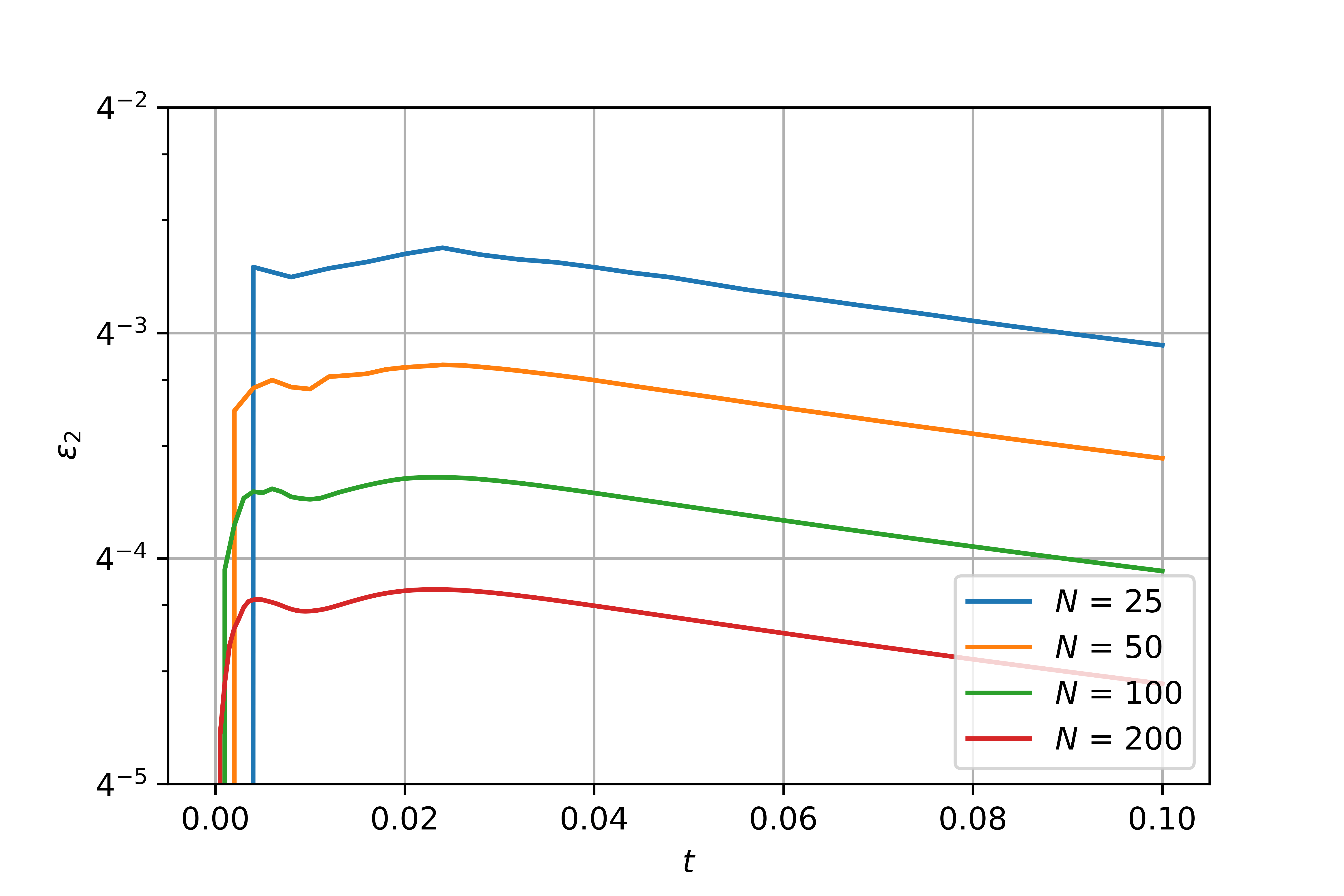} \includegraphics[width=0.4\linewidth]{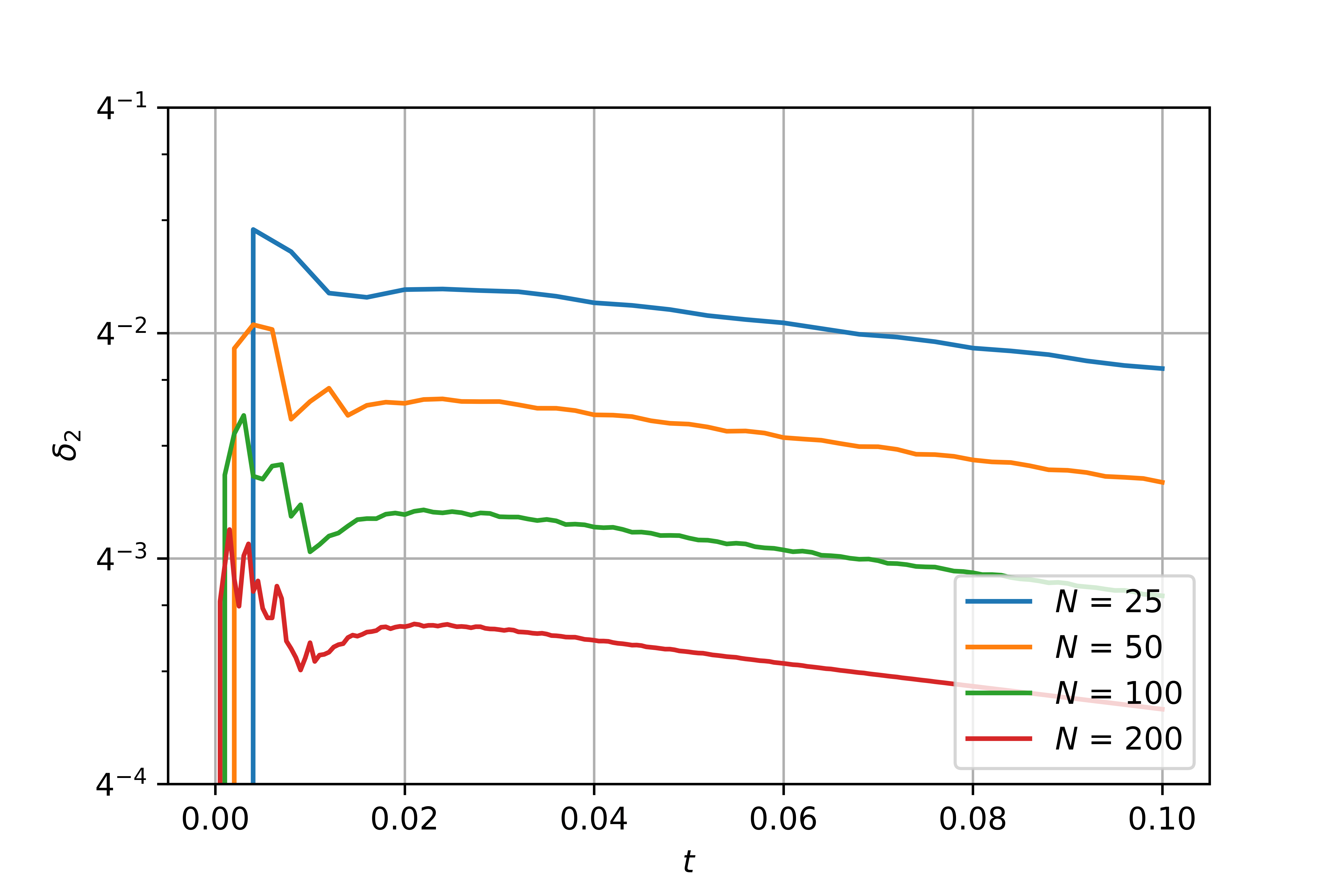}
\caption{Decomposition-composition scheme (\ref{5.18}).}
\label{f-8}
\end{figure}

At $\sigma = 1$, the decomposition-composition scheme (\ref{5.18}) coincides with the scheme (\ref{5.9}), which (see (\ref{3.8})) is connected with the transfer of the diagonal and lower triangular parts of the operator matrix $\bm A$ to a new level in time.
Fig.~\ref{f-8} shows the results of the comparison of the scheme (\ref{5.17}) at $\sigma = 1/2$ with the purely implicit scheme without splitting.
We can say that the solutions with the first order in $\tau$ are close to each other.

\begin{figure}[htbp]
\centering
\includegraphics[width=0.4\linewidth]{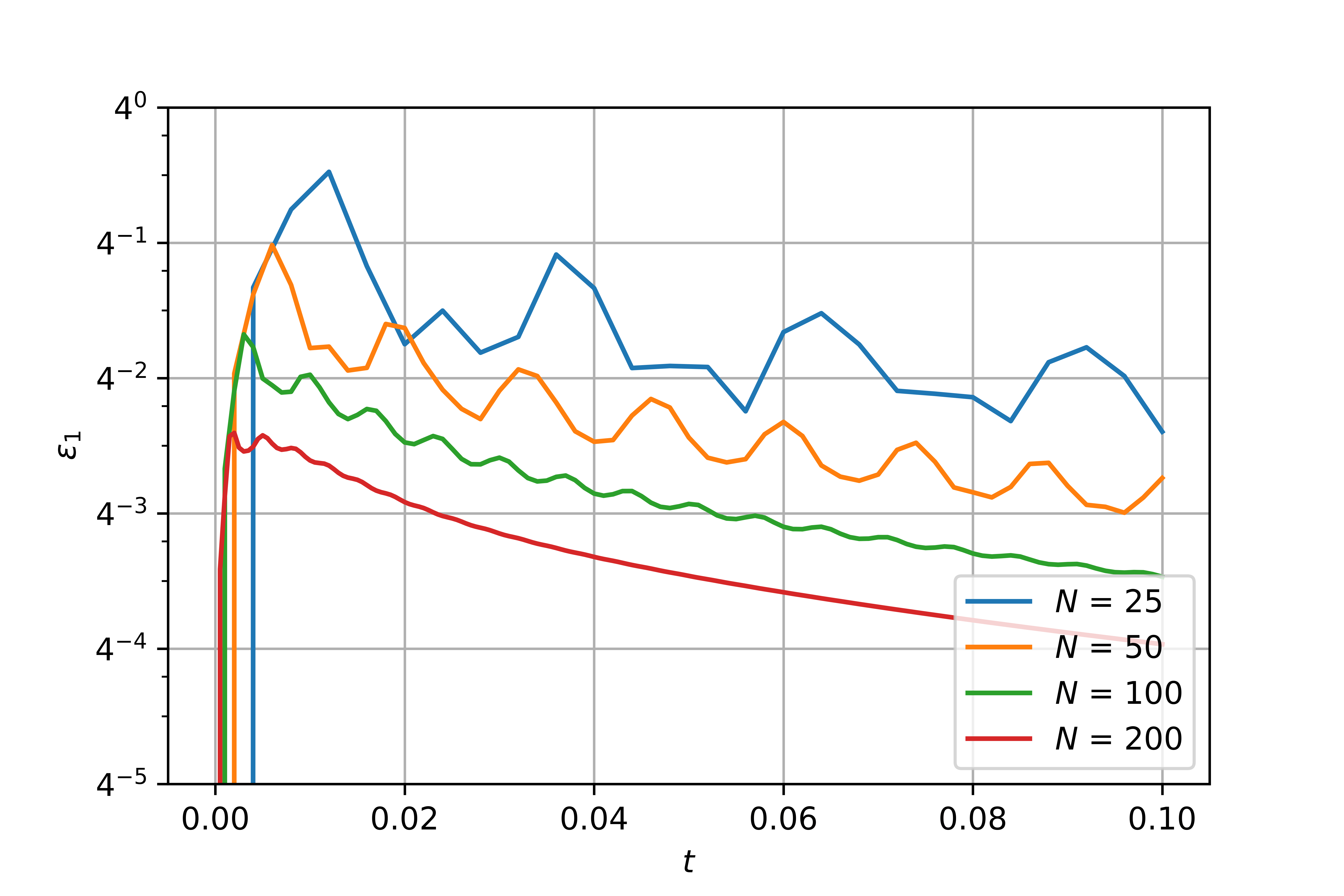} \includegraphics[width=0.4\linewidth]{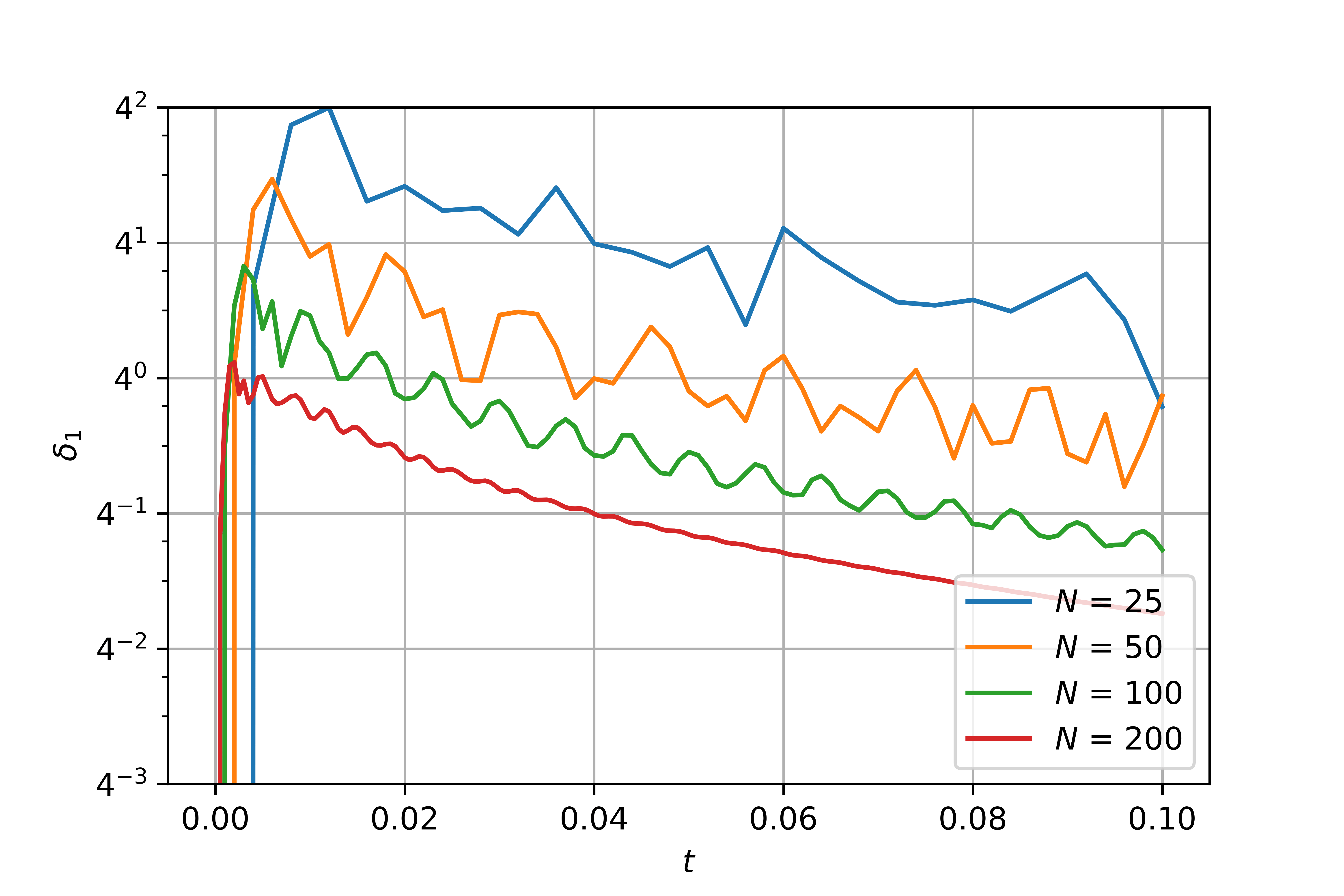} \\
\includegraphics[width=0.4\linewidth]{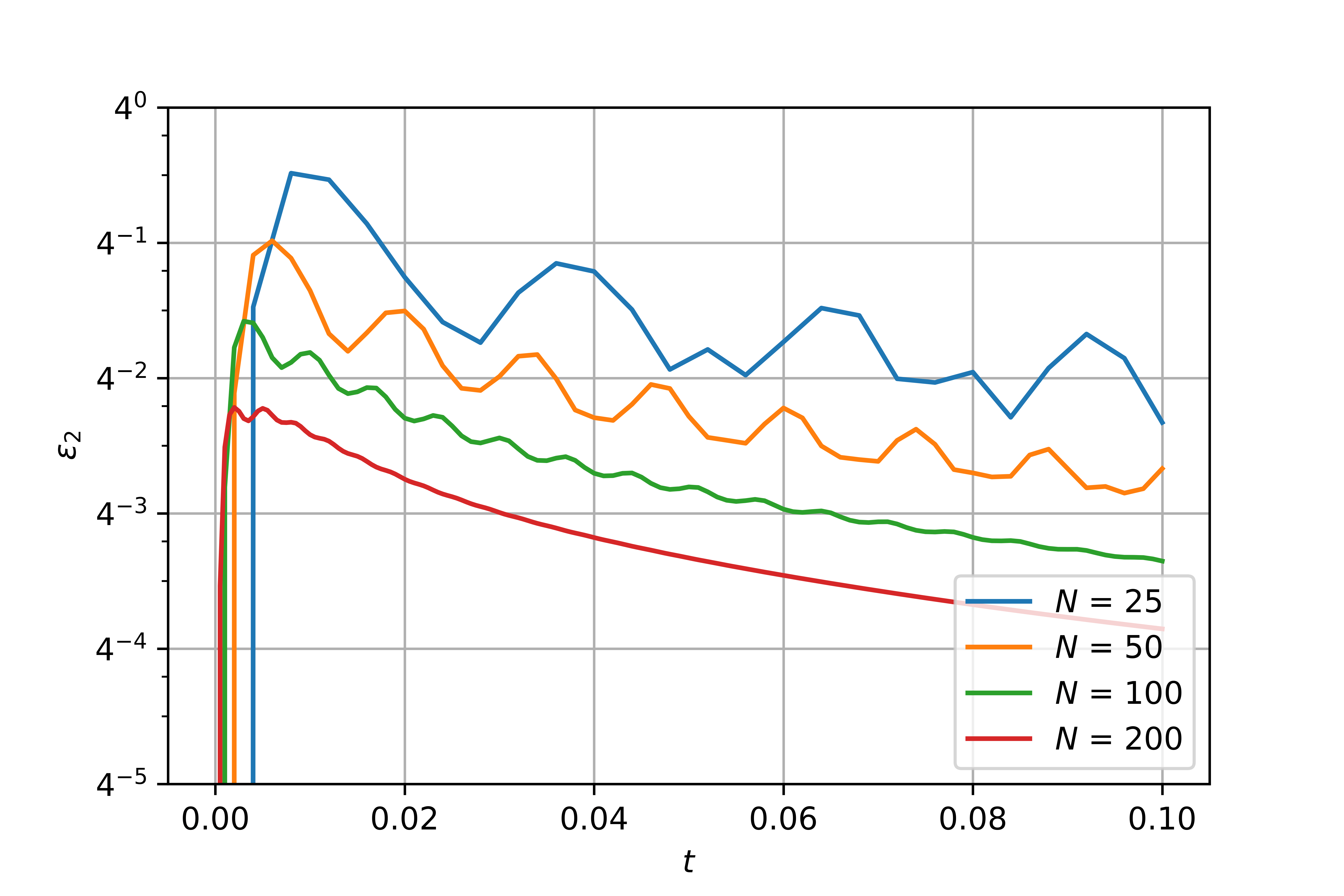} \includegraphics[width=0.4\linewidth]{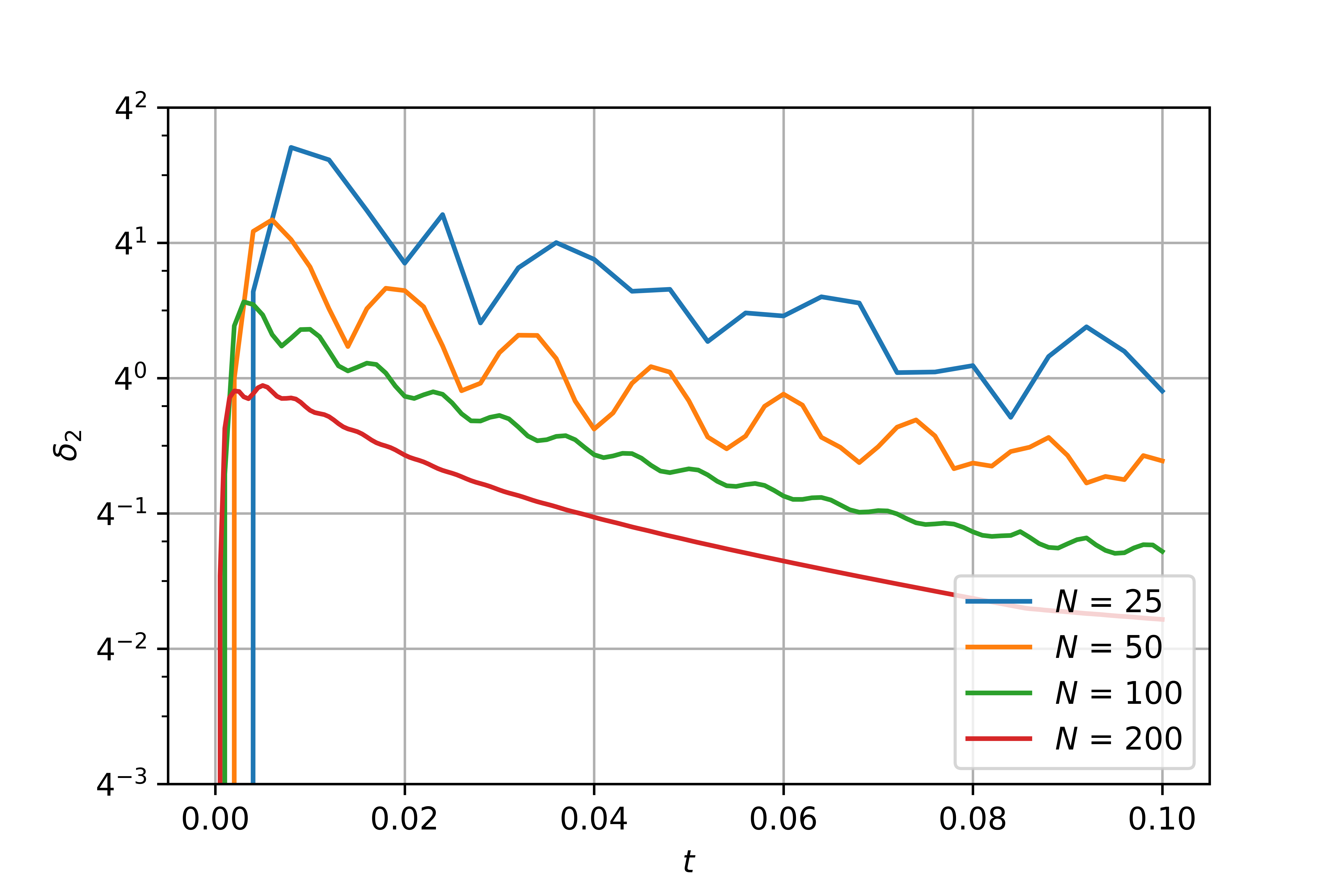}
\caption{Decomposition-composition scheme (\ref{5.19}).}
\label{f-9}
\end{figure}

When decomposing by columns (\ref{5.16}) for the individual components from (\ref{5.17}), we have
\begin{equation}\label{5.19}
\begin{split}
& (1+ \sigma \tau A_{11} ) \frac{y_1^{n+1/2} - y_1^{n}}{\tau } + A_{11} y_1^{n} = 0, \\
& y_2^{n+1/2} = y_2^{n} - \sigma \tau A_{21} (y_1^{n+1/2} - y_1^{n}) - \tau A_{21} y_1^{n} , \\
& (1+ \sigma \tau A_{22} ) \frac{y_2^{n+1} - y_2^{n+1/2}}{\tau } + A_{22} y_2^{n+1/2} = 0 , \\
& y_1^{n+1} = y_1^{n+1/2} - \sigma \tau A_{12} (y_2^{n+1} - y_2^{n+1/2}) - \tau A_{12} y_2^{n+1/2} .
\end{split}
\end{equation}
The computational data for this scheme at $\sigma = 1/2$ are shown in Fig.~\ref{f-9}.
For our test problem, the accuracy of the component-wise splitting scheme for row decomposition (Fig.~\ref{f-8}) is significantly higher than for column decomposition (Fig.~\ref{f-9}).
This is because the stability (and therefore convergence) of the scheme (\ref{5.18}) takes place in a stronger norm.

Our test problem illustrates the possibility of using regularized decomposition-composition schemes.
Given a two-component split (\ref{5.14}), the approximate solution is found from (see (\ref{4.23}))
\begin{equation}\label{5.20}
\frac{\bm y^{n+1} - \bm y^n}{\tau} + \sum_{\alpha = 1}^{2} (\bm I+\sigma \tau \bm A_\alpha )^{-1}  \bm A_\alpha \bm y^n = 0 .
\end{equation}
By applying (\ref{5.20}), the time transition to a new level is provided as follows:
\begin{equation}\label{5.21}
\bm y^{n+1} = \bm y^n + \bm q + \bm g,
\end{equation}
\begin{equation}\label{5.22}
\begin{split}
& (\bm I+\sigma \tau \bm A_1 ) \bm q + \tau \bm A_1 \bm y^n = 0, \\
& (\bm I+\sigma \tau \bm A_2 ) \bm g + \tau \bm A_2 \bm y^n  = 0.
\end{split}
\end{equation}

When decomposing by rows (\ref{5.15}) for $\bm q$ and $\bm g$ from (\ref{5.22}), we obtain
\[
\begin{split}
& (1+ \sigma \tau A_{11} ) q_1 + \tau A_{11} y_1^{n} + \tau A_{12} y_2^{n} = 0 , \quad q_2 = 0 , \\
& (1+ \sigma \tau A_{22} ) g_2 + \tau A_{21} y_1^{n} + \tau A_{22} y_2^{n} = 0 , \quad g_1 = 0 .
\end{split}
\]
Therefore, from (\ref{5.20}) we have $\bm y^{n+1} = \bm y^n + \{q_1, g_2\}$.
Given this, we have (\ref{5.8}) for the approximate solution at the new time level.
Thus, the regularized decomposition-composition scheme with row splitting is algebraically equivalent to the decomposition-composition scheme with separation of the diagonal part of the problem operator (\ref{5.8}).
Therefore, we will not give computational data for the scheme (\ref{5.14}), (\ref{5.15}), (\ref{5.20}).

\begin{figure}[htbp]
\centering
\includegraphics[width=0.4\linewidth]{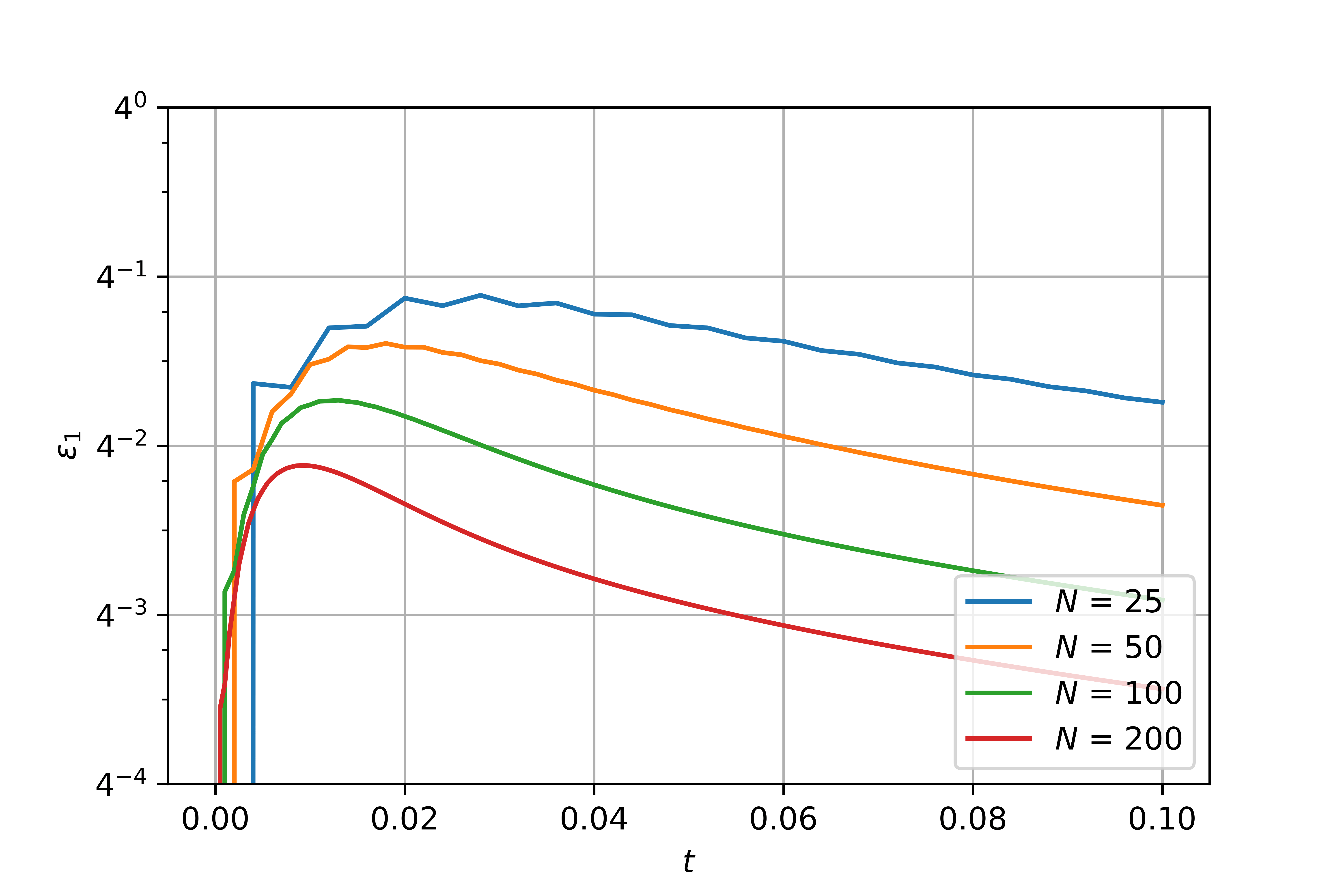} \includegraphics[width=0.4\linewidth]{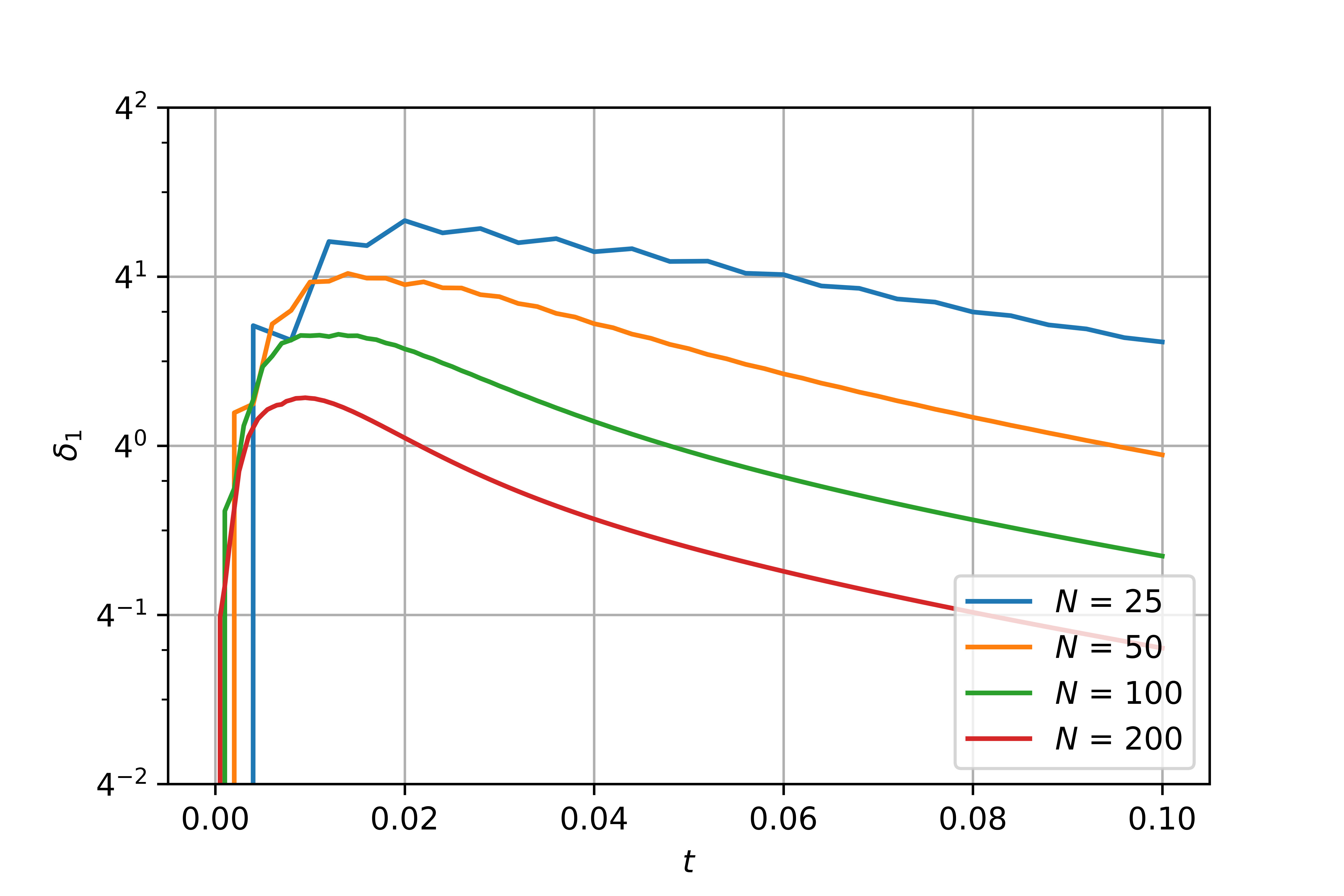} \\
\includegraphics[width=0.4\linewidth]{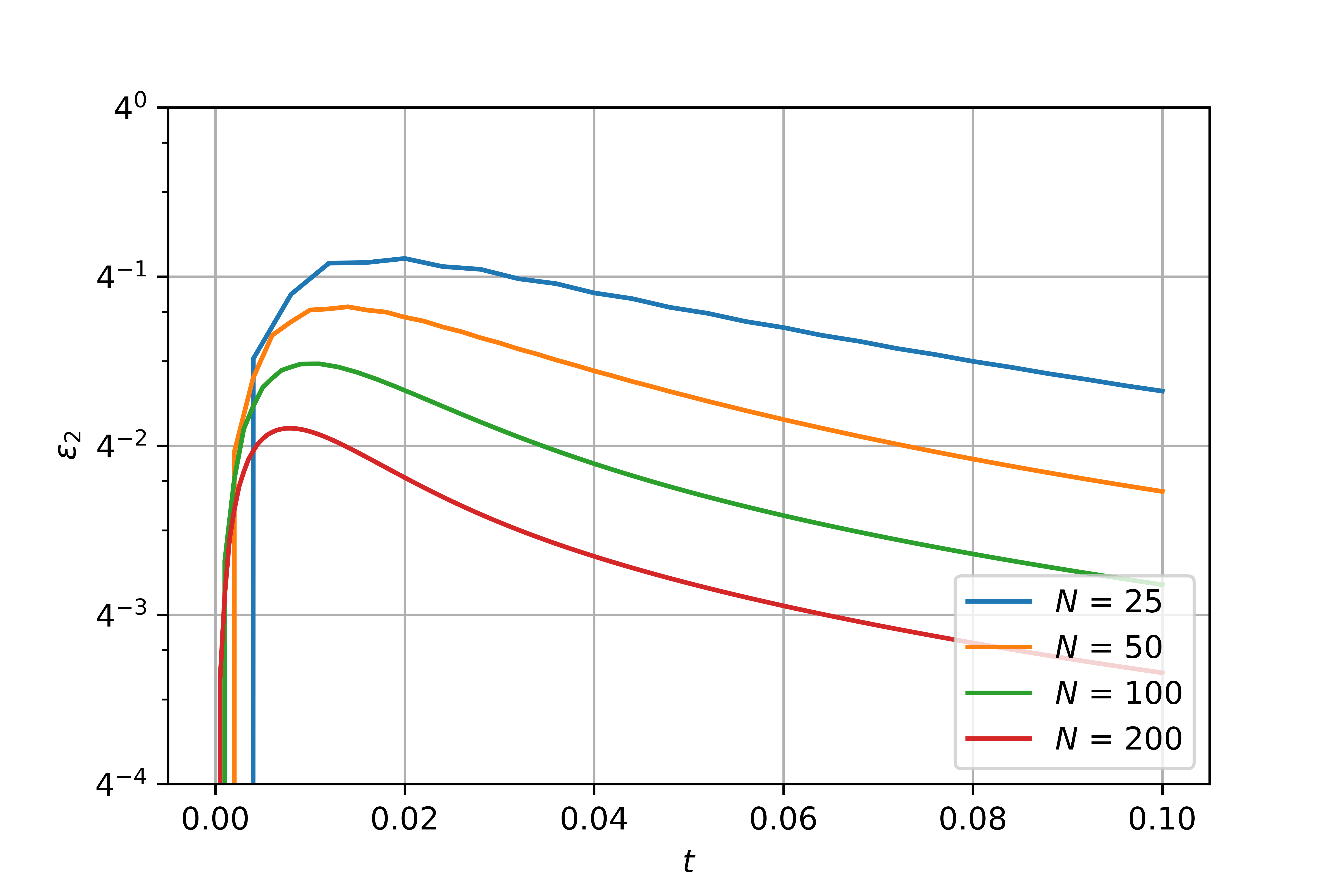} \includegraphics[width=0.4\linewidth]{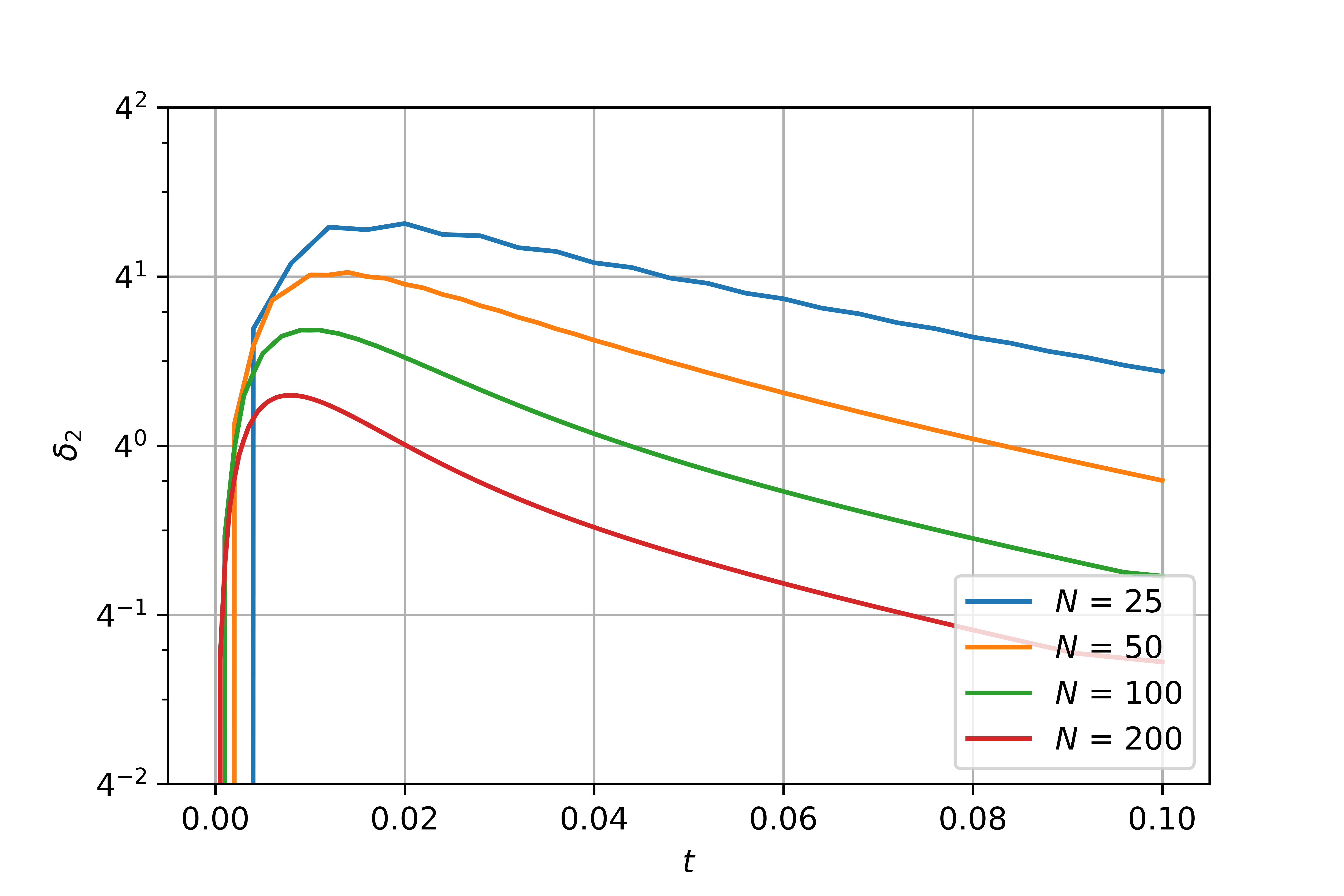}
\caption{Decomposition-composition scheme (\ref{5.23}).}
\label{f-10}
\end{figure}

When decomposing by columns (\ref{5.16}) for $\bm q$ and $\bm g$ we have
\begin{equation}\label{5.23}
\begin{split}
& (1+ \sigma \tau A_{11} ) q_1 + \tau A_{11} y_1^{n} = 0,
\quad q_2 + \tau A_{21} (y_1^{n} + \sigma q_1) = 0 , \\
& (1+ \sigma \tau A_{22} ) g_2 + \tau A_{22} y_2^{n} = 0 ,
\quad g_1 + \sigma \tau A_{12} g_2 + \tau A_{12} y_2^{n} = 0 .
\end{split}
\end{equation}
The results on the accuracy of the scheme (\ref{5.23}) at $\sigma = 1$ are shown in Fig.~\ref{f-10}.
For the test problem, the regularized scheme for column decomposition (Fig.~\ref{f-10}) is inferior to the regularized scheme for row decomposition (Fig.~\ref{f-2}) and comparable to the component-wise splitting scheme (Fig.~\ref{f-9}).

Computational experiments with second-order precision component-wise splitting schemes deserve special attention.
Instead of (\ref{5.14}), a four-component additive representation is used
\begin{equation}\label{5.24}
\bm A = \widetilde{\bm A}_1 + \widetilde{\bm A}_2 + \widetilde{\bm A}_3 + \widetilde{\bm A}_4 ,
\end{equation}
in which
\[
\widetilde{\bm A}_1 = \widetilde{\bm A}_4 = \frac{1}{2} \bm A_1,
\quad \widetilde{\bm A}_2 = \widetilde{\bm A}_3 = \frac{1}{2} {\bm A}_2 .
\]
Using the notation (\ref{5.24}), the component-wise splitting scheme at $\sigma = 1/2$ is written in the following form
\begin{equation}\label{5.25}
\begin{split}
& \left (\bm I +  \frac{1}{4} \tau \bm A_1 \right ) \frac{\bm y^{n+1/4} -  \bm y^{n}}{\tau} + \frac{1}{2} \bm A_1  \bm y^{n} = 0 , \\
& \left (\bm I +  \frac{1}{4} \tau \bm A_2 \right ) \frac{\bm y^{n+1/2} -  \bm y^{n+1/4}}{\tau} + \frac{1}{2} \bm A_2  \bm y^{n+1/4} = 0 , \\
& \left (\bm I +  \frac{1}{4} \tau \bm A_2 \right ) \frac{\bm y^{n+3/4} -  \bm y^{n+1/2}}{\tau} + \frac{1}{2} \bm A_2  \bm y^{n+1/2} = 0 , \\
& \left (\bm I +  \frac{1}{4} \tau \bm A_1 \right ) \frac{\bm y^{n+1} -  \bm y^{n+3/4}}{\tau} + \frac{1}{2} \bm A_1  \bm y^{n+3/4} = 0 .
\end{split}
\end{equation}
Similarly (\ref{5.17}), (\ref{5.18}) from (\ref{5.25}) for individual components we obtain the scheme (\ref{5.11}).
The fact of such equivalence concerning standard alternating direction schemes and locally one-dimensional schemes is well known (see, e.g., \cite{Samarskii1989}).
Therefore, for the two-equation test problem we are considering, it makes no sense to give computational data for (\ref{5.25}).

\section{Conclusions}\label{s-6}

\begin{enumerate}[(1)]

\item The problem of constructing efficient computational algorithms for approximate solution of the Cauchy problem for the system of linear coupled evolution equations of the first order is considered. Decoupling technology is based on additive decomposition (analysis) of the operator matrix of the system of equations and composition (synthesis) of the approximate solution due to particular time approximations. Solving auxiliary more straightforward problems for separate operator summands provides the transition to a new level in time.

\item Decoupling schemes of the first and second order of accuracy are constructed for two-level decomposition of the operator matrix in variants with separation of diagonal and triangular parts of the operator matrix of the problem.
The unconditional stability of two- and three-level decomposition-composition schemes in the corresponding finite-dimensional Hilbert spaces is proved.

\item Multicomponent decomposition based on the additive representation of the unit operator when solving systems of evolutionary equations gives additive splitting of the problem's operator matrix by rows or columns.
The composition of the solution is provided by applying component-wise splitting schemes and regularized additive operator-difference schemes.

\item The performance of the proposed decoupling schemes is demonstrated on a two-dimensional test problem.
We solve the problem for a system of two self- and cross-diffusion equations with discontinuous coefficients using standard piecewise linear finite element approximations over space.

\item We expect the decoupling technique to be successfully applied to various applied problems.
Among the first-priority problems of theoretical analysis are the construction and investigation of decomposition-composition schemes for systems of second-order evolution equations and the generalization of the results to problems with non-self-adjoint operator matrices.

\end{enumerate}

%

\begin{thebibliography}{10}
\expandafter\ifx\csname url\endcsname\relax
  \def\url#1{\texttt{#1}}\fi
\expandafter\ifx\csname urlprefix\endcsname\relax\def\urlprefix{URL }\fi
\expandafter\ifx\csname href\endcsname\relax
  \def\href#1#2{#2} \def\path#1{#1}\fi

\bibitem{KnabnerAngermann2003}
P.~Knabner, L.~Angermann, Numerical Methods for Elliptic and Parabolic Partial
  Differential Equations, Springer Verlag, 2003.

\bibitem{QuarteroniValli}
A.~Quarteroni, A.~Valli, Numerical Approximation of Partial Differential
  Equations, Springer, 2008.

\bibitem{Samarskii1989}
A.~A. Samarskii, The Theory of Difference Schemes, Marcel Dekker, New York,
  2001.

\bibitem{SamarskiiMatusVabischevich2002}
A.~A. Samarskii, P.~P. Matus, P.~N. Vabishchevich, Difference Schemes with
  Operator Factors, Kluwer Academic Pub, 2002.

\bibitem{Ascher2008}
U.~M. Ascher, Numerical Methods for Evolutionary Differential Equations,
  Society for Industrial and Applied Mathematics, 2008.

\bibitem{LeVeque2007}
R.~J. LeVeque, Finite Difference Methods for Ordinary and Partial Differential
  Equations. Steady-State and Time-Dependent Problems, Society for Industrial
  and Applied Mathematics, 2007.

\bibitem{DC}
P.~N. Vabishchevich, Computational decomposition and composition technique for
  approximate solution of nonstationary problems, Journal of Computational and
  Applied Mathematics 451~(116111) (2024) 1--18.

\bibitem{DecSys}
P.~N. Vabishchevich, Decoupling methods for systems of parabolic equations, in:
  I.~Lirkov, S.~Margenov (Eds.), Large-Scale Scientific Computing. LSSC 2021,
  Vol. 13127 of Lecture Notes in Computer Science, Springer, 2022, pp.
  287--294.

\bibitem{Ascher1995}
U.~M. Ascher, S.~J. Ruuth, B.~T.~R. Wetton, Implicit-explicit methods for
  time-dependent partial differential equations, SIAM Journal on Numerical
  Analysis 32~(3) (1995) 797--823.

\bibitem{HundsdorferVerwer2003}
W.~H. Hundsdorfer, J.~G. Verwer, Numerical Solution of Time-Dependent
  Advection-Diffusion-Reaction Equations, Springer Verlag, 2003.

\bibitem{Vabishchevich2020}
P.~N. Vabishchevich, Explicit-implicit schemes for first-order evolution
  equations, Differential Equations 56~(7) (2020) 882--889.

\bibitem{Marchuk1990}
G.~I. Marchuk, Splitting and alternating direction methods, in: P.~G. Ciarlet,
  J.-L. Lions (Eds.), Handbook of Numerical Analysis, Vol. I, North-Holland,
  1990, pp. 197--462.

\bibitem{VabishchevichAdditive}
P.~N. Vabishchevich, Additive Operator-Difference Schemes: Splitting Schemes,
  Walter de Gruyter GmbH, Berlin, Boston, 2013.

\bibitem{okubo2001diffusion}
A.~Okubo, S.~A. Levin, Diffusion and Ecological problems: Modern Perspectives,
  Springer, 2001.

\bibitem{murray2003mathematical}
J.~D. Murray, Mathematical Biology: II: Spatial Models and Biomedical
  Applications, Springer, 2003.

\bibitem{Saad2003}
Y.~Saad, Iterative Methods for Sparse Linear Systems, SIAM, 2003.

\bibitem{vabishchevich2014PTM}
P.~N. Vabishchevich, Three level schemes of the alternating triangular method,
  Computational Mathematics and Mathematical Physics 54~(6) (2014) 953--962.

\bibitem{vabishchevich2023nonunif}
P.~N. Vabishchevich, Operator-difference schemes on non-uniform grids for
  second-order evolutionary equations, Russ. J. Numer. Anal. Math. Modelling
  38~(4) (2023) 267--277.

\bibitem{vabishchevich2023subdomain}
P.~N. Vabishchevich, Subdomain solution decomposition method for nonstationary
  problems, Journal of Computational Physics 472 (2023) 111679.

\bibitem{SamVabGulin}
A.~A. Samarskii, P.~N. Vabishchevich, A.~V. Gulin, Stability of
  operator-difference schemes, Differ. Uravn. 35~(2) (1999) 152--187, in
  Russian.

\bibitem{SamarskiiVabischevich1998}
A.~A. Samarskii, P.~N. Vabishchevich, Regularized additive full approximation
  schemes, Dokl. Akad. Nauk 358 (1998) 461--464, in Russian.

\bibitem{Fryazinov1968}
I.~V. Fryazinov, Economical symmetrized schemes for solving boundary value
  problems for multi-dimensional parabolic equation, Zh. Vychisl. Mat. Mat.
  Fiz. 8 (1968) 436--443, in Russian.

\bibitem{Strang1968}
G.~Strang, On the construction and comparison of difference schemes, SIAM
  Journal on Numerical Analysis 5~(3) (1968) 506--517.

\bibitem{brenner2008mathematical}
S.~C. Brenner, L.~R. Scott, The Mathematical Theory of Finite Element Methods,
  Springer, 2008.

\bibitem{ciarlet2002finite}
P.~G. Ciarlet, The Finite Element Method for Elliptic Problems, SIAM, 2002.

\bibitem{DouglasRachford1956}
J.~J. Douglas, H.~H. Rachford, On the numerical solution of heat conduction
  problems in two and three space variables, Trans. Amer. Math. Soc. 82 (1956)
  421--439.

\bibitem{PeacemanRachford1955}
D.~W. Peaceman, H.~H. Rachford, The numerical solution of parabolic and
  elliptic differential equations, J. SIAM 3 (1955) 28--41.

\end{thebibliography}
%

\end{document}